\documentclass[a4paper,11pt]{article}
\pdfoutput=1
\usepackage{fancyhdr}
\usepackage{eurosym}
\usepackage[english]{babel}
\usepackage[utf8]{inputenc}
\usepackage{ mathrsfs }
\usepackage{amsmath}
\usepackage{amssymb}
\usepackage{amsthm}
\usepackage{bbm}
\usepackage{tikz}
\usetikzlibrary{patterns}
\usepackage{enumerate}
\usepackage{graphicx}
\graphicspath{ {./Figures/} }
\usepackage{enumerate}
\numberwithin{equation}{section}
\usepackage[backend = biber, sorting = nyt, defernumbers=true, style = numeric]{biblatex}
\usepackage{url}
\usepackage[margin = 1.5cm]{geometry}
\usepackage{csquotes}
\usepackage{hyperref}
\usepackage{todonotes}
\usepackage{pgfplots}
\usepackage{orcidlink}
\usepackage{tkz-euclide,subfigure}
\RequirePackage{doi}
\usepackage{microtype}

\usetikzlibrary{decorations.pathreplacing}
\usetikzlibrary{decorations.pathmorphing}

\DeclareMathOperator{\Var}{Var}
\DeclareMathOperator{\Cov}{Cov}

\newcommand{\R}{\mathbb{R}}
\newcommand{\N}{\mathbb{N}}
\newcommand{\Z}{\mathbb{Z}}

\newcommand{\Coh}{C_1}
\newcommand{\ClogW}{\mathcal{C}_{\log}}
\newcommand{\ClogWhat}{\widehat{\mathcal{C}}_{\log}}
\newcommand{\ExpPos}{\gamma}
\newcommand{\ExpPoshat}{\widehat{\gamma}}
\newcommand{\ConstWnW}{C_2}
\newcommand{\Cnull}{C_3}
\newcommand{\Cnullhat}{\widehat{C}_3}
\newcommand{\ExpStart}{C_4}
\newcommand{\ConstStart}{c_2}
\newcommand{\ExpStarthat}{\widehat{C}_4}
\newcommand{\ConstStarthat}{\widehat{c}_2}
\newcommand{\ConstWnWhat}{\widehat{C}_2}
\newcommand{\ptast}{p_t^\ast}
\newcommand{\ptasthat}{\widehat{p}_t^\ast}
\newcommand{\pnast}{p_n^\ast}
\newcommand{\pnasthat}{\widehat{p}_n^\ast}

\theoremstyle{plain}
\newtheorem{theorem}{Theorem}[section]
\newtheorem{definition}[theorem]{Definition}
\newtheorem{proposition}[theorem]{Proposition}
\newtheorem{corollary}[theorem]{Corollary}
\newtheorem{lemma}[theorem]{Lemma}
\theoremstyle{definition}

\newtheorem{remark}[theorem]{Remark}

\allowdisplaybreaks

\makeatletter
\let\blx@rerun@biber\relax
\makeatother

\hypersetup{
     pdfauthor={Xaver Kriechbaum},
    pdfborder={0 0 0},
   pdftitle = {Tightness BRWre}
}

\title{Tightness for branching random walk in time-inhomogeneous random environment}
\author{Xaver Kriechbaum\footnote{Weizmann Institute of Science, Israel. \url{xaver.kriechbaum@weizmann.ac.il}}\ \orcidlink{0009-0005-4470-299X}}
\date{}

\newcommand{\barrierprobUB}{p_n^{\smallfrown,x}(y)}
\newcommand{\barrierprobLB}{p_n^{\smallsmile,\xi_0}(y)}

\bibliography{LitBRWRE}

\begin{document}
\maketitle
\begin{abstract}
We consider a branching random walk in time-inhomogeneous random environment, in which all particles at generation $k$ branch into the same random number of particles $\mathcal{L}_{k+1}\ge 2$, where the $\mathcal{L}_k$, $k\in\N$, are i.i.d., and the increments are standard normal. Let $\mathbb{P}$ denote the law of $(\mathcal{L}_k)_{k\in\N}$, and let $M_n$ denote the position of the maximal particle in generation $n$. We prove that there are $m_n$, which are functions of only $(\mathcal{L}_k)_{k\in\{0,\dots, n\}}$, such that (with regard to $\mathbb{P}$) the sequence $(M_n-m_n)_{n\in\N}$ is tight with high probability.
\end{abstract}
\section{Introduction}
\subsection{Model definition and main result}
We study the maximum of branching random walk in (time-inhomogeneous) random environment (BRWre). Given a sequence $\mathcal{L} :=(\mathcal{L}_k)_{k\in\N}$ of natural numbers, which we refer to as \emph{environment}, we start with a single particle in position 0 at time (respectively generation) $0$. At each time $k\in\N_0$ all particles die after giving birth to $\mathcal{L}_{k+1}$ children, which take independent $\mathcal{N}(0,1)$ distributed jumps relative to their parents.  Given $\mathcal{L}$ let $\mathbf{T}$ denote the rooted tree in which all vertices with distance $k$ from the root are adjacent to exactly $\mathcal{L}_{k+1}$ vertices with distance $k+1$ from the root. For $u\in\mathbf{T}$ denote by $|u|$ the distance to the root and by $V(u)$ the position of the particle.  Denote $M_n := \max_{\{u : |u| = n\}} V(u)$. Given $\mathcal{L}$ we denote the law (respectively expectation) of our process by $\mathbb{P}_{\mathcal{L}}$ (respectively $\mathbb{E}_{\mathcal{L}}$). The law $\mathbb{P}_{\mathcal{L}}$ is called the \emph{quenched} law. We will consider $(\mathcal{L}_k)_{k\in\N}$ an i.i.d.\@ sequence of random variables on a probability space $(\Omega,\mathcal{F},\mathbb{P})$ such that $\mathcal{L}_1\ge 2$ $\mathbb{P}$-a.s.\@ and $\Var_{\mathbb{P}}[\log(\mathcal{L}_1)]<\infty$ and define  $P := \mathbb{P}\otimes \mathbb{P}_{\mathcal{L}}$ and $E = \mathbb{E}\otimes \mathbb{E}_{\mathcal{L}}$. We call $P$ the \emph{annealed} law of our process.

In the case of $(\mathcal{L}_k)_{k\in\N}$ constant the asymptotic behavior of $M_n$ is well understood, see for example Theorems 2 and 4 in  \cite[pp. 5, 9]{OZLN} for proofs that in this case $M_n$ grows with ballistic speed $v$ and that there is a constant $c$ such that $\mathbb{E}[M_n] = vn-c\log(n)+O(1)$. Furthermore, see Theorem 5.15 in \cite[p. 62]{ZhanShi} for a proof that $M_n-vn+c\log(n)$ converges in distribution for appropriate $c$, and a description of the limit.

The aim of our paper is to prove, that for the (random) centering $m_n$ defined in Definitions \ref{Def:Basics} and \ref{Def:pn}, the sequence $(M_n-m_n)_{n\in\N}$ is tight with respect to the annealed measure $P$.  We proceed by introducing some notation and defining the centering $m_n$.
\begin{definition}\label{Def:Basics}
For $k\in\N$, $\vartheta\in\R$, set
\begin{align}
\kappa_k(\vartheta) &:= \log(\mathcal{L}_k)+\frac{\vartheta^2}{2}, \notag\\
\kappa(\vartheta) &:= \mathbb{E}[\kappa_1(\vartheta)], \notag\\
\vartheta^\ast &:={\arg\inf} _{\vartheta>0}\frac{\kappa(\vartheta)}{\vartheta} = \sqrt{2\mathbb{E}[\log(\mathcal{L}_1)]}, \notag\\
K_n &:= \sum_{k=1}^n \kappa_k(\vartheta^\ast),\label{eq:K}\\
m_n &:= (\vartheta^\ast)^{-1}(K_n+\log(p_n)), \notag
\end{align}
where $p_n$ is the probability of a certain barrier event associated with the environmental random walk $(K_k)_{k\le n}$ which is defined in Lemma \ref{Lem:Boundpnastweak}. We note that $m_n$ is a function of the environment $\mathcal{L}$ since $p_n$ and $K_n$ depend on the environment. We note that $m_n$ is a random variable whose law is determined by $\mathbb{P}$.
\end{definition}

Next, we present our main theorem.
\begin{theorem}\label{Theo:Main}
We have that
$
(M_n-m_n)_{n\in\N}
$
is tight with regard to the annealed measure $P$.
\end{theorem}

\begin{remark}
\begin{enumerate}[(a)]
\item For the sake of ease in exposition, our model uses a concrete branching model and Gaussian increments. Our methods probably can be extended to cover non-Gaussian increments as well as the case that $(\mathcal{L}_k)_{k\in\N}$ is an i.i.d.\@ sequence with values in the set of probability measures on $\N$. However, this extension poses various technical issues, which would add a significant burden to an already long paper. Thus, we chose to consider the simplest model, for which new techniques needed to be developed. In Section \ref{sec:difficulties} we quickly sketch which difficulties one faces, when tackling the more general model.
\item It is natural to speculate that $M_n-m_n$ converges in distribution (under either the annealed or quenched laws). While we expect that improvements of our techniques might lead to such a statement, significant challenges remain, and thus this remains an open problem.
\end{enumerate}
\end{remark}

\subsection{Literature}
Our main motivation is \cite{MM_timeinh}. They consider a more general variant of the model we use, in which $(\mathcal{L}_k)_{k\in\N}$ is a sequence of point process laws instead of a sequence of natural numbers. At each time $k\in\N_0$ all particles die after generating children, whose placements (relative to the position of their parent) are distributed according to an i.i.d.\@ copy of a variable of law $\mathcal{L}_{k+1}$. They prove that under some integrability conditions on $\mathcal{L}_1$ and $\kappa_1$, which are analogues of the quantities in \eqref{eq:K} above, there is a constant $\varphi$ such that in $\mathbb{P}$-probability
\[
\lim\limits_{n\to \infty} \mathbb{P}_{\mathcal{L}}\left[M_n- \frac{1}{\vartheta^\ast} K_n \ge -\beta \log n \right] = \begin{cases}
1&\quad\text{if}\quad \beta>\varphi,\\
0&\quad\text{if}\quad \beta<\varphi,
\end{cases}
\]
in particular they prove that in $P$-probability
\[\lim\limits_{n\to \infty} \frac{M_n-\frac{1}{\vartheta^\ast} K_n}{\log(n)} = -\varphi.\] Their proof uses a ballot theorem for a random walk in random environment, see Theorem 3.3 in \cite{MM_timeinh}, which is proved in \cite{MM_BMaboveQRW}; this ballot theorem is used to evaluate the asymptotics of certain barrier probabilities, where the barrier itself is random. The main novel ingredient for the proof of Theorem \ref{Theo:Main} is a variant of this ballot theorem, see \eqref{eq:Structure1}, which lets us compare different starting heights relative to the chosen barrier. The model in \cite{MM_timeinh} goes back to \cite{BiggKyp}. In \cite{HuangLiu} it is proved that $M_n$ grows at a linear speed almost surely. Appendix A.2 in \cite{MM_timeinh} briefly discusses annealed tightness, i.e.\@ the content of Theorem \ref{Theo:Main}. They show that under the additional assumption that the increments of the BRW have bounded support one easily can show tightness of $(M_n-\mathbb{E}_{\mathcal{L}}[M_n])_{n\in\N}$ using the Dekking-Host argument from \cite{DH91}. The canonical way of then achieving tightness even for unbounded support, is to show that $\mathbb{E}[|M_n-\mathbb{E}_{\mathcal{L}}[M_n]|]  = O_n(1)$, however they can only show that it is of order $o_n(\log(n))$.

In \cite{MalleinInterfaces} a variant of the model in \cite{MM_timeinh} is considered. The environment is defined by using a P-uple $(\mathcal{L}_p)_{p\in \{1,\dots,P\}}$ of point process laws and, for fixed $0 = \alpha_0 < \alpha_1<\cdots <\alpha_P =1$, using $\mathcal{L}_p$ for reproduction at times $ n \alpha_{p-1} < k <  n \alpha_{p}$. They prove that under certain assumptions there are real numbers $v$ and $\lambda$ such that $(M_n-nv+\lambda \log(n))_{n\in\N}$ is tight.

In \cite{CernyDrewitz} a space-inhomogeneous branching random walk in random environment is considered. They take $(\xi(x))_{x\in\Z}$ an i.i.d.\@ family of random variables with $0<\mathrm{ess}\,\mathrm{inf}\xi(0)<\mathrm{ess}\,\mathrm{sup}\xi(0)<\infty$ and consider an initial configuration $u_0:\Z\to\N_0$. Given $\xi$ and $u_0$ at each $x\in\Z$ they place $u_0(x)$ particles at time 0, which all move independently according to continuous time simple random walk with jump rate 1. Furthermore, while at site $x$, a particle splits into two with rate $\xi(x)$ independently of everything else. They prove that under the right assumptions there is a velocity $v_0$ and a constant $\overline{\sigma}_{v_0}\in (0,\infty)$ such that the sequence of processes $\big(\frac{M(nt)-v_0nt}{\overline{\sigma}_{v_0}\sqrt{n}}\big)_{t\ge0}$ converges in ``annealed'' distribution to a standard Brownian motion. We proved in \cite{SubsTight} that there are deterministic subsequences $(t_k)_{k\in\N}$ and a function $m_{t_k}^\xi$ of the environment such that $(M_{t_k}-m_{t_k}^\xi)_{k\in\N}$ is (annealed) tight. Another motivation for the study in this paper is that we hope that the techniques used to prove Theorem \ref{Theo:Main} can be adapted to get tightness in the space-inhomogeneous model without the need to consider subsequences. It should be noted that recently in \cite{CerDreOsTight} $\mathbb{P}$-a.s.\@ tightness of $M_n$ around its quenched median has been shown for a space-continuous variant of the model, which is called branching Brownian motion in random environment. The proof uses a combination of analytic and probabilistic techniques. In particular, it is not clear to us if and how the proof could be adapted to cover the space-discrete version of the model from \cite{CernyDrewitz}.

In \cite{Fang} branching random walks allowing time dependence as well as local dependence between siblings are considered. They give assumptions under which $M_n$ is tight around the median of its distribution. The assumptions can be found in Sections 2 and 5 of \cite{Fang}. This result implies $\mathbb{P}$-a.s.\@ quenched tightness for BRWre, provided that the assumptions from \cite{Fang} hold $\mathbb{P}$-almost surely. While the quenched tightness in \cite{Fang} is stronger than the annealed tightness we prove in Theorem \ref{Theo:Main}, our result has the advantage, that we explicitly determine the correct centering $m_n$. Furthermore, the assumptions stated in \cite{Fang} do not cover the case that $\mathcal{L}_1$ has unbounded support. 

There are other studies of inhomogeneous branching random walks. Studies of  branching random walk in deterministic time-inhomogeneous environments can be found in \cite{BovHart}, \cite{FangZei}, \cite{MaillZei}, \cite{Mall},  \cite{doi:10.1080/03605302.2014.972744}, \cite{FOu}. In \cite{PerThLuZei} branching Brownian motion in an environment which is periodic in space has been analyzed, the corresponding PDE has been studied in \cite{Hamel2016}. Branching Brownian motion in deterministic space-inhomogeneous environment was studied in \cite{10.1214/aop/1176991677}, \cite{10.2307/2244200}.

\subsection{Structure of the paper}
In Section \ref{Subsec:MT1} we do some preliminary work, by defining barrier events and stating the many-to-one lemma in our setting. 

In Section \ref{Sec:Strctprf} we give a high-level description of the structure of the proof of Theorem \ref{Theo:Main}.

Section \ref{Sec:Sketch} is devoted to proving Theorem \ref{Theo:Main} while postponing some of the details to Sections \ref{Sec:UB} to \ref{Sec:Finish}. The proof of the theorem is broken down into three parts: an upper bound on the right tail of $M_n-m_n$ (proved in Sections \ref{Subsec:UB} and \ref{Sec:UB}), a lower bound on the right tail of $M_n-m_n$ (proved in Sections \ref{Subsec:LB} and \ref{Sec:DetLB}), and an upper bound on the left tail of $M_n-m_n$ (proved in Sections \ref{Subsec:LTail} and \ref{Sec:Finish}). Finally, in Section \ref{Sec:ProofMainTheo} we prove Theorem \ref{Theo:Main}, i.e.\@ the tightness of $M_n-m_n$, by combining the upper bounds on the right and left tails of $M_n-m_n$.

The rest of the paper contains technical estimates concerning barriers, which were used in Section \ref{Sec:UB}-\ref{Sec:Finish}.

In Section \ref{Sec:MoveStart} we compare the probability that a Brownian motion starting at $y\le x \le c$ stays below a barrier, with the probability that a Brownian motion starting at $x$ stays below the same barrier. 

In Section \ref{Sec:Moveh} we compare the probability to stay below a barrier shifted by a deterministic curve to the probability to stay below the same barrier without the added curve.

Section \ref{Sec:Crude} gives a very rough lower bound for a type of barrier event we frequently use.

These barrier estimates are proved in a rather general setup for the barrier, where certain time and barrier-dependent constants are assumed to be finite. We will need these families of constants to be tight for the environment we use. That this is the case is checked in Section \ref{Sec:Constants} using standard calculations. 

Finally, Section \ref{Sec:BarrComps} combines the results of Sections \ref{Sec:MoveStart} to \ref{Sec:Constants} into the exact barrier estimates we use in the proof of Theorem \ref{Theo:Main}.
\subsection{Notation}
We use $c_\varepsilon$, $C_\varepsilon$ to denote positive constants depending on $\varepsilon$, which will change from line to line and $c$, $C$ to denote positive constants, which will also change from line to line. Named constants will not change from line to line. We use $\xi_0\in \N_{\le0}$ to denote a fixed negative constant, which is smaller than $-e-1$.

Recall that $\mathbf{T}$ denotes the genealogical tree of our BRWre. For $u\in \mathbf{T}$ we use $(u_0,u_1,\dots, u_{|u|})$ to denote the vector of ancestors of the particle. For $u,v\in\mathbf{T}$ we denote by $u\wedge v$ the  last common ancestor of $u$ and $v$.

\noindent
{\bf Acknowledgements}  This project has received funding from the European Research Council (ERC) under the European Union's Horizon 2020 research and innovation programme (grant agreement No. 692452).

Thanks to Ofer Zeitouni for suggesting the problem and for many useful discussions. We also thank an anonymous referee for many helpful comments improving the presentation of this paper.

\section{Preliminaries: Barrier Events and a Many-to-One Lemma} \label{Subsec:MT1}

We first give a rather general definition of barrier events and probabilities.
\begin{definition}[Barrier Events]\label{Def:Barr}
For $I \subseteq \R$, $t:= \max I$, $(Z_s)_{s\in I}$ a real valued process, $f: I \to \R$ a function, $y\in \R$ and $J\subseteq \R$ an interval, define
\begin{align*}
\mathcal{B}_{I,f}^{y,J}(Z_{\cdot}) := \{\forall_{s \in I}\, y+Z_s +f(s) \le 0, y+Z_{t}+f(t) \in J\}.
\end{align*}
If $f=0$ or $y=0$ we suppress them from the notation. Also set $J_x := [x-1,x]$ for $x \in\R$. If $J =\R$ we suppress it from notation. We call events of this form barrier events.
\end{definition}

\begin{remark}
We emphasize that $\mathcal{B}_{I,f}^{y,J}$ is the probability that $y+Z_\cdot$ stays \emph{below} the barrier given by $-f$ and ends up in the interval $J-f(t)$. Often estimates of barrier events like the ballot theorem are stated in terms of $Z_\cdot$ staying above a barrier. However, this change of sign is common in the study of the maximum of BRW (see $A_v$ in \cite[p. 11]{OZLN}, $X_n(\delta)$ in \cite[p. 15]{MM_timeinh} and Lemma 2.4 in \cite[p. 9]{FOu}).
\end{remark}

\begin{definition}\label{Def:pn}
For $n\in\N$, let $(X_n)_{n\in\N}$ be i.i.d.\@ with $X_1\sim \mathcal{N}(\vartheta^\ast,1)$ and set
\begin{align*}
S_n &:= \sum_{k=1}^n X_k\sim \mathcal{N}(n\vartheta^\ast, n),\quad B_n := S_n-n\vartheta^\ast,\quad
W_n := K_n/\vartheta^\ast-n\vartheta^\ast.
\end{align*}We can realize $(B_n)_{n\in\N}$ as the values at integer times of a Brownian motion $(B_s)_{s\ge0}$ independent of $\mathcal{L}$. Furthermore we define $W_s$ for all $s\ge 0$ by linear interpolation. For $s\ge 0$ set $T_s := \vartheta^\ast(B_s-W_s)$. 

Finally, denote by $\mathbb{P}_{\mathcal{L}}^k$ the measure using the environment $(\mathcal{L}_{j+k})_{j\in\N}$. We also introduce the notation $S_n^{(k)} := \sum_{j=k+1}^{k+n} X_j$ and $T_n^{(k)} = \vartheta^\ast S_n^{(k)}-\sum_{j=k+1}^{k+n} \kappa_j(\vartheta^\ast)$.
\end{definition}
\begin{remark}\label{Rem:LiL}
Since $(\mathcal{L}_k)_{k\in\N}$ is i.i.d., so is $(W_{k+1}-W_k)_{k\in\N}$. Furthermore, by choice of $\vartheta^\ast$ we have $\mathbb{E}[W_2-W_1] = 0$ and finally, since we assume that $\Var[\log(\mathcal{L}_1)]<\infty$, we have $\Var[W_2-W_1]<\infty$. Thus we can apply the law of iterated logarithms for $(W_k)_{k\in\N}$, a version of which can be found in Theorem 5.17 in \cite{YuvPerBM}.

We also note that under the annealed law $T_n$ is the sum of i.i.d.\@ random variables with mean zero and finite variance. Thus $T$ behaves similarly to a Brownian motion.
\end{remark}

Now we have defined all random variables necessary to state the many-to-one lemma for our model.

\begin{lemma}[{Lemma 2.1 equation (2.4) and equation (2.5) in \cite{MM_timeinh}}]\label{Lem:MT1}
For any $n\in\N$ and $f:\R^n\to \R_{\ge0}$ non-negative, we have $\mathbb{P}$-a.s.
\begin{align*}
\mathbb{E}_{\mathcal{L}}\left[\sum_{|u| = n}f(\vartheta^\ast V(u_j)-K_j, j\le n)\right] &= \mathbb{E}_{\mathcal{L}}[e^{-T_n} f(T_1,\dots, T_n)].
\end{align*}

We also have, remembering the definition $S_n^{(k)} = X_{k+1}+\dots+X_{n+k}$, the time-shifted version of this
\[
\mathbb{E}_{\mathcal{L}}^k\left[\sum_{|u| = n}f(\vartheta^\ast V(u_j)-K_j, j\le n)\right] = \mathbb{E}_{\mathcal{L}}^k\left[e^{-\vartheta^\ast S_n^{(k)}+\sum_{j=k+1}^{k+n} \kappa_j(\vartheta^\ast)} f(T_1^{(k)},\dots, T_n^{(k)})\right].
\]
\end{lemma}
\begin{proof}
We only need to show that our definition of $S_n$ has the same distribution as the $S_n$ from \cite{MM_timeinh}.  They take
\[
\widetilde{S}_n := \sum_{k=1}^n \widetilde{X}_k,
\]
with $(\widetilde{X}_k)_{k\in\N}$ an independent sequence of random variables such that
\begin{equation}
\mathbb{P}_{\mathcal{L}}\left[\widetilde{X}_k\le x \right] = \mathbb{E}_{\mathcal{L}}\left[\sum_{\ell \in L_k}\mathbf{1}_{\ell \le x} e^{\vartheta^\ast \ell-\kappa_k(\vartheta^\ast)}\right]. \label{eq:Sngleich}
\end{equation}
Here, due to the simplified model we consider, $L_k$ is a set of $\mathcal{L}_k$ independent random variables with law $\mathcal{N}(0,1)$. Thus taking Definition 1.1 into account and letting $Z\sim\mathcal{N}(0,1)$, \eqref{eq:Sngleich} simplifies to
\begin{align*}
\mathbb{P}_{\mathcal{L}}\left[\widetilde{X}_k\le x \right] &= \mathcal{L}_k \mathbb{E}_{\mathcal{L}}\left[ \mathbf{1}_{ Z\le x} e^{\vartheta^\ast Z-\log(\mathcal{L}_k)-\frac{(\vartheta^\ast)^2}{2}}\right]= \mathbb{E}_{\mathcal{L}}\left[\mathbf{1}_{Z\le x} e^{\vartheta^\ast Z-\frac{(\vartheta^\ast)^2}{2}}\right]\\
&= \int_{-\infty}^x \frac{1}{\sqrt{2\pi}} e^{-\frac{z^2}{2}} e^{\vartheta^\ast z-\frac{(\vartheta^\ast)^2}{2}}\;\mathrm{d}z = \int_{-\infty}^x \frac{1}{\sqrt{2\pi}} e^{-\frac{(z-\vartheta^\ast)^2}{2}}\;\mathrm{d}z\\
&=  \mathcal{N}(\vartheta^\ast,1)((-\infty,x]).
\end{align*}
This shows that $\widetilde{X}_k \stackrel{d}{=} X_k$ from Definition \ref{Def:pn}, which implies $S_n\stackrel{d}{=}\widetilde{S}_n$.
\end{proof}

We recall that $\xi_0$ is a fixed negative integer.
\begin{lemma}\label{Lem:Boundpnastweak}
For $t\ge 0$ define
\begin{align}
p_t &:= \mathbb{P}_{\mathcal{L}}\left[\mathcal{B}_{[0,t]}^{\xi_0, J_{\xi_0}}( T_{\cdot}/\vartheta^\ast)\right]. \label{Def:pt}
\end{align}
Then there is a $C>0$ such that  $\mathbb{P}$-a.s.
\[
\limsup_{n\to\infty} \frac{|\log(p_{n})|}{\log(n)} \le C. 
\]
\end{lemma}

\begin{remark}
Because of Theorem 3.3 in \cite{MM_timeinh} we expect that in $\mathbb{P}$-probability as $n\to \infty$ we have 
\begin{equation*}
\log(p_n)/\log(n)\to -\lambda, 
\end{equation*} with $\lambda$ as in equation (3.3) of \cite{MM_timeinh}. To see this we can take $y =(\vartheta^\ast)^{-1}$, $x = -\xi_0/\vartheta^\ast$, $k = k' = 0$ in equation (3.4) of  \cite{MM_timeinh} to get that with high probability (as $n\to\infty$) we have $p_n \le n^{-\lambda+\varepsilon}$, i.e.\@ $\log(p_n)/\log(n) \le -\lambda +\varepsilon$. We can't extract an analogous lower bound from equation (3.4) in \cite{MM_timeinh}, since the infimum only is taken over $[a_n,x_n]$ with $\liminf_{n\to \infty} a_n/\log(n)>0$. Lemma \ref{Lem:Boundpnastweak} provides a lower bound, we note however that the lower bound we achieve is far from exact, i.e.\@ is not close to $-\lambda$.
\end{remark}

The following lemma allows us to factorize the probability of barrier events of Markov processes into two factors. This will be useful since we mostly do barrier computations for time intervals $[0,t]$ but want to apply them also for time intervals $[0,t_1]$ and $[t_1,t]$, which equations \eqref{eq:BarrSplit1} and \eqref{eq:BarrSplit2} enable. The proof of Lemma \ref{Lem:BarrSplit} is a straightforward application of the Markov property, which we omit.

\begin{lemma}\label{Lem:BarrSplit}
Let $I\subseteq\R$, $t:= \max I$, $(Z_s)_{s\in I}$ be a real valued  and translation invariant Markov process on some measure space $(\widetilde{\Omega},\widetilde{\mathcal{F}},\widetilde{\mathbb{P}})$, and $z_1,z_2\in \R$. Take $t_0\in I$ and set $I_1 := \{s\in I : s\le t_0\}$, $I_2 := \{s\in I : s\ge t_0\}$. We have that for $x_0\in \R$,
\begin{align}
\widetilde{\mathbb{P}}\left[\mathcal{B}_{I}^{z_1,J_{z_2}}(Z_{\cdot})\right] &\ge \widetilde{\mathbb{P}}\left[\mathcal{B}_{I_1}^{z_1+1,J_{x_0}}(Z_{\cdot})\right] \inf_{x\in J_{{x_0}-1}}\widetilde{\mathbb{P}}\left[\mathcal{B}_{I_2}^{x,J_{z_2}}(Z_{\cdot}-Z_{t_0})\right], \label{eq:BarrSplit1}\\
\widetilde{\mathbb{P}}\left[\mathcal{B}_{I}^{z_1,[z_2-2,z_2]}(Z_{\cdot})\right] &\ge \widetilde{\mathbb{P}}\left[\mathcal{B}_{I_1}^{z_1,J_{x_0}}(Z_{\cdot})\right] \widetilde{\mathbb{P}}\left[\mathcal{B}_{I_2}^{x_0,J_{z_2}}(Z_{\cdot}-Z_{t_0})\right].\label{eq:BarrSplit2}
\end{align}
\end{lemma}

\section{Structure of the proof} \label{Sec:Strctprf}
The first step in the proof of Theorem \ref{Theo:Main} is an upper bound on $\mathbb{P}_{\mathcal{L}}\left[M_n-m_n\ge -y \right]$ for $y<0$. We can get a good enough upper bound for $\mathbb{P}_{\mathcal{L}}\left[M_n-m_n \ge \xi_0+\log(n)^2 \right]$ by applying Lemma \ref{Lem:MT1}.  Thus we need to consider $y\in [-\log(n)^2,\xi_0]$. This is done using a first moment approach. For this, we introduce the barrier $\left((\vartheta^\ast)^{-1} K_k+\frac{k}{n\vartheta^\ast}\log(p_n)-y\right)_{k\in \{0,\dots,n\}}$ and take $\tau(y)$ to be the first time at which any particle in our system breaks this barrier. If no particle breaks the barrier until time $n$, we define $\tau(y)$ to be equal to $n$. It is immediate that 
\begin{align}
&\mathbb{P}_{\mathcal{L}}\left[M_n-m_n\in [\xi_0-y-1,\xi_0-y] \right] \notag\\
&\hskip2.15cm\le \mathbb{P}_{\mathcal{L}}\left[\tau(y)<n \right]+\mathbb{P}_{\mathcal{L}}\left[\tau(y) = n, M_n-m_n\in [\xi_0-y-1,\xi_0-y] \right]\notag\\
&\hskip2.15cm\le  \mathbb{P}_{\mathcal{L}}\left[\tau(y)<n \right]+\mathbb{P}_{\mathcal{L}}\left[\bigcup_{u:|u| = n} \mathcal{B}_{\{0,\dots, n\}, -(\vartheta^\ast)^{-1}(\frac{\cdot}{n}\log(p_n)+K_\cdot)}^{y,J_{\xi_0}}(V(u_\cdot)) \right]. \label{eq:StrctBl}
\end{align}
The second inequality of \eqref{eq:StrctBl} uses that if $\tau(y) = n$ no particle breaches the barrier before the last step while if $M_n-m_n \in [\xi_0-y-1,\xi_0-y]$ there is a particle which ends up in $[m_n +\xi_0-y-1,m_n +\xi_0-y]$. Both combined mean that there is a particle which satisfies the barrier event $\mathcal{B}_{\{0,\dots, n\}, - (\vartheta^\ast)^{-1}(\frac{\cdot}{n}\log(p_n)+K_\cdot)}^{y,J_{\xi_0}}(V(u_\cdot))$.

Applying the Markov inequality and the many-to-one Lemma \ref{Lem:MT1} to the second summand in \eqref{eq:StrctBl} yields that
\begin{equation}
\mathbb{P}_{\mathcal{L}}\left[M_n-m_n\in J_{\xi_0-y}\right] \le \mathbb{P}_{\mathcal{L}}\left[\tau(y)<n\right]+Cp_n^{-1} e^{\vartheta^\ast y}\mathbb{P}_{\mathcal{L}}\left[\mathcal{B}_{\{0,\dots,n\},\frac{\cdot}{n\vartheta^\ast}\log(p_n)}^{y,J_{\xi_0}}(T_\cdot/\vartheta^\ast)\right]. \label{eq:Mnmn2parts}
\end{equation}
The heart of the matter is now to prove that for all $\varepsilon>0$, there are $c_\varepsilon, C_\varepsilon$ such that
\begin{equation}
\liminf_{n\to \infty}\mathbb{P}\left[\frac{\mathbb{P}_{\mathcal{L}}\left[\mathcal{B}_{\{0,\dots,n\},\frac{\cdot}{n\vartheta^\ast}\log(p_n)}^{y,J_{\xi_0}}(T_{\cdot}/\vartheta^\ast)\right]}{p_n} \le C_\varepsilon |y|^{c_\varepsilon}\right] \ge 1-\varepsilon.  \label{eq:Structure1}
\end{equation}
In Figure \ref{fig:UpBarr} we illustrate an event similar to the one in \eqref{eq:Structure1}. The main difference is that for \eqref{eq:Structure1} the thin black Brownian motion sample would need to stay below the dashed line instead of below the thick black line. We recall that $p_n = \mathbb{P}_{\mathcal{L}}\left[\mathcal{B}_{[0,n]}^{\xi_0, J_{\xi_0}}(T_{\cdot}/\vartheta^\ast)\right]$. Equation (3.4) in \cite{MM_timeinh} roughly corresponds to our equation \eqref{eq:Structure1}. In \cite[Eq. (3.4)]{MM_timeinh} they bound the growth in $n$ of $\mathbb{P}_{\mathcal{L}}\left[\mathcal{B}_{[0,n]}^{y,J_{\xi_0}}(T_{\cdot}/\vartheta^\ast)\right]$ from above (in $\mathbb{P}$-probability and uniformly in $y$ with $|y|$ small enough). In particular, they bound the growth of $p_n$.  In contrast equation \eqref{eq:Structure1} controls how much starting lower increases the probability to stay below the barrier.  We sketch below how equation \eqref{eq:Structure1} is proved. Handling the second summand in \eqref{eq:Mnmn2parts} is immediate once \eqref{eq:Structure1} is established, the first summand can be handled using \eqref{eq:Structure1} and standard techniques (these techniques are for example sketched in Remark 9 of \cite[p. 16]{OZLN}). In the actual proof we need a slightly more complicated barrier, which has an additional (deterministic) curve added, this does only marginally change the proof of \eqref{eq:Structure1}, however we need to then remove this curve after we have moved the starting point. This is mostly an application of Girsanov.

The next step in proving Theorem \ref{Theo:Main} is a lower bound on $\mathbb{P}_{\mathcal{L}}\left[M_n-m_n\ge \xi_0-y-1 \right]$ for $y\in [-\log(n)^2,\xi_0]$. This is done using a second moment approach, for this we define $Z_n(y)$ to be the number of particles  for which $(V(u_k)-(\vartheta^\ast)^{-1}K_k)_{k\in \{0,\dots, n\}}$ stays below $(\frac{k}{n\vartheta^\ast}\log(p_n)-y)_{k\in\{0,\dots,n\}}$ and ends up in $[(\vartheta^\ast)^{-1}\log(p_n)+\xi_0-y-1,(\vartheta^\ast)^{-1}\log(p_n)+\xi_0-y] = J_{(\vartheta^\ast)^{-1}\log(p_n)+\xi_0-y}$. Using Cauchy-Schwartz yields that
\[
\mathbb{P}_{\mathcal{L}}[M_n-m_n\ge \xi_0-y-1] \ge \mathbb{P}_{\mathcal{L}}[Z_n(y)\ge 1]\ge \frac{\mathbb{E}_{\mathcal{L}}[Z_n(y)]^2}{\mathbb{E}_{\mathcal{L}}[Z_n(y)^2]}
\]
and a lower bound on the right tail can be achieved by bounding $\mathbb{E}_{\mathcal{L}}[Z_n(y)]$ and $\mathbb{E}_{\mathcal{L}}[Z_n(y)^2]$. The first is done similar to the calculation in the upper bound on the right tail of $M_n-m_n$, using that for all $\varepsilon>0$ there are $c_\varepsilon, C_\varepsilon$ such that
\begin{equation}
\liminf_{n\to\infty} \mathbb{P}\left[ \frac{\mathbb{P}_{\mathcal{L}}\left[\mathcal{B}_{\{0,\dots,n\},\frac{\cdot}{\vartheta^\ast n}\log(p_n)}^{y,J_{\xi_0}}(T_{\cdot}/\vartheta^\ast)\right]}{p_n}\ge C_\varepsilon |y|^{-c_\varepsilon} \right] \ge1-\varepsilon \label{eq:Structure2}
\end{equation} instead of \eqref{eq:Structure1}. In Figure \ref{fig:DownBarr} we illustrate an event similar to the one in \eqref{eq:Structure2}. The main difference is that for \eqref{eq:Structure2} the thin black Brownian motion sample would need to stay below the dashed line instead of below the thick black line. For the bound on the second moment we use $\Lambda_k(y)$ to denote the number of pairs of particles in generation $n$, which are both counted in $Z_n(y)$, whose last common ancestor is in generation $k$. By definition
\[
Z_n(y)^2 = \sum_{k=0}^{n-1} \Lambda_k(y)+Z_n(y).
\]
Thus we need to prove an upper bound on $\mathbb{E}_{\mathcal{L}}[\Lambda_k(y)]$, we note that this proceeds along similar lines as \cite[Section 4.1]{MM_timeinh}. We quickly sketch how to get the upper bound. For this, we let $w := u\wedge v$ be the $k$-generational common ancestor of a pair $u$, $v$ counting for $\Lambda_k(y)$. Then we condition on the position $V(w)$ and use that given $V(w)$ and the environment $\mathcal{L}$, both $(V(u_j)-V(w))_{j\in \{k,\dots, n\}}$ and $(V(v_j)-V(w))_{j\in \{k,\dots,n\}}$ are independent and identically distributed to write
\begin{align*}
&\mathbb{P}_{\mathcal{L}}[ u,\ v\ \text{are counted by } \Lambda_k(y)] \\
= &\int_{-\infty}^0 \mathbb{P}_{\mathcal{L}}[\mathcal{B}_{\{k,\dots, n\}, -\frac{\cdot-k}{\vartheta^\ast n}\log(p_n) }^{r,J_{m_n+\xi_0-y}}(V(u_\cdot)-V(w)-(K_\cdot-K_k)/\vartheta^\ast)]^2\\
&\qquad \cdot \mathbb{P}_{\mathcal{L}}[\mathcal{B}_{\{0,\dots, k\}, -\frac{\cdot}{n\vartheta^\ast}\log(p_n) }^{y,\R}(V(w_\cdot)- K_j/\vartheta^\ast)| V(w_k)-K_k/\vartheta^\ast  \in dr].
\end{align*}The two factors are then handled by using Lemma \ref{Lem:MT1} and \eqref{eq:Structure1}.   In the proof it will be helpful to use a (downwards) curved barrier instead of the straight barrier used here to force $V(w)$ to be distanced from $K_k/\vartheta^\ast$.

The upper bound on the left tail of $M_n-m_n$ follows from the lower bound on the right tail by using the tree structure and the fact that quenched the subtrees below two different vertices in generation $k$ are independent of each other (for the homogeneous situation this argument can be found on the first half of page 17 in \cite{OZLN}).

Having upper bounds on both the left and the right tail of $M_n-m_n$ then quickly allows to prove Theorem \ref{Theo:Main}, i.e.\@ tightness of $(M_n-m_n)_{n\in\N}$.

We finish this section by sketching how we prove \eqref{eq:Structure1}, in our opinion this is the most interesting part of our paper. We ignore the random part and the $\log$-term of the barrier in this sketch, the difficulties added by these are either already handled by \cite{MM_BMaboveQRW} or minor.\footnote{Of course without the random part and the $\log$-term in the barrier we could just use the reflection principle to get an exact result. Because of this, the strategy we use is not needed in the homogeneous setting.}We also switch to continuous time, which can be justified using standard arguments. Thus let $(B_s)_{s\ge0}$ be a Brownian motion with respect to $\mathbb{P}_{\mathcal{L}}$. As a first step we use a correlation inequality for Brownian bridges to prove that
\begin{equation}
\mathbb{P}_{\mathcal{L}}\left[\mathcal{B}_{[0,n]}^{\xi_0,J_{\xi_0}}(B_{\cdot})\right] \ge \mathbb{P}_{\mathcal{L}}\left[\mathcal{B}_{[0,y^2]}^{\xi_0,[Cy,cy]}(B_{\cdot}) \right]\mathbb{P}_{\mathcal{L}}\left[\mathcal{B}_{[y^2,n]}^{\xi_0,J_{\xi_0}}(B_{\cdot})\right]\label{eq:SketchMain0}
\end{equation}
From \cite{MM_BMaboveQRW} we know that there is a $\gamma>0$ such that $\mathbb{P}_{\mathcal{L}}\left[\mathcal{B}_{[0,y^2]}^{\xi_0,[Cy,cy]}(B_{\cdot})\right] \ge C|y|^{-2\gamma}$ with high probability, such that we approximately get
\begin{equation}
\mathbb{P}_{\mathcal{L}}\left[\mathcal{B}_{[0,n]}^{\xi_0,J_{\xi_0}}(B_{\cdot})\right] \ge C|y|^{-2\gamma}\mathbb{P}_{\mathcal{L}}\left[\mathcal{B}_{[y^2,n]}^{\xi_0,J_{\xi_0}}(B_{\cdot})\right]. \label{eq:SketchMain1}
\end{equation}
On the other hand by monotonicity we get
\begin{equation}
\mathbb{P}_{\mathcal{L}}\left[\mathcal{B}_{[0,n]}^{y,J_{\xi_0}}(B_{\cdot})\right] \le \mathbb{P}_{\mathcal{L}}\left[\mathcal{B}_{[y^2,n]}^{y,J_{\xi_0}}(B_{\cdot})\right]. \label{eq:SketchMain2}
\end{equation}
By using \eqref{eq:SketchMain0}, \eqref{eq:SketchMain1} and \eqref{eq:SketchMain2} we reduce the proof of \eqref{eq:Structure1}  to proving an upper bound on 
$
\mathbb{P}_{\mathcal{L}}\left[\mathcal{B}_{[y^2,n]}^{y,J_{\xi_0}}(B_{\cdot})\right]/\mathbb{P}_{\mathcal{L}}\left[\mathcal{B}_{[y^2,n]}^{\xi_0,J_{\xi_0}}(B_{\cdot})\right].
$
Using the Markov property for Brownian motion at time $y^2$ we can write for $x\in \{\xi_0,y\}$
\[
\mathbb{P}_{\mathcal{L}}\left[\mathcal{B}_{[y^2,n]}^{x,J_{\xi_0}}(B_{\cdot})\right] = \int_{-\infty}^0 g_{x,y^2}(z) \mathbb{P}_{\mathcal{L}}\left[\forall_{s\in [y^2,n]} z+(B_s-B_{y^2}) \le 0, z+(B_n-B_{y^2}) \in J_{\xi_0} \right] \;\mathrm{d}z,
\]
where we use $g_{\mu,\sigma^2}(z)$ to denote the density of a normal distribution with mean $\mu$ and variance $\sigma^2$. For $0\ge z \ge y\log(|y|)$ we know that $g_{y,y^2}(z)/g_{\xi_0,y^2}(z) \le Cy^c$ and we can compare the integral over this region for $x = \xi_0$ to the same integral for $x = y$. The region $z\ge -C\log(t)$ is negligible for $C$ big enough, since we know from \cite{MM_BMaboveQRW} that $\left[\mathcal{B}_{[0,n]}^{y,J_{\xi_0}}(B_{\cdot})\right]$ is polynomial in $t$. The region $z\in [-C\log(t),y\log(|y|)]$ requires extra care, we handle this with an inductive scheme in which we recursively drop the barrier for a time interval depending on $z$ and again split up the possible locations after this time interval into three regions. A more detailed overview of this inductive scheme is given in Section \ref{Sec:MainOS}.

\subsection{Remarks on convergence in distribution} \label{Sec:ConviD}
The natural next question is whether $M_n-m_n$ converges in distribution. We think that this is the case and can be proved by improving our techniques and adopting the approach from \cite{ConviD}. We set $\mathbb{P}_{\mathcal{L}}\left[\mathcal{B}_{\{0,\dots, n\},\frac{\cdot}{n\vartheta^\ast}\log(p_n)}^{-z,J_{\xi_0}}(T_{\cdot}/\vartheta^\ast) \right]/p_n =:f_n(z)$, $z\in \R$, $n\in\N$.  One in our opinion sufficient improvement of  \eqref{eq:Structure1}, is to show that there exists a $f:\R\to\R$ such that for all big enough $z\in \R$, $a_l\to \infty$
 \begin{align}
 \lim\limits_{n\to \infty} f_n(z) &= f(z), \label{eq:Condf1iD}\\
 \lim\limits_{l\to \infty} \frac{f(z+a_l)}{f(a_l)} &= g(z),\qquad \lim\limits_{z\to\infty} f(z) = \infty. 
\label{eq:Condf2iD}
  \end{align}   
If one were able to establish these equalities, it is likely that the procedure from \cite{ConviD} can also be applied to our situation, although we have not checked this carefully.
\subsection{Difficulty in extending to the model from \cite{MM_timeinh}} \label{sec:difficulties}
As mentioned above, in \cite{MM_timeinh} the environment $(\mathcal{L}_k)_{k\in\N}$ is a sequence of point process laws instead of a sequence of natural numbers. As a consequence we'd need to change the definition of $\kappa_n(\vartheta)$ and $v$ from Definition \ref{Def:Basics} as follows
\begin{equation*}
\kappa_n(\vartheta) =\log \mathbb{E}_{\mathcal{L}}\left[\sum_{\ell \in L_n} e^{\vartheta \ell}\right],\quad
v = \inf_{\vartheta>0} \frac{\mathbb{E}[\kappa_1(\vartheta)]}{\vartheta},
\end{equation*}
where $L_n\sim \mathcal{L}_n$. Furthermore, we'd need to implicitly define $\vartheta^\ast>0$ via $\vartheta^\ast \mathbb{E}[\kappa_1(\vartheta^\ast)]-\mathbb{E}[\kappa_1(\vartheta^\ast)] = 0$. This has already been done in \cite[p. 3]{MM_timeinh} and we refer the reader there for more details on these changes and the assumptions, which need to be made. Beyond this, the more general model introduces two further difficulties.
\begin{enumerate}[(i)]
\item We recall that for a particle $u$ in generation $n$ we use $u_j$ to denote the $j$-th generational ancestor of $u$ and $V(u)$ to denote the position of $u$. For $u^{(1)}$, $u^{(2)}$ two different descendants of a particle $u$ the two random variables $(V(u^{(1)}_j)-V(u))_{j\ge |u|}$, $V((u^{(2)}_j)-V(u))_{j\ge |u|}$ won't be independent in the more general model. The main place in which we use this independence is in the upper bound of $\mathbb{E}[\Lambda_k^2(y)]$ in Section \ref{Subsec:Zn2UBDet}, which is used for the second moment calculation in Section \ref{Subsec:Zn2UB}. There we consider pairs of particles $u^{(1)}$, $u^{(2)}$ in generation $n$ for which the last common ancestor $u := u^{(1)}\wedge u^{(2)}$ is in generation $k$. Then we condition on everything happening up to generation $k$, in particular $V(u) = V(u^{(1)}_k) = V(u^{(2)}_k)$ is measurable with regard to this conditioning. We use that under this conditioning for any function $f: \R^{n-k}\to \{0,1\}$ the events $\left\{f\left((V(u^{(1)}_{k+1}),\dots, V(u^{(1)}_n) \right) = 1\right\}$, $\left\{f\left((V(u^{(2)}_{k+1}),\dots, V(u^{(2)}_n) \right) = 1\right\}$ are independent.

To fix this lack of independence, we can instead condition on everything happening up to generation $k+1$, such that $(V(u^{(1)}_{k+1}), V(u^{(2)}_{k+1}))$ is measurable with regard to the conditioning. This change will slightly complicate the calculations in Section \ref{Subsec:Zn2UBDet}, but won't be a big deal. This was already relevant and dealt with for Section 4.2 in \cite{MM_timeinh}, see for example their equation (4.5). 
\item The main difficulty in extending our result is that in the general model $S_n$ and $B_n$ from Definition \ref{Def:pn} won't be Gaussian random variables. As in Section 2 in \cite{MM_timeinh} we'd need to define $S_n = \sum_{k=1}^n X_k$, $(X_k)_{k\in\N}$ independent with 
\[
\mathbb{P}_{\mathcal{L}}[X_k \le x] = \mathbb{E}_{\mathcal{L}}\left[\sum_{\ell\in L_n} \mathbf{1}_{\ell\le x} e^{\vartheta^\ast \ell-\kappa_n(\vartheta^\ast)}\right],
\]
which won't be Gaussian in general. In Sections \ref{Sec:MoveStart} and \ref{Sec:Moveh} we use that $(B_k)_{k\in\N}$ is a Brownian motion evaluated at integer times. Namely, we use a correlation inequality for Brownian bridges, Lemma \ref{Lem: Association}, the explicit form of the Gaussian density and the Girsanov theorem, in particular we use that there is a simple way of tilting normal random variables in a way only changing the mean. Similar issues appear when talking about BRW in homogeneous environment and have been handled there, see \cite{ConviD}. The way Section 6 in  \cite{MM_BMaboveQRW} tackled this difficulty, is by approximating $(B_k)_{k\in\N}$ (and the barrier $(W_k)_{k\in\N}$) by Brownian motion, using a variant of the KMT approximation for independent sums of random variables proved in \cite{Sakhanenko}. This approach is somewhat more delicate in our situation for the following reason. Applying KMT on all of $[0,n]$ at once introduces an error of size roughly $\log(n)^c$, which then changes the starting point relative to the barrier by $\log(n)^c$. This changes the probability by a factor, which is polynomial in $\log(n)$. Since \cite{MM_BMaboveQRW} cares about 
\[
\lim\limits_{n\to\infty} \frac{-\log\left(\mathbb{P}_{\mathcal{L}}\left[\forall_{k\le n} x+B_k\ge W_k, B_n-W_n \in [a\sqrt{n},b\sqrt{n}\right]\right)}{\log(n)}
\]
they can ignore such a multiplicative factor. However, for \eqref{eq:Structure1} we need to understand how changes of $O(1)$ in the starting position change the probability and thus an error of size $\log(n)^c$ in the starting point is far to rough for our results.  This can be circumvented by splitting $[0,n]$ into dyadic intervals $I_1$, \dots, $I_{\log_2(n)}$ and using KMT on each of these intervals individually, but doing so is somewhat technical.
\end{enumerate}
\section{Proof of Theorem \ref{Theo:Main}}\label{Sec:Sketch}
In this section, we prove Theorem \ref{Theo:Main} leaving quite some details and proofs to the later sections.

\subsection{Upper Bound on the Right Tail of $M_n-m_n$}\label{Subsec:UB}

In this section, we prove the following theorem, postponing certain auxiliary estimates to Sections \ref{Sec:UB} and \ref{Sec:BarrComps}.
\begin{theorem}\label{Theo:MainUB}
For all $\varepsilon>0$, there exists a $y_\varepsilon<0$ such that
\[
\liminf_{n\to \infty} \mathbb{P}\left[\mathbb{P}_{\mathcal{L}}[M_n-m_n\ge -y_{\varepsilon}]\le \varepsilon\right]\ge 1-\varepsilon.
\]
\end{theorem}

The proof of Theorem \ref{Theo:MainUB} uses two ingredients. Lemma \ref{Lem:Ausreißer} controls the probability that $M_n-m_n$ is large and Proposition \ref{Prop:EndingUp} controls small to medium deviations of $M_n-m_n$.
\begin{lemma}[{Lemma 2.3 in \cite{MM_timeinh}}]\label{Lem:Ausreißer}
For any $y<0$, we have $\mathbb{P}$-a.s.
\[
\mathbb{P}_{\mathcal{L}}\left[\exists u\in \mathbf{T} :  V(u) +y> (\vartheta^\ast)^{-1} K_{|u|}\right]\le e^{\vartheta^\ast y}.
\]
\end{lemma}
\begin{proposition}\label{Prop:EndingUp}
For all $\varepsilon>0$, there are $c_\varepsilon, C_\varepsilon>0$ such that 
\[
\liminf_{n\to \infty} \mathbb{P}\left[\bigcap_{y\in [-\log(n)^2,\xi_0]\cap \Z}\left\{\mathbb{P}_{\mathcal{L}}[M_n-m_n \in [\xi_0-y-1,\xi_0-y]]\le C_\varepsilon e^{\vartheta^\ast y}|y|^{c_\varepsilon}\right\}\right]\ge 1-\varepsilon.
\]
\end{proposition}

\begin{proof}[Proof of Theorem \ref{Theo:MainUB} assuming Proposition \ref{Prop:EndingUp}]
Take $y_\varepsilon \in\Z\cap(-\infty,\xi_0]$ such that for $c_\varepsilon, C_\varepsilon$ as in Proposition \ref{Prop:EndingUp} $C_{\varepsilon}\sum_{y\in (-\infty,y_\varepsilon]\cap \Z}^{\infty} |y|^{c_{\varepsilon}}e^{-\vartheta^\ast y}\le \varepsilon/2$. Then, by Proposition \ref{Prop:EndingUp}
\begin{align}
\liminf_{n\to \infty}\mathbb{P}\left[\mathbb{P}_{\mathcal{L}}\left[M_n-m_n\in [\xi_0-y_\varepsilon, \xi_0+\log(n)^2\right]\le \varepsilon/2\right]\ge 1-\varepsilon. \label{eq:TyingTogether1}
\end{align}
On the other hand, from Lemma \ref{Lem:Ausreißer}, $\mathbb{P}$-a.s.,
\begin{align}
\mathbb{P}_{\mathcal{L}}\left[M_n-m_n\ge \xi_0+\log(n)^2\right]\le e^{-\vartheta^\ast \log(n)^2-\vartheta^\ast \xi_0-\log(p_n)}. \label{eq:zwischen}
\end{align}
Combining \eqref{eq:zwischen} and Lemma \ref{Lem:Boundpnastweak} yields that $\mathbb{P}$-a.s.\@ $\mathbb{P}_{\mathcal{L}}\left[M_n-m_n\ge \xi_0+\log(n)^2\right]\to 0$. In particular, this holds in probability, which implies
\begin{align}
\lim_{n\to \infty} \mathbb{P}\left[\mathbb{P}_{\mathcal{L}}\left[M_n-m_n\ge \xi_0+\log(n)^2\right]\le \varepsilon/2\right]=1. \label{eq:TyingTogether2}
\end{align}
Combining \eqref{eq:TyingTogether1} and \eqref{eq:TyingTogether2} finishes the proof.
\end{proof}

\begin{figure}
\includegraphics{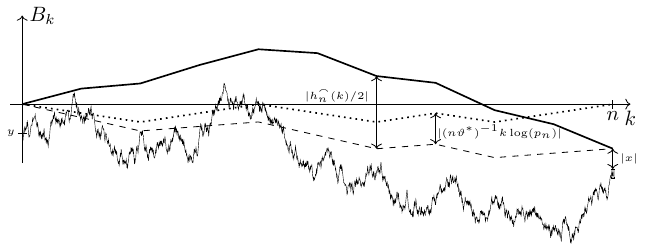}
\centering
\caption{The event in \eqref{Def:barrprobUB}. Drawn are $(K_k/\vartheta^\ast)_{k\le n}$ (dotted line), $(K_k/\vartheta^\ast+\frac{k}{n\vartheta^\ast} \log(p_n))_{k\le n}$ (dashed line), $(K_k/\vartheta^\ast+\frac{k}{n\vartheta^\ast}\log(p_n)-h_n^{\smallfrown}(k)/2)_{k\le n}$ (thick black line) and a sample path of $(y+B_s)_{s\le n}$ realizing the barrier event $\mathcal{B}_{\{0,\dots, n\},m_{n,h/2}^{\smallfrown}}^{y,J_x}(T_\cdot/\vartheta^\ast)$ (thin black line). We have the parameters $n = 10$, $\mathcal{L}_1\sim\mathrm{Unif}(\{2,3\})$, $(\mathcal{L}_k)_{k\le 10} = (2,2,3,3,2,2,3,2,3,3)$.}
\label{fig:UpBarr}
\end{figure}

Thus it remains to prove Proposition \ref{Prop:EndingUp}, which uses barrier computations. For this, we need the following definition, for which we recall \eqref{eq:K}, \eqref{Def:pt} and Definition \ref{Def:Barr}. 
\begin{definition}\label{Def:tau}
For $t>0$, $y,x<0$ and $n\in \N$, define
\begin{align}
h_t^{\smallfrown}&:[0,t]\to \R,\  s\mapsto -(((1+s)\wedge(1+t-s))^{1/6}-1), \notag\\
m_{t,h}^{\smallfrown}(s) &:= h_t^{\smallfrown}(s)-\frac{s}{t\vartheta^\ast}\log(p_t), \notag\\
\barrierprobUB &:= \mathbb{P}_{\mathcal{L}}\left[\mathcal{B}_{\{0,\dots,n\}, m_{n,h/2}^{\smallfrown}}^{y,J_x}( T_{\cdot}/\vartheta^\ast)\right], \label{Def:barrprobUB}
\end{align}
See Figure \ref{fig:UpBarr} for an illustration of the event in \eqref{Def:barrprobUB}. For $u$ with $|u| = n$, let
\[
\tau(u,y) := \min\left\{k\in \{0,\dots, n\} : V(u_k)-(\vartheta^\ast)^{-1} K_k+m_{n,h/2}^{\smallfrown}(k)+y+1\ge 0\right\}
\]
and $\tau(y) := \min_{u: |u| = n} \tau(u,y)\wedge n$. The dependency on $n$ of $\tau(y)$ is omitted from notation. Finally, define
\begin{align*}
A_n(y,x) &:= \bigcup_{u : |u|=n} \mathcal{B}_{\{0,\dots,n\}, m_{n,h/2}^{\smallfrown}}^{y,J_x}(V(u_\cdot)-(\vartheta^\ast)^{-1}K_\cdot).
\end{align*}
\end{definition}

It follows directly from the definitions that for all $y<0$
\begin{align}
\mathbb{P}_{\mathcal{L}}\left[M_n-m_n\in [\xi_0-y-1,\xi_0-y]\right]&\le \mathbb{P}_{\mathcal{L}}[\tau(y)<n]+\mathbb{P}_{\mathcal{L}}\left[A_n(y,\xi_0)\right]. \label{eq:Split}
\end{align}

It will be useful to be able to compare $\barrierprobUB$ from \eqref{Def:barrprobUB} to $p_n$ from \eqref{Def:pt}. The next Proposition \ref{Prop:BarrComp} allows just that. Proposition \ref{Prop:BarrComp} will be proved in Section \ref{Sec:BarrComps}, the term $|y|^{c_\varepsilon}|x-1|^{c_\varepsilon}$ in Proposition \ref{Prop:BarrComp} comes from moving the start- and endpoint. Going from discrete to continuous time and removing the curve of the barrier only costs a constant factor. We recall that $\xi_0$ is negative.

\begin{proposition}\label{Prop:BarrComp}
For all $\varepsilon>0$, there are $c_\varepsilon, C_\varepsilon>0$ such that
\begin{align*}
\liminf_{n\to \infty}\mathbb{P}\left[\bigcap_{x\in [-\log(n),\xi_0]\cap\Z\atop y\in [-\log(n)^{2},\xi_0]\cap \Z}\left\{\frac{\barrierprobUB}{p_n} \le C_\varepsilon|y|^{c_\varepsilon}|x-1|^{c_\varepsilon}\right\}\right]\ge 1-\varepsilon.
\end{align*} 
\end{proposition}

By \eqref{eq:Split}, proving Proposition \ref{Prop:EndingUp} is reduced to bounding $\mathbb{P}_{\mathcal{L}}[A_n(y,\xi_0)]$ and $\mathbb{P}_{\mathcal{L}}[\tau(y)<n]$, which we do in Propositions \ref{Prop:ControlAn} and \ref{Prop:ControlTau}. 
\begin{proposition}\label{Prop:ControlAn}
For all $\varepsilon>0$, there are $C_\varepsilon, c_\varepsilon>0$ such that 
\[
\liminf_{n\to \infty} \mathbb{P}\left[\bigcap_{y\in [-\log(n)^{2},\xi_0]\cap \Z}\left\{\mathbb{P}_{\mathcal{L}}[A_n(y,\xi_0)]\le C_\varepsilon e^{\vartheta^\ast y}|y|^{c_\varepsilon}\right\}\right]\ge 1-\varepsilon.
\] 
\end{proposition}

\begin{proof}[Proof of Proposition \ref{Prop:ControlAn} assuming Proposition \ref{Prop:BarrComp}]
Applying Lemma \ref{Lem:MT1} yields that $\mathbb{P}$-a.s.
\begin{align*}
\mathbb{P}_{\mathcal{L}}[A_n(y,\xi_0)] &= \mathbb{E}_{\mathcal{L}}\left[e^{-T_n};\mathcal{B}_{\{0,\dots,n\}, m_{n,h/2}^{\smallfrown}}^{y,J_{\xi_0}}( T_{\cdot}/\vartheta^\ast)\right]\le p_n^{-1}e^{\vartheta^\ast(1+y-\xi_0)}p_n^{\smallfrown,\xi_0}(y).
\end{align*}
Applying Proposition \ref{Prop:BarrComp} and pulling the factor $|\xi_0-1|^{c_\varepsilon}e^{\vartheta^\ast(1-\xi_0)}$ into the $C_\varepsilon$, which is possible, since $\xi_0$ is a fixed (negative) constant, finishes the proof.
\end{proof}

The proof of the forthcoming Proposition \ref{Prop:ControlTau} appears in Section \ref{Subsec:ProofControlTau}.
\begin{proposition}\label{Prop:ControlTau}
For all $\varepsilon>0$, there are $c_\varepsilon, C_\varepsilon>0$ such that 
\[
\liminf_{n\to\infty}\mathbb{P}\left[\bigcap_{y\in [-\log(n)^{2},\xi_0]\cap \Z}\left\{\mathbb{P}_{\mathcal{L}}[\tau(y)<n]\le C_\varepsilon e^{\vartheta^\ast y} |y|^{c_\varepsilon}\right\}\right]\ge 1-\varepsilon.
\]
\end{proposition}

\begin{proof}[Proof of Proposition \ref{Prop:EndingUp} assuming Proposition \ref{Prop:ControlTau}]
We combine \eqref{eq:Split}, Proposition \ref{Prop:ControlAn} and Proposition \ref{Prop:ControlTau}.
\end{proof}
\subsection{Lower Bound on the Right Tail of $M_n-m_n$}\label{Subsec:LB}
\begin{figure}
\includegraphics{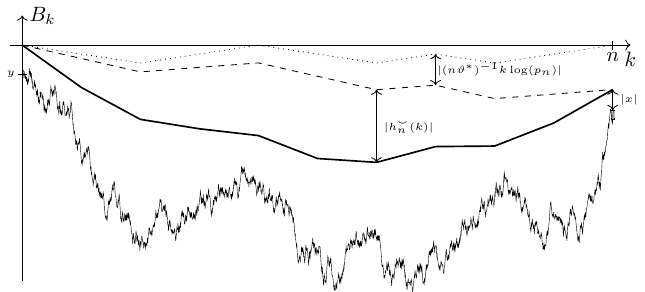}
\centering
\caption{The event in \eqref{Def:LBRTpt}. Drawn are $(K_k/\vartheta^\ast)_{k\le n}$ (dotted line), $(K_k/\vartheta^\ast+\frac{k}{n\vartheta^\ast} \log(p_n))_{k\le n}$ (dashed line), $(K_k/\vartheta^\ast+\frac{k}{n\vartheta^\ast}\log(p_n)-h_n^{\smallsmile}(k))_{k\le n}$ (thick black line) and a sample path of $(y+B_s)_{s\le n}$ realizing the barrier event $\mathcal{B}_{\{0,\dots, n\},m_{n,h/2}^{\smallsmile}}^{y,J_x}(T_\cdot/\vartheta^\ast)$. We have the parameters $n = 10$, $\mathcal{L}_1 \sim \mathrm{Unif}(\{2,3\})$, $(\mathcal{L}_k)_{k\le 10} = (2,2,3,3,2,2,3,2,3,3)$. }
\label{fig:DownBarr}
\end{figure}

\begin{theorem}\label{Theo:UBRT}
For all $\varepsilon>0$, there are $c_\varepsilon, C_\varepsilon>0$ such that for $y\le \xi_0$,
\[
\liminf_{n\to \infty} \mathbb{P}\left[\mathbb{P}_{\mathcal{L}}[M_n-m_n\ge \xi_0-y-1] \ge e^{\vartheta^\ast y} C_\varepsilon |y|^{-c_\varepsilon}\right]\ge 1-\varepsilon.
\]
\end{theorem}
We will prove Theorem \ref{Theo:UBRT} using a second moment approach. For this, we need the following definition, for which we recall the Definitions \ref{Def:Basics}, \ref{Def:Barr} and \eqref{Def:pt}.
\begin{definition}\label{Def:BasicsLB}
For $t>0$, $s\in [0,t]$,  $y<0$ and $n\in\N$, define
\begin{align}
h_t^{\smallsmile}(s) &:= ((1+s)\wedge (t-s+1))^{1/6}-1, \notag\\
m_{t,h,y}^{\smallsmile}(s) &:= y+h_t^{\smallsmile}(s)-\frac{s}{t\vartheta^\ast}\log(p_t),\quad m_{t,h}^{\smallsmile}(s) := m_{t,h,0}^{\smallsmile}(s),\notag\\
\barrierprobLB &:= \mathbb{P}_{\mathcal{L}}\left[\mathcal{B}_{\{0,\dots,n\}, m_{n,h,y}^{\smallsmile}}^{J_{\xi_0}}( T_{\cdot}/\vartheta^\ast)\right], \label{Def:LBRTpt}\\
Z_n(y) &:= \sum_{|u| = n}\mathbf{1}_{\mathcal{B}_{\{0,\dots,n\}, m_{n,h,y}^{\smallsmile}}^{J_{\xi_0}}( V(u_\cdot)-(\vartheta^\ast)^{-1}K_\cdot)}.\notag
\end{align}
See Figure \ref{fig:DownBarr} for an illustration of the event in \eqref{Def:LBRTpt}
\end{definition}

The random variable $Z_n$ counts the number of particles, which stay below the barrier $-m_{n,h,y}^{\smallsmile}(j)+(\vartheta^\ast)^{-1}K_{\cdot}$ and end up in $[m_n+\xi_0-y-1,m_n+\xi_0-y]$, in particular we have
\begin{align}
\mathbb{P}_{\mathcal{L}}[M_n-m_n\ge \xi_0-y-1] &\ge \mathbb{P}_{\mathcal{L}}[Z_n(y)\ge 1]\ge \frac{\mathbb{E}_{\mathcal{L}}[Z_n(y)]^2}{\mathbb{E}_{\mathcal{L}}[Z_n(y)^2]}\label{eq:CSS},
\end{align}
where the last step uses the Cauchy-Schwartz inequality. We will establish the following two Propositions \ref{Prop:LBZny}, \ref{Prop:Zn2UB}. We prove Proposition \ref{Prop:LBZny} in Section \ref{Subsec:LBoZny}. We sketch the proof of Proposition \ref{Prop:Zn2UB} in Section \ref{Subsec:Zn2UB} and fill in the detail in Section \ref{Subsec:Zn2UBDet}.

\begin{proposition}\label{Prop:LBZny}
For all $\varepsilon>0$, there are $c_\varepsilon, C_\varepsilon>0$ such that for all $y\le \xi_0$,
\[
\liminf_{n\to \infty} \mathbb{P}\left[\mathbb{E}_{\mathcal{L}}[Z_n(y)]\ge C_\varepsilon e^{\vartheta^\ast y}|y|^{-c_\varepsilon}\right]\ge 1-\varepsilon.
\]
\end{proposition}

\begin{proposition}\label{Prop:Zn2UB}
For all $\varepsilon>0$, there are $c_\varepsilon, C_\varepsilon>0$ such that for all $y\le \xi_0$,
\[
\liminf_{n\to \infty} \mathbb{P}\Big[\mathbb{E}_{\mathcal{L}}[Z_n(y)^2]\le C_\varepsilon e^{\vartheta^\ast y} |y|^{c_\varepsilon}+\mathbb{E}_{\mathcal{L}}[Z_n(y)]\Big]\ge 1-\varepsilon.
\]
\end{proposition}
\begin{proof}[Proof of Theorem \ref{Theo:UBRT} assuming Propositions \ref{Prop:LBZny} and \ref{Prop:Zn2UB}]

Combining \eqref{eq:CSS} and Proposition \ref{Prop:Zn2UB} yields that there are $c_{\varepsilon/2}, C_{\varepsilon/2}>0$ such that
\[
1-\varepsilon/2 \le \liminf_{n\to\infty} \mathbb{P}\left[\mathbb{P}_{\mathcal{L}}[M_n-m_n\ge \xi_0-y-1] \ge\frac{\mathbb{E}_{\mathcal{L}}[Z_n(y)]}{1+C_{\varepsilon/2} e^{\vartheta^\ast y} |y|^{c_{\varepsilon/2}}/\mathbb{E}_{\mathcal{L}}[Z_n(y)]}\right].
\]
Applying Proposition \ref{Prop:LBZny} yields that there are $C_{\varepsilon/2}'$, $c_{\varepsilon/2}'$ such that 
\begin{align*}
1-\varepsilon \le \liminf_{n\to\infty} \mathbb{P}\left[\mathbb{P}_{\mathcal{L}}[M_n-m_n\ge \xi_0-y-1] \ge \frac{C_{\varepsilon/2}'e^{\vartheta^\ast y}|y|^{-c_{\varepsilon/2}'}}{1+C_{\varepsilon/2}C_{\varepsilon/2}'e^{2\vartheta^\ast y} |y|^{c_{\varepsilon/2}+c_{\varepsilon/2}'}} \right]
\end{align*}
However, $C_{\varepsilon/2}C_{\varepsilon/2}'e^{2\vartheta^\ast y} |y|^{c_{\varepsilon/2}+c_{\varepsilon/2}'}$ is bounded uniformly in $y<0$, such that the last display implies the statement of Theorem \ref{Theo:UBRT}.
\end{proof}
\subsubsection{Lower Bound on $\mathbb{E}_{\mathcal{L}}[Z_n(y)]$}\label{Subsec:LBoZny}

In this section, we sketch the proof of Proposition \ref{Prop:LBZny}, it closely resembles the proof of Proposition \ref{Prop:ControlAn}, with the forthcoming Proposition \ref{Prop:barrprobLBLB} replacing  Proposition \ref{Prop:BarrComp}.

The proof has two steps, first applying the many-to-one Lemma \ref{Lem:MT1} gives that
\begin{equation}
\mathbb{E}_{\mathcal{L}}[Z_n(y)] \ge e^{\vartheta^\ast y-\vartheta^\ast \xi_0}p_n^{-1} \barrierprobLB\qquad\mathbb{P}\text{-a.s.} \label{eq:LBZNyMT1}
\end{equation}
Similar to Proposition \ref{Prop:BarrComp} we prove a lower bound on $\barrierprobLB/p_n$, recall \eqref{Def:pt} and \eqref{Def:LBRTpt}. The following Proposition \ref{Prop:barrprobLBLB} will be proved in Section \ref{Sec:BarrComps}.
\begin{proposition}\label{Prop:barrprobLBLB}
For all $\varepsilon>0$, there are $c_\varepsilon, C_\varepsilon>0$ such that for $y\le \xi_0$,
\[
\liminf_{n\to \infty}\mathbb{P}\left[\frac{\barrierprobLB}{p_n}\ge C_\varepsilon |y|^{-c_\varepsilon}\right]\ge 1-\varepsilon.
\]
\end{proposition}
\begin{proof}[Proof of Proposition \ref{Prop:LBZny} assuming Proposition \ref{Prop:barrprobLBLB}]
We combine \eqref{eq:LBZNyMT1} with Proposition \ref{Prop:barrprobLBLB}.
\end{proof}
\begin{remark}
The $|y|^{-c_\varepsilon}$ in Proposition \ref{Prop:barrprobLBLB} is far from optimal, with more effort one can prove that the statement holds (as one expects) with  $|y|^{c_\varepsilon}$ instead, however we do not need that.
\end{remark}

\subsubsection{Upper Bound on $\mathbb{E}_{\mathcal{L}}[Z_n(y)^2]$}\label{Subsec:Zn2UB} In this section, we outline the proof of Proposition \ref{Prop:Zn2UB}. 

The random variable $Z_n(y)^2$ roughly counts the number of pairs of particles $u,v$ in generation $n$, for which both particles stay below the barrier $-m_{n,h,y}^{\smallsmile}(\cdot)+(\vartheta^\ast)^{-1}K_\cdot$ and end up in $[m_n+\xi_0-y,m_n+\xi_0-y+1]$. We partition with regard to $|u\wedge v|$, which leads to the following definition.
\begin{definition}\label{Def:Lky}
For $n\in\N$ and $y\le \xi_0$, we define
\begin{align*}
\Lambda_k(y) &:= \sum_{|u| = k} \sum_{|u^1| = |u^2| = n\atop u^1\wedge u^2 = u} \prod_{i=1}^2\mathbf{1}_{\mathcal{B}_{\{0,\dots,n\}, m_{n,h,y}^{\smallsmile}(\cdot)}^{J_{\xi_0}}(V(u_\cdot^i)-(\vartheta^\ast)^{-1}K_\cdot)}.
\end{align*}
\end{definition}
We have that
\begin{equation}
Z_n(y)^2 = \sum_{k=0}^n \Lambda_k(y) =\sum_{k=0}^{n-1}\Lambda_k(y)+Z_n(y), \label{eq:SplitZn2Lk}
\end{equation}
where the second equality uses that $\Lambda_n(y) = Z_n(y)$. Thus all we need to do to prove Proposition \ref{Prop:Zn2UB}, is to bound $\sum_{k=0}^{n-1}\mathbb{E}_{\mathcal{L}}[\Lambda_k(y)]$.

To do this, we introduce some additional notation. We recall Definitions \ref{Def:Basics}, \ref{Def:BasicsLB} and the convention that $T_n^{(k)} =\vartheta^\ast\sum_{j=k+1}^{k+n} X_j-\kappa_j(\vartheta^\ast)$.
\begin{definition}\label{Def:fkqk}
For  $n\in\N$, $k\in \{0,\dots, n-1\}$, $y\le \xi_0$ and $x\in \R$ define
\begin{align}
f_k(x) &:= \mathbb{E}_{\mathcal{L}}^k\left[\sum_{|u| = n-k}\mathbf{1}_{\mathcal{B}_{\{0,\dots,n-k\},m_{n,h,y}^{\smallsmile}(k+\cdot)}^{x,J_{\xi_0}}(V(u_\cdot)-(\vartheta^\ast)^{-1}K_{k+\cdot})}\right], \notag\\
q_k(x) &:= \mathbb{P}_{\mathcal{L}}^k\Bigg[\mathcal{B}_{\{0,\dots,n-k\}, m_{n,h,y}^{\smallsmile}(k+\cdot)}^{x-K_k/\vartheta^\ast, J_{\xi_0}}(T_{\cdot}^{(k)}/\vartheta^\ast)\Bigg]. \label{Def:qk}
\end{align}
Both $q_k(x)$ and $f_k(x)$ depend on $n$ and $y$, which we omit from notation.
\end{definition}
 By conditioning on $\mathcal{F}_{k} := \sigma(u, V(u) : |u|\le k)$ we will prove the following upper bound on $\mathbb{E}_{\mathcal{L}}[\Lambda_k(y)]$ in Section \ref{Subsec:PropLky}.

\begin{proposition}\label{Prop:Lky}

For $n\in \N$, $k\in \{0,\dots, n-1\}$ and $y\le \xi_0$, we have that $\mathbb{P}$-a.s.
\begin{align*}
\mathbb{E}_{\mathcal{L}}[\Lambda_k(y)]&\le c e^{\vartheta^\ast y-2\vartheta^\ast \xi_0} p_n^{\frac{k}{n}-2}e^{-\vartheta^\ast h_n^{\smallsmile}(k)}\barrierprobLB \max_{x\ge0} e^{-\vartheta^\ast x}q_k(K_k/\vartheta^\ast-m_{n,h,y}^{\smallsmile}(k)-x).
\end{align*}

In particular, we have that
\begin{equation}
\mathbb{E}_{\mathcal{L}}[\Lambda_k(y)]\le ce^{\vartheta^\ast y-2\vartheta^\ast \xi_0}p_n^{\frac{k}{n}-2}e^{-\vartheta^\ast h_n^{\smallsmile}(k)}\barrierprobLB. \label{eq:Lkeasy}
\end{equation}
\end{proposition}
It will be useful to derive an upper bound on $\barrierprobLB/p_n$, recall \eqref{Def:pt} and \eqref{Def:LBRTpt}. This is done in the following proposition, which is proved in Section \ref{Sec:BarrComps}.
\begin{proposition}\label{Prop:BarrierUBLB}
For all $\varepsilon>0$, there are $c_\varepsilon, C_\varepsilon>0$ such that 
\[
\liminf_{n\to \infty}\mathbb{P}\left[\bigcap_{y\in [-\log(n)^2,\xi_0]\cap\Z \atop }\left\{\frac{\barrierprobLB}{p_n}\le C_\varepsilon |y|^{c_\varepsilon}\right\}\right]\ge 1-\varepsilon. 
\]
\end{proposition}
\begin{remark}
Proposition \ref{Prop:BarrierUBLB} also holds, with minor changes in the proof, when one replaces $J_{\xi_0}$ in \eqref{Def:LBRTpt} by $[\xi_0-2,\xi_0]$, this version will be used in Section \ref{Subsec:Lkhard} for the proof of Lemma \ref{Lem:UpperBoundqk}.
\end{remark}
The next two lemmata allow us to bound $\sum_{k=0}^{n-1} \mathbb{E}_{\mathcal{L}}[\Lambda_k]$.
\begin{lemma}\label{Lem:Lkeasy}
For all $\varepsilon>0$, there are $c_\varepsilon, C_\varepsilon>0$ such that for all $y\le \xi_0$,
\begin{equation}
\liminf_{n\to \infty}\mathbb{P}\left[\sum_{k=\lceil\log(n)^7\rceil}^{n-1}\mathbb{E}_{\mathcal{L}}[\Lambda_k(y)]\le C_\varepsilon e^{\vartheta^\ast y}|y|^{c_\varepsilon}\right]\ge 1-\varepsilon. \label{eq:ln6n1}
\end{equation}
\end{lemma}
\begin{proof}[Proof of Lemma \ref{Lem:Lkeasy} assuming Proposition \ref{Prop:BarrierUBLB}]
Applying Proposition \ref{Prop:BarrierUBLB} to equation \eqref{eq:Lkeasy} of Proposition \ref{Prop:Lky} yields that there are $c_\varepsilon, C_\varepsilon>0$ such that for $y\le \xi_0$ 
\begin{align}
\liminf_{n\to \infty}\mathbb{P}\left[\bigcap_{k\in [\log(n)^7,n-1]\cap \Z}\mathbb{E}_{\mathcal{L}}[\Lambda_k(y)]\le C_\varepsilon e^{\vartheta^\ast y} p_n^{\frac{k}{n}-1}e^{-\vartheta^\ast h_n^{\smallsmile}(k)}|y|^{c_\varepsilon}\right]\ge 1-\varepsilon.\label{eq:Lkln6n}
\end{align}
By Lemma \ref{Lem:Boundpnastweak} there is a $C>0$ such that $\mathbb{P}$-a.s.
$
\limsup_{n\to\infty}\frac{|\log(p_n)|}{\log(n)}\le C $, which together with \eqref{eq:Lkln6n} implies that there are $c_\varepsilon, C_\varepsilon>0$ such that for $y\le \xi_0$ 
\begin{align}
1-\varepsilon &\le \liminf_{n\to \infty}\mathbb{P}\left[\sum_{k=\lceil\log(n)^7\rceil}^{n-1} \mathbb{E}_{\mathcal{L}}[\Lambda_k(y)] \le C_\varepsilon e^{\vartheta^\ast y}|y|^{c_\varepsilon} \sum_{k=\lceil\log(n)^7\rceil}^{n-1}e^{-\vartheta^\ast h_n^{\smallsmile}(k)}e^{C \log(n)\left(1-\frac{k}{n}\right)}\right]. \label{eq:L116last}
\end{align}

We can bound $\sup_n\sum_{k=\lceil\log(n)^7\rceil}^{n-1}e^{-\vartheta^\ast h_n^{\smallsmile}(k)}e^{C \log(n)\left(1-\frac{k}{n}\right)}<\infty$, which together with \eqref{eq:L116last} finishes the proof.
\end{proof}

Handling $k<\log(n)^7$ in \eqref{eq:ln6n1} requires a more nuanced argument because the gain from $h_n^{\smallsmile}(k)$ in \eqref{eq:L116last} is not sufficient. We prove the following lemma in Section \ref{Subsec:Lkhard}.
\begin{lemma}\label{Lem:Lkhard}
For all $\varepsilon>0$, there are $c_\varepsilon, C_\varepsilon>0$ such that for all $y\le \xi_0$,
\[
\liminf_{n\to \infty} \mathbb{P}\left[\sum_{k=0}^{\lfloor \log(n)^7\rfloor} \mathbb{E}_{\mathcal{L}}[\Lambda_k(y)] \le C_\varepsilon e^{\vartheta^\ast y}|y|^{c_\varepsilon}\right]\ge 1-\varepsilon.
\]
\end{lemma}
\begin{proof}[Proof of Proposition \ref{Prop:Zn2UB}]
We combine \eqref{eq:SplitZn2Lk} with Lemma \ref{Lem:Lkeasy} and Lemma \ref{Lem:Lkhard}.
\end{proof}

\subsection{Bound on the Left Tail of $M_n-m_n$}\label{Subsec:LTail}

The aim of this section is to prove the following theorem.
\begin{theorem}\label{Theo:Main2}
For all $\varepsilon>0$, there is a $y_\varepsilon\le  0$ such that 
\[
\liminf_{n\to\infty}\mathbb{P}\left[\mathbb{P}_{\mathcal{L}}\left[M_n-m_n\le y_{\varepsilon}\right]\le \varepsilon\right]\ge 1-\varepsilon.
\]
\end{theorem}

The idea is, to cut the tree at depth $l$, use that there are at least $2^l$ particles in generation $l$, and then to finish the argument by applying the lower bound on the right tail proved in Theorem \ref{Theo:UBRT}. 

Slightly more detailed, we note that, similarly to Lemma \ref{Lem:Ausreißer}, we get the following lemma.
\begin{lemma}\label{Lem:MinBRW}
For all $\varepsilon>0$, there are $c_\varepsilon', c_\varepsilon''>0$ such that for all $l\in\N$, 
\[
\mathbb{P}\left[\mathbb{P}_{\mathcal{L}}\left[\min_{u: |u| = l} V(u)\le -c_\varepsilon' l\right]\le e^{-c_{\varepsilon}''l}\right]\ge 1-\varepsilon.
\] 
\end{lemma}
Then by cutting the tree at depth $l$ one can make the calculation that there are $c_\varepsilon, c_{\varepsilon}', c_{\varepsilon}''$ and $C_\varepsilon$ such that
\begin{align*}
\liminf_{n\to \infty}\mathbb{P}\left[\mathbb{P}_{\mathcal{L}}[M_n-m_{n-l}\le y]\le e^{-c''_\varepsilon l}+\left(1-C_\varepsilon |y+c_\varepsilon'l|^{-c_\varepsilon}\right)^{2^l}\right]\ge 1-\varepsilon,
\end{align*}
and by choosing $l,y$ right, but importantly independent of $n$, we get that $e^{-c''_\varepsilon l}+\left(1-C_\varepsilon |y+c_\varepsilon'l|^{-c_\varepsilon}\right)^{2^l}\le \varepsilon$, which allows us to conclude Theorem \ref{Theo:Main2}. In the final calculation we can't use $m_{n-l}$, but will instead need to use time-shifted objects, also it will be necessary to prove, that $m_{n}-m_{n-l}$ isn't too big. All of this is done in Section \ref{Sec:Finish}.

\subsection{Proof of Theorem \ref{Theo:Main} assuming Theorems \ref{Theo:MainUB} and \ref{Theo:Main2}}\label{Sec:ProofMainTheo}
Let $\varepsilon>0$ and $y\ge0$. Set 
\[
A_{n,\varepsilon}(y) := \left\{\mathbb{P}_{\mathcal{L}}[|M_n-m_n|\ge y] \le \varepsilon/2 \right\}.
\]
By applying Theorems \ref{Theo:MainUB} and \ref{Theo:Main2} there is a $y_\varepsilon\in\R$ such that
\begin{equation}
\limsup_{n\to \infty} \mathbb{P}\left[A_{n,\varepsilon}(y_\varepsilon)^c\right] \le \varepsilon/2. \label{eq:UBPAne}
\end{equation}
Furthermore, we have that
\begin{align}
P[|M_n-m_n| \ge y_\varepsilon] &= \mathbb{E}\left[\mathbb{P}_{\mathcal{L}}\left[|M_n-m_n|\ge y_\varepsilon\right]\mathbf{1}_{A_{n,\varepsilon}(y_\varepsilon)}\right]\notag\\
&\qquad+\mathbb{E}\left[\mathbb{P}_{\mathcal{L}}\left[|M_n-m_n|\ge y_\varepsilon\right]\mathbf{1}_{A_{n,\varepsilon}(y_\varepsilon)^c}\right]\notag\\
&\le \varepsilon/2+\mathbb{P}\left[A_{n,\varepsilon}(y_\varepsilon)^c\right]\label{eq:UBPMnmn}
\end{align}
Combining \eqref{eq:UBPAne} and \eqref{eq:UBPMnmn} yields that
\[
\limsup_{n\to\infty}P[|M_n-m_n| \ge y_\varepsilon]\le \varepsilon,
\]
which gives that $(M_n-m_n)_{n\ge0}$ is tight w.r.t.\@ $P$.\qed
\section{Details for the Proof of Theorem \ref{Theo:MainUB} -- Upper Bound on the Right Tail of $M_n-m_n$}\label{Sec:UB}
\subsection{Proof of Proposition \ref{Prop:ControlTau}}\label{Subsec:ProofControlTau}

We recall Definition \ref{Def:tau}. Proposition \ref{Prop:ControlTau} follows directly from the following two lemmata.

\begin{lemma}\label{Lem:tauklein}
For all $\varepsilon>0$, there is a $C_\varepsilon>0$ such that 
\[
\liminf_{n\to \infty} \mathbb{P}\left[\bigcap_{y\in [-\log(n)^{2},\xi_0]\cap\Z}\left\{\mathbb{P}_{\mathcal{L}}\left[\tau(y)\le n-\log(n)^7\right]\le  C_\varepsilon e^{\vartheta^\ast(y+1)}\right\}\right]\ge 1-\varepsilon.
\]
\end{lemma}

\begin{lemma}\label{Lem:MainWork}
For all $\varepsilon>0$, there are $C_\varepsilon, c_\varepsilon>0$ such that 
\begin{align*}
\liminf_{n\to \infty} \mathbb{P}\left[\bigcap_{y\in [-\log(n)^{2},\xi_0]\cap\Z}\left\{\mathbb{P}_{\mathcal{L}}\left[n-\log(n)^7<\tau(y)< n\right]\le C_\varepsilon e^{\vartheta^\ast y}|y|^{c_\varepsilon}\right\}\right]\ge 1-\varepsilon.
\end{align*}
\end{lemma}

\begin{proof}[Proof of Lemma \ref{Lem:tauklein}]
For all $y<0$ and all $k\le n-\log(n)^7$ we have that $\mathbb{P}$-a.s.
\begin{align}
\mathbb{P}_{\mathcal{L}}[\tau(y) = k] &\le \mathbb{P}_{\mathcal{L}}\left[\sum_{u : |u| = k} \mathbf{1}_{\{V(u)-(\vartheta^\ast)^{-1} K_k+m_{n,h/2}^{\smallfrown}(k)+y+1\ge 0\}}\ge 1 \right]\notag\\
&\stackrel{L.\ref{Lem:MT1}}{\le} \mathbb{E}_{\mathcal{L}}\left[e^{-T_k}; T_k/\vartheta^\ast+m_{n,h/2}^{\smallfrown}(k)+y+1\ge 0\right]\notag\\
&\le e^{\vartheta^\ast(y+1)}e^{\vartheta^\ast h_n^{\smallfrown}(k)/2}p_n^{-\frac{k}{n}}. \label{eq:tkmt1}
\end{align}

Combining Lemma \ref{Lem:Boundpnastweak} with \eqref{eq:tkmt1} gives that there is a $C_{\varepsilon}>0$ such that
\begin{align}
&\liminf_{n\to\infty}\mathbb{P}\left[\forall_{y\in [-\log(n)^{2},\xi_0]\cap\Z}\, \mathbb{P}_{\mathcal{L}}[\tau(y)\le n-\log(n)^7]\le  e^{\vartheta^\ast(y+1)}\hskip-3pt\sum_{k=1}^{\lfloor n-\log(n)^7\rfloor}\hskip-3pte^{\vartheta^\ast h_n^{\smallfrown}(k)/2}n^{C_\varepsilon\frac{k}{n}}\right]\notag
\\&\ge 1-\varepsilon. \label{eq:taunklein}
\end{align}
From $ \sup_{n\in \N}\sum_{k=1}^{\lfloor n-\log(n)^7\rfloor} e^{\vartheta^\ast h_n^{\smallfrown}(k)/2}n^{C_\varepsilon \frac{k}{n}}<\infty$ for all $C_\varepsilon>0$  and \eqref{eq:taunklein} we can conclude the claim of the lemma.
\end{proof}
The rest of this section is devoted to the proof of Lemma \ref{Lem:MainWork}. First, we need an additional definition, after this we will give a brief overview over the structure of the proof of Lemma \ref{Lem:MainWork}. We recall the definition $T_n^{(k)} =\vartheta^\ast\sum_{j=k+1}^{k+n} X_j-\kappa_j(\vartheta^\ast)$.

\begin{definition}\label{Def:tkhard}
For $n\in\N$, $k\le n$ and $y\le 0$ we set
\begin{align}
N_k(y) &:= \Bigg\{ u \in \mathbf{T}, |u| =n-k: \forall_{j\le n-k-1}\, V(u_j)-(\vartheta^\ast)^{-1}K_j+m_{n,h/2}^{\smallfrown}(j)+y+1\le 0, \notag\\
&\qquad V(u_{n-k})-(\vartheta^\ast)^{-1} K_{n-k}+m_{n,h/2}^{\smallfrown}(n-k)+y+1\ge 0\Bigg\},\label{Def:Nky}\\
\widetilde{T}_j(y) &:= T_j/\vartheta^\ast+m_{n,h/2}^{\smallfrown}(j)+y+1, \notag\\
q_k(y) &:= \mathbb{P}_{\mathcal{L}}\left[\forall_{j\le n-k-1}\,\widetilde{T}_j(y)\le 0, \widetilde{T}_j(y)\ge 0\right]. \notag
\end{align}
The dependence on $n$ is omitted from the notation.

For $l\in\N_0$, $x<0$ we define
\begin{align}
\widetilde{q}_k(l) &:=\mathbb{P}_{\mathcal{L}}\left[\mathcal{B}_{\{0,\dots,n-k-1\}, m_{n,h/2}^{\smallfrown}}^{y+1,J_{-l}}(T_{\cdot}/\vartheta^\ast)\right],\label{Def:qktild}\\
q_{k,\text{end}}(x) &:=\mathbb{P}_{\mathcal{L}}^{n-k-1}\Bigg[\mathcal{B}_{\{0,\dots,k+1\}, m_{n,h/2}^{\smallfrown}(n-k-1+\cdot)}^{x-m_{n,h/2}^{\smallfrown}(n-k-1), J_{\xi_0}}(T_{\cdot}^{(n-k-1)}/\vartheta^\ast) \Bigg].\label{Def:qkend}
\end{align}
The dependence on $n$ and $y$ is omitted from the notation, we emphasize, that $q_{k,\text{end}}(x)$ doesn't depend on $y$. 
\end{definition}

From Definition \ref{Def:tau} and \eqref{Def:Nky} it is immediate that $\mathbb{P}_{\mathcal{L}}\left[\tau(y) = n-k\right]\le \mathbb{P}_{\mathcal{L}}\left[|N_k(y)|\ge 1\right]$. This together with the Markov inequality and the many-to-one Lemma \ref{Lem:MT1} yields that

\begin{equation}
\mathbb{P}_{\mathcal{L}}\left[\tau(y) = n-k\right] \le p_n^{-\frac{n-k}{n}}e^{\vartheta^\ast(y+1)+\vartheta^\ast h_n^{\smallfrown}(k)/2} q_k(y),\label{eq:qk}
\end{equation}
where we used that $h_n^{\smallfrown}(n-k) = h_n^{\smallfrown}(k)$.  The idea is to decompose $q_k(y)$ according to the position of $\widetilde{T}_{n-k-1}$. Lemmata \ref{Lem:UBDeltaT} and \ref{Lem:qkyqk'} will establish, that 
\begin{equation}
q_k(y)\le \sum_{l=0}^\infty e^{-c_\varepsilon l^2} \widetilde{q}_k(l), \label{eq:qkyStrukt}
\end{equation}
i.e.\@ that $\widetilde{T}_{n-k-1} \approx l$ costs an $e^{-c_\varepsilon l^2}$ factor. The term $\widetilde{q}_k(l)$ then needs to be compared to $p_n^{-1}$, recall \eqref{Def:pt}, for this we want to use the barrier computation of Proposition \ref{Prop:BarrComp}. To do this we need to compute an upper bound on $\widetilde{q}_k(l)(p_n^{\smallfrown,-l}(y))^{-1}$, recall \eqref{Def:barrprobUB}. This is done in Corollary \ref{Cor:ApplBarr1}, which is proved using Lemma \ref {Lem:BarrSplit}. Corollary \ref{Cor: pnqk} handles the summation over $l$ in \eqref{eq:qkyStrukt}. Finally, we need one more barrier computation, which is stated in Lemma \ref{Lem:Ga} and proved in Section \ref{Sec:BarrComps}.
 
\begin{lemma}\label{Lem:UBDeltaT}
For all $\varepsilon>0$, there are $c, C>0$ such that for all $x\ge 0$, $k\in\N$,
\begin{align*}
\mathbb{P}_{\mathcal{L}}\left[(T_{n-k}-T_{n-k-1})/\vartheta^\ast \ge x\right]\le Ce^{-cx^2}.
\end{align*}
\end{lemma}
\begin{proof}
We have that
\begin{align*}
\mathbb{P}_{\mathcal{L}}\left[(T_{n-k}-T_{n-k-1})/\vartheta^\ast\ge x\right] &= \mathbb{P}_{\mathcal{L}}\left[ X_{n-k}-(\vartheta^\ast)^{-1}\log(\mathcal{L}_{n-k})-\frac{(\vartheta^\ast)}{2}\ge x\right]\\
&= \mathbb{P}_{\mathcal{L}}\left[ X_1-\frac{\vartheta^\ast}{2}\ge x+(\vartheta^\ast)^{-1}\log(\mathcal{L}_{n-k})\right]\\
&\le e^{-\frac{(x+(\vartheta^\ast)^{-1}\log(\mathcal{L}_{n-k})-\vartheta^\ast/2)^2}{2}}\le Ce^{-c x^2}
\end{align*}
for some $c,C>0$, since $\log(\mathcal{L}_{n-k})\ge 0$. The second to last step has used that $X_1-\vartheta^\ast/2$ is $\mathcal{N}(\vartheta^\ast/2,1)$ distributed.
\end{proof}
\begin{lemma}\label{Lem:qkyqk'}
For all $\varepsilon>0$, there is a $c_{\varepsilon,1}>0$ such that
\[
\liminf_{n\to \infty} \mathbb{P}\left[\forall_{y\in [-\log(n)^{2},\xi_0]\cap \Z}\forall_{k\in \{1,\dots, \lfloor\log(n)^7\rfloor\}}\, q_k(y)\le \sum_{l=0}^\infty \widetilde{q}_k(l)Ce^{-c_{\varepsilon,1} l^2}\right]\ge 1-\varepsilon.
\]
\end{lemma}
\begin{proof}
This proof is mainly partitioning with respect to the value of $\widetilde{T}_{n-k-1}$ and using that by Lemma \ref{Lem:UBDeltaT} having $\widetilde{T}_{n-k-1} \approx -l$ and $\widetilde{T}_{n-k}\ge 0$ costs roughly $e^{-cl^2}$, where considering $\widetilde{T}_k$ instead of $T_k/\vartheta^\ast$ will only change the $c$.

We have that $\mathbb{P}$-a.s.\@ for all $y\in [-\log(n)^{2},\xi_0]\cap \Z$, $k\in\{1,\dots,\lfloor \log(n)^7\rfloor\}$
\begin{align}
q_k(y) &= \sum_{l=0}^\infty \mathbb{P}_{\mathcal{L}}\left[\forall_{j\le n-k-1}\,\widetilde{T}_j\le 0, \widetilde{T}_{n-k}\ge 0,\widetilde{T}_{n-k-1}\in [-l-1,-l]\right]\notag\\
&\le \sum_{l=0}^\infty \mathbb{P}_{\mathcal{L}}\left[\forall_{j\le n-k-1}\,\widetilde{T}_j\le 0, \widetilde{T}_{n-k-1}\in [-l-1,-l],\widetilde{T}_{n-k}-\widetilde{T}_{n-k-1} \ge l\right]\notag\\
&\stackrel{L. \ref{Lem:UBDeltaT}, D. \ref{Def:tkhard}}{\le} \sum_{l=0}^\infty \widetilde{q}_k(l)Ce^{-c\left(\log(p_n)/(n\vartheta^\ast)+1/2( k^{1/6}- (k-1)^{1/6})+l\right)^2}. \label{eq:ImUBqk}
\end{align}
The statement of the lemma follows from \eqref{eq:ImUBqk} and Lemma \ref{Lem:Boundpnastweak}.
\end{proof}

We recall the definitions \eqref{Def:barrprobUB}, \eqref{Def:qktild} and \eqref{Def:qkend}.
\begin{corollary}\label{Cor:ApplBarr1}
Fix $k\le n$. For $0\le l\le |\xi_0|$, we have that
\begin{equation}
\widetilde{q}_k(l)\inf_{x\in [\xi_0-2,-1]} q_{k,\text{end}}(x) \le p_n^{\smallfrown,\xi_0}(y).\label{eq:AplBar11}
\end{equation}
For $l\ge |\xi_0|$ we have that
\begin{equation}
\widetilde{q}_k(l)\inf_{x\in [\xi_0-2,-1]}q_{k,\text{end}}(x) \le p_n^{\smallfrown,-l}(y).\label{eq:AplBar12}
\end{equation}
\end{corollary}
\begin{proof}

First, let $l\le|\xi_0|$.  By applying \eqref{eq:BarrSplit1} from Lemma \ref{Lem:BarrSplit} for $Z_{\cdot} = T_{\cdot}/\vartheta^\ast+m_{n,h/2}^{\smallfrown}$, $z_1 = y$, $z_2 = \xi_0$, $I= \{0,\dots,n\}$, $t_0 = n-k-1$ and $x_0 =-l$ we get, that  
\begin{align}
p_n^{\smallfrown,\xi_0}(y)  &\ge\widetilde{q}_k(l)\inf_{x\in J_{x_0-1}}\mathbb{P}_{\mathcal{L}}\left[\mathcal{B}_{\{n-k-1,\dots,n\},m_{n,h/2}^{\smallfrown}-m_{n,h/2}^{\smallfrown}(n-k-1)}^{x,J_{\xi_0}}( (T_\cdot-T_{n-k-1})/\vartheta^\ast )\right]\notag\\
&= \widetilde{q}_k(l) \inf_{x\in[-l-2,-l-1]} q_{k,\text{end}}(x), \label{eq:ErstSchrit1}
\end{align}
where the last step follows from \eqref{Def:qkend} and a time shift.  Since $l\le |\xi_0|$ equation \eqref{eq:ErstSchrit1} implies \eqref{eq:AplBar11}.

The case $l\ge|\xi_0|$, i.e.\@ proving \eqref{eq:AplBar12}, is analogous, but we take $z_2 =-l$ instead, thus we get
\begin{align*}
p_n^{\smallfrown,-l}(y) &\ge  \widetilde{q}_k(l)\inf_{x\in[-l-2,-l-1]}\mathbb{P}_{\mathcal{L}}^{n-k-1}\left[\mathcal{B}_{\{0,\dots,k+1\},m_{n,h/2}^{\smallfrown}(n-k-1+\cdot)}^{x-m_{n,h/2}^{\smallfrown}(n-k-1),J_{-l}}(T_{\cdot}/\vartheta^\ast)\right]\\
&=\widetilde{q}_k(l)\inf_{x\in[\xi_0-2,\xi_0-1]}\mathbb{P}_{\mathcal{L}}^{n-k-1}\left[\mathcal{B}_{\{0,\dots,k+1\},m_{n,h/2}^{\smallfrown}(n-k-1+\cdot)}^{x-\xi_0-l-m_{n,h/2}^{\smallfrown}(n-k-1),J_{-l}}(T_{\cdot}/\vartheta^\ast)\right]\\
&\ge \widetilde{q}_k(l)\inf_{x\in[\xi_0-2,\xi_0-1]} q_{k,\text{end}}(x),
\end{align*}
where the last step used that since $l\ge |\xi_0|$ we have $\xi_0+l\ge 0$. Now, observing that $[\xi_0-2,\xi_0-1] \subseteq [\xi_0-2,-1]$ finishes the proof.
\end{proof}
\begin{corollary}\label{Cor: pnqk}
For all $\varepsilon>0$, there are $C_\varepsilon, c_\varepsilon>0$ such that  
\begin{align*}
&\liminf_{n\to \infty}\mathbb{P}\left[\forall_{y\in [-\log(n)^{2},\xi_0]\cap\Z}\forall_{k\in\{1,\dots, \lfloor \log(n)^7\rfloor\}}\, p_n^{-1}q_k(y) \le \left(\inf_{x\in [\xi_0-2,-1]} q_{k,\text{end}}(x)\right)^{-1}C_\varepsilon|y|^{c_\varepsilon}\right]\\
&\qquad\ge 1-\varepsilon.
\end{align*}
\end{corollary}
\begin{proof}
To shorten the displays in this proof we set 
\[
I_n := ([-\log(n)^2,\xi_0]\cap \Z) \times \{1,\dots, \lfloor \log(n)^7\rfloor\}.
\]
Furthermore, in this proof $\sum_{j=a}^b$ is to be read as $\sum_{j = \lceil a\rceil}^{\lfloor b\rfloor}$ in case that $a,b$ aren't integers.

With $c_{\varepsilon,1}$ from Lemma \ref{Lem:qkyqk'} we have that
\begin{equation}
\liminf_{n\to \infty} \mathbb{P}\left[\forall_{(y,k)\in I_n}\, p_n^{-1}q_k(y)\le \sum_{l=0}^\infty p_n^{-1}\widetilde{q}_k(l)Ce^{-c_{\varepsilon,1} l^2}\right]\ge 1-\varepsilon. \label{eq:seed}
\end{equation}

We first handle the summation from $l =0$ to $|\xi_0|$ in \eqref{eq:seed}. By \eqref{eq:AplBar11} we have that for all $(y,k)\in I_n$
\begin{align}
\sum_{l=0}^{|\xi_0|} p_n^{-1}\widetilde{q}_k(l)Ce^{-c_{\varepsilon,1} l^2} &\le \left(\inf_{x\in [\xi_0-2,-1]} q_{k,\text{end}}(x)\right)^{-1}\sum_{l=0}^{|\xi_0|} p_n^{\smallfrown,\xi_0}(y)p_n^{-1}Ce^{-c_{\varepsilon,1} l^2}. \label{eq:lxi0}
\end{align}
By applying Proposition \ref{Prop:BarrComp} and using that $e^{-c_{\varepsilon,1} l^2}$ is summable, \eqref{eq:lxi0} implies that there are $c_\varepsilon$, $C_\varepsilon$ such that
\begin{align}
&1-\varepsilon/3\notag\\
& \le \liminf_{n\to\infty}\mathbb{P}\left[\forall_{(y,k)\in I_n}\,\sum_{l=0}^{|\xi_0|} p_n^{-1}\widetilde{q}_k(l)Ce^{-c_{\varepsilon,1} l^2}\le \left(\inf_{x\in [\xi_0-2,-1]} q_{k,\text{end}}(x)\right)^{-1}C_\varepsilon|y|^{c_\varepsilon}\right] \label{eq:qktild1}
\end{align}

Now we handle the summation from $|\xi_0|$ to $\sqrt{\log(n)\lambda/c_{\varepsilon,1}}$ in \eqref{eq:seed}, where $\lambda>0$ is such that $p_n^{-1}n^{-\lambda}$ converges to zero $\mathbb{P}$-a.s.\@ according to Lemma \ref{Lem:Boundpnastweak}.

 A completely parallel argument to the proof of \eqref{eq:qktild1}, but using the summability of $e^{-c_{\varepsilon,1}l^2}l^{c_\varepsilon}$,  yields that there are $c_\varepsilon, C_\varepsilon>0$ such that
\begin{align}
&1-\varepsilon/3 \notag\\\
&\le \liminf_{n\to\infty}\mathbb{P}\left[\forall_{(y,k)\in I_n}\,\hskip-0.575cm\sum_{l=|\xi_0|}^{\sqrt{\log(n)\lambda/c_{\varepsilon,1}}} \hskip-0.4cmp_n^{-1}\widetilde{q}_k(l)Ce^{-c_{\varepsilon,1} l^2}\le \left(\inf_{x\in [\xi_0-2,-1]} q_{k,\text{end}}(x)\right)^{-1}\hskip-0.15cmC_\varepsilon|y|^{c_\varepsilon}\right] \label{eq:qktild2}
\end{align}

Finally, we handle $l> \sqrt{\log(n)\lambda/c_{\varepsilon,1}}$ in \eqref{eq:seed}. Since $\widetilde{q}_k(l)$ is a probability it is smaller than $1$, which yields that 

\begin{align}
\sum_{l= \sqrt{\log(n)\lambda/c_{\varepsilon,1}}}^{\infty} p_n^{-1}\widetilde{q}_k(l)Ce^{-c_{\varepsilon,1} l^2} &\le C_\varepsilon p_n^{-1}e^{-\lambda \log(n)} \label{eq:qktild3}
\end{align}

Since $p_n^{-1}e^{-\lambda \log(n)}\to 0$ in probability, equations \eqref{eq:seed}, \eqref{eq:qktild1}, \eqref{eq:qktild2} and \eqref{eq:qktild3} finish the proof.
\end{proof}

We need one more Lemma for the proof of Lemma \ref{Lem:MainWork}, the proof of this lemma involves barrier computations and thus is done in Section \ref{Sec:BarrComps}. We recall the Definition \eqref{Def:qkend} of $q_{k,\text{end}}(x)$.
\begin{lemma}\label{Lem:Ga}
For all $\varepsilon>0$, there is a $C_\varepsilon>0$ such that 
\begin{align*}
\liminf_{n\to\infty}\mathbb{P}\left[\sum_{k=1}^{\lfloor \log(n)^7\rfloor} e^{\vartheta^\ast h_n^{\smallfrown}(k)/2}\left(\inf_{x\in [\xi_0-2,-1]} q_{k,\text{end}}(x)\right)^{-1} \le C_\varepsilon\right]\ge 1-\varepsilon.
\end{align*}
\end{lemma}

We have now collected all ingredients necessary to finish the proof of Lemma \ref{Lem:MainWork}
\begin{proof}[Proof of Lemma \ref{Lem:MainWork} assuming Lemma \ref{Lem:Ga}]

Applying \eqref{eq:qk} yields that $\mathbb{P}$-a.s.\@ for all $y\in [-\log(n)^{2},\xi_0]\cap\Z$
\begin{equation}
\mathbb{P}_{\mathcal{L}}\left[n-\log(n)^7<\tau(y)<n\right] \le \sum_{k=1}^{\lfloor \log(n)^7\rfloor} e^{\vartheta^\ast(y+1)} p_n^{-1}q_k(y)e^{\vartheta^\ast h_n^{\smallfrown}(k)/2}. \label{eq:tynln}
\end{equation}
Combining Corollary \ref{Cor: pnqk} and Lemma \ref{Lem:Ga} implies that there are $C_\varepsilon, c_\varepsilon>0$ such that 
\begin{equation}
1-\varepsilon \le \liminf_{n\to\infty} \mathbb{P}\left[\forall_{y\in [-\log(n)^{2},\xi_0]\cap\Z}\,\sum_{k= 1}^{\lfloor \log(n)^7\rfloor} p_n^{-1}q_k(y)e^{\vartheta^\ast h_n^{\smallfrown}(k)/2}\le C_\varepsilon|y|^{c_\varepsilon}\right]\label{eq:sumtynln}
\end{equation}
Equations \eqref{eq:tynln} and \eqref{eq:sumtynln} together imply Lemma \ref{Lem:MainWork}.
\end{proof}

\section{Details for Theorem \ref{Theo:UBRT} -- Lower Bound on the Right Tail of $M_n-m_n$ }\label{Sec:DetLB}

\subsection{Proof of Proposition \ref{Prop:Zn2UB}} \label{Subsec:Zn2UBDet}
\subsubsection{Proof of Proposition \ref{Prop:Lky}}\label{Subsec:PropLky}
In this section, fix $n\in\N$, $k\in\{0,\dots, n-1\}$ and $y\le \xi_0$. We recall Definitions \ref{Def:BasicsLB} and \ref{Def:fkqk}. Furthermore, define
\begin{align*}
\widetilde{V}(u,j) &:= V(u_j)-K_j/\vartheta^\ast+m_{n,h,y}^{\smallsmile}(j),\\
A_k(u) &:= \left\{\forall_{j\le k}\, \widetilde{V}(u, j) \le 0\right\},\\
\widetilde{T}_j &:= T_j/\vartheta^\ast+m_{n,h,y}^{\smallsmile}(j),\\
\mathcal{F}_k &:= \sigma(u,V(u) : |u|\le k).
\end{align*}
We will use two results on $f_k$ and $q_k$ in the proof of Proposition \ref{Prop:Lky}. First, as a direct consequence of the many-to-one Lemma \ref{Lem:MT1} applied to $\mathbb{P}_{\mathcal{L}}^k$ we have that for all $x\in \R$
\begin{equation}
f_k(x)\le ce^{\vartheta^\ast x+\vartheta^\ast y-\vartheta^\ast \xi_0-K_k}p_n^{-1}q_k(x)\qquad\mathbb{P}\text{-a.s.} \label{eq:UBqk}
\end{equation}
Furthermore, by conditioning on $\mathcal{F}_k$ it follows directly from the definitions, that for all $y<0$
\begin{equation}
\barrierprobLB =  \mathbb{E}_{\mathcal{L}}\left[q_{k}((T_k/\vartheta^\ast+K_k)/\vartheta^\ast);\forall_{j\le k}\, \widetilde{T}_j\le 0\right]\qquad \mathbb{P}\text{-a.s.} \label{eq:qkInter}
\end{equation}

\begin{proof}[Proof of Proposition \ref{Prop:Lky}]
We first use the tree structure by conditioning on $\mathcal{F}_k$, which gives that
\begin{align*}
\mathbb{E}_{\mathcal{L}}[\Lambda_k(y)|\mathcal{F}_k]&\le c\sum_{|u| = k} \mathbf{1}_{A_k(u)}f_{k}(V(u))^2.
\end{align*}
The $\mathbf{1}_{A_k(u)}$ just means that $u$ stays below the barrier up until time $k$, and the $f_k(V(u))^2$ is an upper bound for the number of pairs of particles, which descend from $u$ for which both particles stay below the barrier in $\{k,\dots, n\}$ and end up at the right place. In this calculation, we use that given $V(u)$ and $\mathcal{L}$ the events that two different particles behave this way are independent.

Applying \eqref{eq:UBqk} to bound $f_k(V(u))$ and taking expectation gives that $\mathbb{P}$-a.s.
\[
\mathbb{E}_{\mathcal{L}}[\Lambda_k(y)] \le c\Gamma_k\mathbb{E}_{\mathcal{L}}\left[\sum_{|u| = k}\mathbf{1}_{A_k(u)}\left(e^{\vartheta^\ast V(u)-K_k}q_k(V(u))\right)^2\right],
\]
where $\Gamma_k := e^{2\vartheta^\ast y-2\vartheta^\ast \xi_0}p_n^{-2}$.

If we wouldn't have the square, applying the many-to-one Lemma \ref{Lem:MT1} and \eqref{eq:qkInter} would allow us to bound $\mathbb{E}_{\mathcal{L}}[\Lambda_k]$ by $c\Gamma_k\barrierprobLB$; so we want to pull the factor $e^{\vartheta^\ast V(u)-K_k}q_k(V(u))$ out of the expectation. For this, we bound it by $\max_{x\ge0} e^{-\vartheta^\ast x-\vartheta^\ast m_{n,h,y}^{\smallsmile}(k)}q_k(K_k/\vartheta^\ast-m_{n,h,y}^{\smallsmile}(k)-x)$, which gives that $\mathbb{P}$-a.s.
\begin{align*}
\mathbb{E}_{\mathcal{L}}[\Lambda_k(y)] &\le c\Gamma_ke^{-\vartheta^\ast y-\vartheta^\ast h_n^{\smallsmile}(k)}p_n^{\frac{k}{n}}\max_{x\ge 0} e^{-\vartheta^\ast x}q_k(K_k/\vartheta^\ast-m_{n,h,y}^{\smallsmile}(k)-x)\\
&\qquad\mathbb{E}\left[\sum_{|u| = k} \mathbf{1}_{A_k(u)}e^{\vartheta^\ast V(u)-K_k}q_k(V(u))\right].
\end{align*}
Applying Lemma \ref{Lem:MT1} then gives that $\mathbb{P}$-a.s.
\begin{align*}
\mathbb{E}_{\mathcal{L}}[\Lambda_k(y)] &\le  c\Gamma_k e^{-\vartheta^\ast y-\vartheta^\ast h_n^{\smallsmile}(k)}p_n^\frac{k}{n}\mathbb{E}_{\mathcal{L}}\left[q_k((T_k/\vartheta^\ast+K_k)/\vartheta^\ast);\forall_{j\le k}\, \widetilde{T}_j\le 0\right]\\
&\qquad\max_{x\le 0}e^{-\vartheta^\ast x}q_k(K_k/\vartheta^\ast-m_{n,h,y}^{\smallsmile}(k)-x),
\end{align*}
Applying \eqref{eq:qkInter} finishes the proof of Proposition \ref{Prop:Lky} .\qedhere
\end{proof}

\subsubsection{Proof of Lemma \ref{Lem:Lkhard}}\label{Subsec:Lkhard}
We recall \eqref{Def:pt} and Definitions \ref{Def:BasicsLB}, \ref{Def:fkqk}. In order to derive an upper bound on $\sum_{k=0}^{\lfloor\log(n)^7\rfloor}\mathbb{E}_{\mathcal{L}}[\Lambda_k(y)]$ we need an upper bound on $p_n^{-1}\max_{x\ge0}e^{-\vartheta^\ast x}q_k(K_k/\vartheta^\ast-m_{n,h,y}^{\smallsmile}(k)-x)$. For this, we need one more definition.
\begin{definition}\label{Def:pstark}
Set
\begin{align*}
p_{\text{start},k} := \inf_{x\in [\xi_0-1,0]}\mathbb{P}_{\mathcal{L}}\left[\mathcal{B}_{\{0,\dots,k\},m_{n,h}^{\smallsmile}}^{\xi_0,J_x}(T_{\cdot}/\vartheta^\ast)\right].
\end{align*}
\end{definition}

The proof of Lemma \ref{Lem:Lkhard} is split into two further lemmata. The proof of Lemma \ref{Lem:Handlepstart} involves barrier computations and thus is done in Section \ref{Sec:BarrComps}.
\begin{lemma}\label{Lem:UpperBoundqk}
For all $\varepsilon>0$, there is a $C_\varepsilon>0$ such that
\begin{align*}
&\liminf_{n\to \infty} \mathbb{P}\left[\bigcap_{k\le \lfloor\log(n)^7\rfloor}\left\{ p_n^{-1} \max_{x\ge0} e^{-\vartheta^\ast x}q_k(K_k/\vartheta^\ast-m_{n,h,y}^{\smallsmile}(k)-x) \le C_\varepsilon (p_{\text{start},k})^{-1}\right\} \right]\\
&\quad\ge 1-\varepsilon.
\end{align*}
\end{lemma}

\begin{lemma}\label{Lem:Handlepstart}
For all $\varepsilon>0$, there is a  $C_\varepsilon>0$ such that
\[
\liminf_{n\to\infty} \mathbb{P}\left[\sum_{k=0}^{\lfloor\log(n)^7\rfloor} e^{-\vartheta^\ast h_n^{\smallsmile}(k)}(p_{\text{start},k})^{-1}\le C_\varepsilon\right]\ge 1-\varepsilon.
\]
\end{lemma}

\begin{proof}[Proof of Lemma \ref{Lem:Lkhard} assuming Lemma \ref{Lem:UpperBoundqk}, \ref{Lem:Handlepstart} ]
Plugging Lemma \ref{Lem:UpperBoundqk} and Lemma \ref{Lem:Handlepstart} into Proposition \ref{Prop:Lky}, yields that there exists a $C_{\varepsilon/2}>0$ such that
\begin{align}
1-\varepsilon/2 &\le \liminf_{n\to \infty}\mathbb{P}\left[\sum_{k=0}^{\lfloor\log(n)^7\rfloor} \mathbb{E}_{\mathcal{L}}[\Lambda_k(y)] \le C_{\varepsilon/2} e^{\vartheta^\ast y}p_n^{-1}\barrierprobLB\right],\label{eq:sumkqk}
\end{align}
where $e^{-2\vartheta^\ast \xi_0}$ has been pulled into the $C_{\varepsilon/2}$, which is fine, since $\xi_0$ is a (negative) constant.

Applying Proposition \ref{Prop:BarrierUBLB} yields that there are $C_{\varepsilon/2}, c_{\varepsilon/2}>0$ such that
\begin{equation}
1-\varepsilon/2\le \liminf_{n\to \infty}\mathbb{P}\left[ p_n^{-1}\barrierprobLB \le C_{\varepsilon/2} |y|^{c_{\varepsilon/2}}\right] \label{eq:pnbarrLB}.
\end{equation}

Combining \eqref{eq:sumkqk} and \eqref{eq:pnbarrLB} yields the claim of Lemma \ref{Lem:Lkhard}.
\end{proof}

To prepare the proof of Lemma \ref{Lem:UpperBoundqk} we need a statement parallel to Corollary \ref{Cor:ApplBarr1}. The proof is analogous, but using \eqref{eq:BarrSplit2} instead of \eqref{eq:BarrSplit1}. We recall \eqref{Def:LBRTpt} and \eqref{Def:qk}.
\begin{corollary}\label{Cor:ApplBarrSplit2}
Fix $k\le n$. For $x\ge |\xi_0|$ we have that
\begin{align}
q_{k}\left(K_{k}/\vartheta^\ast-m_{n,h,y}^{\smallsmile}(k)-x\right)\le \left(\sum_{j=0}^1 p_n^{\smallsmile,\xi_0-j}(-\lceil x\rceil)\right)\left(p_{\text{start},k}\right)^{-1}. \label{eq:Maxqk}
\end{align}

For $0\le x\le |\xi_0|$, we have that
\begin{align}
q_k(K_{k}/\vartheta^\ast-m_{n,h,y}^{\smallsmile}(k)-x)\le \left(\sum_{j=0}^1 p_n^{\smallsmile,\xi_0-j}(\xi_0)\right)\left(p_{\text{start},k}\right)^{-1}. \label{eq:Maxqk'}
\end{align}
\end{corollary}
\begin{proof}
We begin with the proof of \eqref{eq:Maxqk}. We apply \eqref{eq:BarrSplit2} for $Z_{\cdot} = T_{\cdot}/\vartheta^\ast+m_{n,h}^{\smallsmile}$, $z_1 = -\lceil x\rceil$, $z_2 =\xi_0$, $I = \{0,\dots,n\}$, $t_0 =k$ and $x_0 = -x$ and get that
\begin{align*}
\left(\sum_{j=0}^1 p_n^{\smallsmile,\xi_0-j}( -\lceil x\rceil)\right) &\ge\mathbb{P}_{\mathcal{L}}\left[\mathcal{B}_{\{0,\dots,k\},m_{n,h}^{\smallsmile}}^{ -\lceil x\rceil,J_{x_0}}(T_{\cdot}/\vartheta^\ast)\right]\mathbb{P}_{\mathcal{L}}\left[\mathcal{B}_{\{k,\dots,n\},m_{n,h}^{\smallsmile}-m_{n,h}^{\smallsmile}(k)}^{x_0, J_{\xi_0}}((T_\cdot-T_k)/\vartheta^\ast)\right] \\
&=\mathbb{P}_{\mathcal{L}}\left[\mathcal{B}_{\{0,\dots,k\},m_{n,h}^{\smallsmile}}^{ -\lceil x\rceil,J_{-x}}(T_{\cdot}/\vartheta^\ast)\right]q_k(K_k/\vartheta^\ast-m_{n,h,y}^{\smallsmile}(k)-x)\\
&\ge \mathbb{P}_{\mathcal{L}}\left[\mathcal{B}_{\{0,\dots,k\},m_{n,h}^{\smallsmile}}^{ \xi_0,J_{\xi_0+\lceil x\rceil-x}}(T_{\cdot}/\vartheta^\ast)\right]q_k(K_k/\vartheta^\ast-m_{n,h,y}^{\smallsmile}(k)-x),
\end{align*}
where the last step used that $-\lceil x\rceil \le \xi_0$ and the second equality uses the definition \eqref{Def:qk} of $q_k$ as well as 
\[
m_{n,h,y}^{\smallsmile}(k+\cdot)-m_{n,h,y}^{\smallsmile}(k) = m_{n,h}^{\smallsmile}(k+\cdot)-m_{n,h}^{\smallsmile}(k).
\] Since $\xi_0+\lceil x\rceil-x\in [\xi_0-1,0]$ this allows us to conclude \eqref{eq:Maxqk}.

In the case $x\le |\xi_0|$ we take $z_1 =\xi_0$ instead. The application of \eqref{eq:BarrSplit2} yields that
\[
\left(\sum_{j=0}^1 p_n^{\smallsmile,\xi_0-j}( \xi_0)\right) \ge\mathbb{P}_{\mathcal{L}}\left[\mathcal{B}_{\{0,\dots,k\},m_{n,h}^{\smallsmile}}^{ \xi_0,J_{-x}}(T_{\cdot}/\vartheta^\ast)\right]q_k(K_k/\vartheta^\ast-m_{n,h,y}^{\smallsmile}(k)-x)
\]
and since $-x\in [\xi_0-1,0]$ this allows us to conclude \eqref{eq:Maxqk'}.
\end{proof}

\begin{proof}[Proof of Lemma \ref{Lem:UpperBoundqk}]
Define $d_0 = 0$, $d_1 = -\xi_0$, $d_2 = \log(n)^2$ and $d_3 = \infty$ and
\[
\pi_{j,k} := \max_{x\in [d_{j-1},d_j]} e^{-\vartheta^\ast x}p_n^{-1}q_k(K_k/\vartheta^\ast-m_{n,h,y}^{\smallsmile}(k)-x),\qquad j = 1,2,3,\quad k\le \lfloor\log(n)^7\rfloor.
\] By definition we have $\max_{x\ge 0}e^{-\vartheta^\ast x}p_n^{-1}q_k(K_k/\vartheta^\ast-m_{n,h,y}^{\smallsmile}(k)-x)= \max_{j\le 3}\{\pi_{j,k}\}$, thus it suffices to prove that for $j\in\{1,2,3\}$, $\varepsilon>0$, there is a $C_{\varepsilon/3}>0$ such that 
\begin{equation}
\liminf_{n\to \infty} \mathbb{P}\left[\bigcap_{k\le \lfloor\log(n)^7\rfloor}\left\{\pi_{j,k} \le C_{\varepsilon/3} (p_{\text{start},k})^{-1}\right\} \right]\ge 1-\varepsilon/3. \label{eq:Condpj}
\end{equation}

The case $j = 1$:
By applying \eqref{eq:Maxqk'} we have that $\mathbb{P}$-a.s.\@ for all $k\le \lfloor\log(n)^7\rfloor$
\begin{align}
\pi_{1,k} \le p_n^{-1}\left(\sum_{j=0}^1p_n^{\smallsmile,\xi_0-j}(\xi_0-1)\right)(p_{\text{start},k})^{-1} \label{eq:ApplMaxqk'}
\end{align}
Applying Proposition \ref{Prop:BarrierUBLB} yields that there is a $C_{\varepsilon/3}>0$ such that 
\begin{equation}
1-\varepsilon/3 \le \liminf_{n\to \infty}\mathbb{P}\left[p_n^{-1}\left(\sum_{j=0}^1p_n^{\smallsmile,\xi_0-j}(\xi_0-1)\right) \le  C_{\varepsilon/3}\right] \label{eq:ApplUBLB1}
\end{equation}
equations \eqref{eq:ApplMaxqk'} and \eqref{eq:ApplUBLB1} yield \eqref{eq:Condpj} for $j=1$.

The case $j = 2$:
By applying \eqref{eq:Maxqk} we have that $\mathbb{P}$-a.s.\@ for all $k\le \lfloor\log(n)^7\rfloor$
\begin{align}
\pi_{2,k} \le \max_{-\xi_0\le x\le \log(n)^2} e^{-\vartheta^\ast x}\left(\sum_{j=0}^1p_n^{\smallsmile,\xi_0-j}(-\lceil x\rceil)\right)\left(p_{\text{start},k}\right)^{-1}. \label{eq:ApplMaxqk}
\end{align}
Applying Proposition \ref{Prop:BarrierUBLB} yields that there are $c_{\varepsilon/3}, C_{\varepsilon/3}>0$ such that 
\begin{align}
1-\varepsilon/3 &\le \liminf_{n\to \infty}\mathbb{P}\left[\forall_{x\in [-\log(n),\xi_0]\cap\Z}\, p_n^{-1}\left(\sum_{j=0}^1p_n^{\smallsmile,\xi_0-j}(x)\right)  \le C_{\varepsilon/3} |x|^{c_{\varepsilon/3}}  \right], \label{eq:ApplUBLB2}
\end{align}
and since $\max_{x\ge 0} e^{-\vartheta^\ast x}\lceil x\rceil^{c_{\varepsilon/3}}\le C_{\varepsilon/3}<\infty$ equations \eqref{eq:ApplMaxqk} and \eqref{eq:ApplUBLB2} yield \eqref{eq:Condpj} for $j =2$.

The case $j = 3$:
By bounding $q_k(x)\le 1$ we have that $\pi_{3,k}\le p_n^{-1}e^{-\log(n)^2}$, which $\mathbb{P}$-a.s.\@ converges to 0 by Lemma \ref{Lem:Boundpnastweak}. This gives \eqref{eq:Condpj} for $j=3$.

Thus we have proved \eqref{eq:Condpj} for $j\in\{1,2,3\}$ which finishes the proof.
\end{proof}

\section{Proof of Theorem \ref{Theo:Main2} -- Cutting the Tree}\label{Sec:Finish}
We recall Definitions \ref{Def:Basics}, \ref{Def:Barr} and \ref{Def:pn}. Furthermore, we recall the convention that under $\mathbb{P}_{\mathcal{L}}^l$ we write $S_n = X_{l+1}+\dots+X_{n+l}$ to be shifted in time by $l$, without this being reflected in the notation. We also extend this convention to $T_n$, i.e.\@ under $\mathbb{P}_{\mathcal{L}}^l$ we have $T_n = \vartheta^\ast \sum_{j = l+1}^{n+l} \left(X_j-\kappa_j(\vartheta^\ast)\right)$.
\begin{proof}[Proof of Lemma \ref{Lem:MinBRW}]
Since our jumps are symmetric, we have that $\min_{u:|u| = n} V(u) \stackrel{d}{=} -M_n$. Thus by Lemma \ref{Lem:Ausreißer} we have that $\mathbb{P}$-a.s.
\begin{align*}
\mathbb{P}_{\mathcal{L}}\left[M_n<-2K_n/\vartheta^\ast\right]\le e^{-K_n}.
\end{align*}
It now suffices to prove, that there are $c_\varepsilon'$, $c_{\varepsilon}''$ such that for all $n\in\N$
\begin{align*}
\mathbb{P}[K_n\ge c_{\varepsilon}'' n] &\ge 1-\varepsilon,\\
\mathbb{P}[K_n\le c_{\varepsilon}'n] &\ge 1-\varepsilon.
\end{align*}

This is true since $K_n/n$ converges $\mathbb{P}$-a.s.\@ because the $\kappa_k(\vartheta^\ast)$ are i.i.d.\@ with finite variance.
\end{proof}

\begin{definition}
For $n,l\in\N$, set
\begin{align*}
K_{n}^{(l)} &:= \sum_{j=l+1}^{n+l} \kappa_j(\vartheta^\ast),\\
p_n^{(l)} &:= \mathbb{P}_{\mathcal{L}}^{l}\left[\forall_{s\le n}\, \xi_0+T_s/\vartheta^\ast\le 0, \xi_0+T_n/\vartheta^\ast\in [\xi_0-1,\xi_0]\right],\\
m_n^{(l)} &:= K_{n}^{(l)}/\vartheta^\ast+\log(p_{n}^{(l)})/\vartheta^\ast.
\end{align*}
These are the by $l$ time-shifted analogues of $K_n$, $p_n$ and $m_n$.
\end{definition}
We note that the following version of Theorem \ref{Theo:UBRT} holds, since we only make a statement about the distribution of $\mathbb{P}_{\mathcal{L}}[\dots]$, which is invariant under time-shifts
\begin{corollary}\label{Kor:TheoTimeshift}
For all $\varepsilon>0$, there are $c_\varepsilon>0$, $C_{\varepsilon} \in (0,1]$ such that for $y\le \xi_0$,
\[
\liminf_{n\to \infty} \mathbb{P}\left[\mathbb{P}_{\mathcal{L}}^l[M_{n-l}-m_{n-l}^{(l)}\ge \xi_0-y-1] \ge e^{\vartheta^\ast y} C_\varepsilon |y|^{-c_\varepsilon}\right]\ge 1-\varepsilon.
\]
\end{corollary}
\begin{lemma}\label{Lem:EndPrelim1}
For all $\varepsilon>0$, there are $c_\varepsilon$, $C_\varepsilon$, $c_\varepsilon'$, $c_\varepsilon''>0$ such that for $l\in\N$ and $y\in [-c_\varepsilon' l+1,0]$,
\begin{align*}
&\liminf_{n\to \infty} \mathbb{P}\left[\mathbb{P}_{\mathcal{L}}[M_n-m_{n-l}^{(l)} \le y]\le e^{-c_{\varepsilon}''l}+\left(1-C_\varepsilon|\xi_0-c_{\varepsilon}'l-y-1|^{-c_\varepsilon}e^{\vartheta^\ast (\xi_0-c_{\varepsilon}'l-y-1))}\right)^{2l}\right]\\
&\ge 1-\varepsilon
\end{align*}
\end{lemma}
\begin{proof}
We have that $\mathbb{P}$-a.s.
\begin{align*}
&\mathbb{P}_{\mathcal{L}}\left[M_n-m_{n-l}^{(l)} \le y\right]\\
&\le \mathbb{P}_{\mathcal{L}}[M_l\le -c_{\varepsilon}' l]+\mathbb{P}_{\mathcal{L}}\left[M_l\ge -c_{\varepsilon}'l, \forall_{v: |v|= l}\, \max_{u: |u| = n, u_l = v } V(u)-m_{n-l}^{(l)} \le c_\varepsilon'l+y\right]\\
&\le  \mathbb{P}_{\mathcal{L}}[M_l\le -c_{\varepsilon}' l]+\mathbb{P}_{\mathcal{L}}^l\left[M_{n-l}-m_{n-l}^{(l)} \le c_\varepsilon' l+y\right]^{2l}\\
&=\mathbb{P}_{\mathcal{L}}[M_l\le -c_{\varepsilon}' l]+\left(1-\mathbb{P}_{\mathcal{L}}^l\left[M_{n-l}-m_{n-l}^{(l)} \ge c_\varepsilon' l+y\right]\right)^{2l}.
\end{align*}
Now applying Lemma \ref{Lem:MinBRW} and Corollary \ref{Kor:TheoTimeshift} yields the claim.
\end{proof}

\begin{lemma}\label{Lem:Prelim2}
For all $\varepsilon>0$, there are $l_\varepsilon\in\N$, $y_\varepsilon<0$ such that 
\[
\liminf_{n\to \infty}\mathbb{P}\left[\mathbb{P}_{\mathcal{L}}[M_n-m_{n-l_\varepsilon}^{(l_\varepsilon)} \le y_\varepsilon] \le \varepsilon\right]\ge 1-\varepsilon.
\]
\end{lemma}
\begin{proof}
Choosing $y = -c_\varepsilon' l+1$ in Lemma \ref{Lem:EndPrelim1} yields that
\begin{align*}
1-\varepsilon&\le \liminf_{n\to \infty}\mathbb{P}\Bigg[\mathbb{P}_{\mathcal{L}}\left[M_n-m_{n-l}^{(l)}  \le -c_{\varepsilon}' l+1\right] \le e^{-c_{\varepsilon}'' l}+\left(1-C_\varepsilon |\xi_0-2|^{-c_\varepsilon}e^{\vartheta^\ast (\xi_0-2)}\right)^{2l}\Bigg].
\end{align*}
Since $0<C_{\varepsilon} \le1$ we have $\left(1-C_\varepsilon |\xi_0-2|^{-c_\varepsilon}e^{\vartheta^\ast (\xi_0-2)}\right) \in [0,1)$. Thus, we can choose $l$ such that $e^{-c_{\varepsilon}'' l}+\left(1-C_\varepsilon |\xi_0-2|^{-c_\varepsilon}e^{\vartheta^\ast (\xi_0-2)}\right)^{2l}\le \varepsilon$, which yields the claim.
\end{proof}

\begin{lemma}\label{Lem:ContrOffs}
For all $\varepsilon>0$, $l\in\N_{\ge2}$ there is a $C_{\varepsilon,l}>0$ such that,
\begin{align}
\liminf_{n\to \infty}\mathbb{P}\left[ \left|K_n-K_{n-l}^{(l)}\right|\le C_{\varepsilon,l}\right]&\ge 1-\varepsilon , \label{eq:Kntimeshift}\\
\liminf_{n\to \infty}\mathbb{P}\left[\log(p_n)-\log(p_{n-l}^{(l)}) \le C_{\varepsilon,l}\right] &\ge 1-\varepsilon. \label{eq:pntimeshift}
\end{align}
\end{lemma}
Furthermore, for all $\varepsilon>0$, $l\in\N_{\ge2}$ there is a $C_{\varepsilon,l}>0$
\begin{align*}
\liminf_{n\to \infty} \mathbb{P}\left[m_n-m_{n-l}^{(l)}\le C_{\varepsilon,l}\right]\ge 1-\varepsilon
\end{align*}

\begin{proof}[Proof of \eqref{eq:Kntimeshift}]
The first statement is immediate since by definition $K_n-K_{n-l}^{(l)} = K_l$ and $K_l/l\to 0$ $\mathbb{P}$-a.s.

For the proof of the second statement we define for $n\in\N$, $l\le n$, $x\le 0$
\begin{align*}
p_{n-l}^{(l)}(x) &:= \mathbb{P}_{\mathcal{L}}^l\left[\forall_{s\le n-l} x+T_s/\vartheta^\ast\le0, \xi_0+T_n/\vartheta^\ast\in [\xi_0-1,\xi_0]\right],\\
\widetilde{p}_{\text{start},l} &:= \inf_{x\in[\xi_0-1,0]} \mathbb{P}\left[\mathcal{B}_{[0,l]}^{\xi_0,J_x}(T_\cdot/\vartheta^\ast) \right].
\end{align*}
As in Corollary \ref{Cor:ApplBarrSplit2} applying \eqref{eq:BarrSplit2} yields that for $\xi_0\le x\le 0$ 
\begin{equation}
\begin{aligned}
p_{n-l}^{(l)}(x) &\le \mathbb{P}_{\mathcal{L}}\left[\mathcal{B}_{[0,n]}^{\xi_0, [\xi_0-2,\xi_0]}(T_\cdot/\vartheta^\ast)\right](\widetilde{p}_{\text{start},l})^{-1},\quad \text{for}\ \xi_0 \le x\le 0.
\end{aligned}\label{eq:TimeShiftToSpaceShift}
\end{equation}
We use $g_{\mu,\sigma^2}(x)$ to denote the Gaussian density with mean $\mu\in\R$ and variance $\sigma^2\ge 0$. We use the Markov property at time $l$ to get that
\begin{equation}\label{eq:pnUBFor7}
\begin{aligned}
p_n &= \int_{-\infty}^{0} g_{\xi_0-W_l,l}(z) p_{n-l}^{(l)}(z)\;\mathrm{d}z \\
&\stackrel{\eqref{eq:TimeShiftToSpaceShift}}{\le} \int_{-\infty}^{\xi_0} g_{\xi_0-W_l,l}(x) p_{n-l}^{(l)}(x)\;\mathrm{d}x+(\widetilde{p}_{\text{start},l})^{-1} \int_{\xi_0}^0 g_{\xi_0-W_l,l}(x) \mathbb{P}_{\mathcal{L}}\left[\mathcal{B}_{[0,n]}^{\xi_0,[\xi_0-2,\xi_0]}(T_{\cdot}/\vartheta^\ast)\right] \;\mathrm{d}x.
\end{aligned}
\end{equation}
Similarly to Proposition \ref{Prop:BarrComp} we see that\footnote{Technically we use a time-shifted version of the result. For these the constants bounded in Section \ref{Sec:Constants} will depend on $(W_{n+l}-W_l)_{n\in\N}$. For our application, where we fix $l$, the constants will still be bounded with high probability.  We note however that we can't bound the constants for all $l$ simultanously.} for all $\varepsilon>0$ there are $c_{\varepsilon,l}, C_{\varepsilon,l}, C_\varepsilon>0$ such that 
\begin{align}
1-\varepsilon/3&\le \liminf_{n\to \infty}\mathbb{P}\left[\bigcap_{x\in [-\log(n)^2,\xi_0]\cap \Z} \left\{\frac{ \inf_{z \in [x-1,x]} p_{n-l}^{(l)}(x)}{p_{n-l}^{(l)}} \le C_{\varepsilon,l} |x|^{c_{\varepsilon,l}}\right\}\right], \label{eq:StartMoveTimeShift}\\
1-\varepsilon/3&\le \liminf_{n\to\infty} \mathbb{P}\left[\left\{ \frac{\mathbb{P}_{\mathcal{L}}\left[\mathcal{B}_{[0,n]}^{\xi_0, [\xi_0-2,\xi_0]}(T_{\cdot}/\vartheta^\ast)\right]}{p_n}\right\}\le C_\varepsilon \right]. \label{eq:EndBlowUp}
\end{align}

Furthermore,  by Lemma \ref{Lem:Boundpnastweak} there is\footnote{Actually we have shifted the starting point by one compared to the precise statement of Lemma \ref{Lem:Boundpnastweak}, this will not make a difference in the proof of the lemma we give in Section \ref{Sec:Constants}.} is a $C>0$ such that $\mathbb{P}$-a.s.\@ we have 
\begin{equation} \label{eq:plstartbig}
\limsup_{l\to \infty} \frac{|\widetilde{p}_{\text{start},l}|}{\log(l)} \le C.
\end{equation}
Finally, since $g_{\xi_0-W_l,l}(x)\le \frac{1}{\sqrt{2\pi l}}$ for all $x\in\R$ we can follow from \eqref{eq:EndBlowUp} and \eqref{eq:plstartbig} that 
\[
1-(2/3)\varepsilon\le \liminf_{n\to\infty}\mathbb{P}\left[ (\widetilde{p}_{\text{start},l})^{-1} \int_{\xi_0}^0 g_{\xi_0-W_l,l}(x) \mathbb{P}_{\mathcal{L}}\left[\mathcal{B}_{[0,n]}^{\xi_0,[\xi_0-2,\xi_0]}(T_{\cdot}/\vartheta^\ast)\right] \;\mathrm{d}x \le p_n/2\right].
\]
Combining the last display with \eqref{eq:pnUBFor7} and \eqref{eq:StartMoveTimeShift} while using Gaussian tailbounds and the fact that $\int_{-\infty}^{\xi_0} g_{\xi_0-W_l,l}(x)C|x|^c\;\mathrm{d}x<\infty$ for all $c,C$ yields that there exists $C_{\varepsilon,l}>0$ such that
\begin{equation} \label{eq:pnpntimeshiftFast}
1-\varepsilon \le \liminf_{n\to\infty}\mathbb{P}\left[ p_n\le C_{\varepsilon,l} p_{n-l}^{(l)}+e^{-\frac{(\log(n)^2+W_l-\xi_0)^2}{2}} \right].
\end{equation}
Another application of Lemma \ref{Lem:Boundpnastweak}, yields that $\mathbb{P}$-a.s.\@ we have $\lim\limits_{n\to\infty} p_n^{-1}e^{-\frac{(\log(n)^2-W_l+\xi_0)^2}{2}} = 0$, which together with \eqref{eq:pnpntimeshiftFast} implies that there exists $C_{\varepsilon,l}>0$ such that 
\[
1-\varepsilon \le \liminf_{n\to\infty}\mathbb{P}\left[ p_n\le C_{\varepsilon,l} p_{n-l}^{(l)}\right].
\]
Rearranging the event in the last display yields \eqref{eq:pntimeshift}.

The third claim follows directly from the first two.
\end{proof}

\begin{proof}[Proof of Theorem \ref{Theo:Main2}]
By Lemma \ref{Lem:Prelim2} there are $l_\varepsilon\in\N$, $z_\varepsilon<0$ such that
\begin{align}
1-\varepsilon/2&\le \liminf_{n\to \infty}\mathbb{P}\left[\mathbb{P}_{\mathcal{L}}[M_n-m_{n-l_\varepsilon}^{(l_\varepsilon)} \le z_\varepsilon]\le \varepsilon\right].\label{eq:FinPrf1}
\end{align}
Since for all $y<0$
\[
\mathbb{P}_{\mathcal{L}}[M_n-m_n \le y] =\mathbb{P}_{\mathcal{L}}[M_n-m_{n-l_\varepsilon}^{(l_\varepsilon)} \le m_n-m_{n-l_\varepsilon}^{(l_\varepsilon)}+y],
\]
Lemma \ref{Lem:ContrOffs} implies that there is a $C_{\varepsilon,l_{\varepsilon}}>0$ such that for all $y<0$
\begin{equation}
1-\varepsilon/2 \le \liminf_{n\to\infty}\mathbb{P}\left[\mathbb{P}_{\mathcal{L}}[M_n-m_n \le y] \le \mathbb{P}_{\mathcal{L}}[M_n-m_{n-l_\varepsilon}^{(l_\varepsilon)} \le y+C_{\varepsilon,l_{\varepsilon}}]\right]. \label{eq:FinPrf2}
\end{equation}
Taking $y_\varepsilon := z_\varepsilon-C_{\varepsilon,l_{\varepsilon}}$ in \eqref{eq:FinPrf2} and combining it with \eqref{eq:FinPrf1} yields the claim of the theorem.
\end{proof}

\section{Shift of the Starting Point -- Preparation for Section \ref{Sec:BarrComps}}\label{Sec:MoveStart}
In this section, let $B$ be a Brownian motion starting at $x$ with regard to $\mathbb{P}_x$. For $t\ge 0$ let $W^{(t)}: [0,t]\to \R$ be continuous with $W_0 = 0$ and let $h_t:[0,t]\to \R$ be a continuous function with $h_t(0) = h_t(t) = 0$. We will use $W$ and $h$ as a shorthand for $(W^{(t)})_{t\ge0}$, $(h_t)_{t\ge0}$, we will suppress the upper index of $W^{(t)}_s$ and instead write $W_s$. Additional assumptions on $h$ and $W$ are stated in Section \ref{Subsec:DefSum} together with further definitions, in Section \ref{Sec:Constants} we check that the assumptions hold for the specific $W$ we use in the first half of the paper. In Section \ref{Sec:ResMoveStart} we will state the results of Section \ref{Sec:MoveStart}. Sections \ref{Sec:NebOS}, \ref{Sec:MainThLB}, \ref{Sec:MainOS} contain the proofs of Theorems \ref{Theo:NebOS}, \ref{Theo:MainThLB} and \ref{Sec:MainOS} respectively.

\subsection{Definitions and Assumptions}\label{Subsec:DefSum}
We recall Definition \ref{Def:Barr}.
\begin{definition}\label{Def:pthAllg}
For $y_0, y,x \le0$, define
\begin{align*}
p_{t,h}^{(y_0)}(y) &:= \mathbb{P}_{y}\left[\mathcal{B}_{[0,t],h_t(\cdot)-W_\cdot}^{J_{y_0}}(B_\cdot)\right],\\
p_{t,h}^{(y_0)}(y,x) &:= \mathbb{P}_{y}\left[\mathcal{B}_{[x,t],h_t(\cdot)-W_\cdot}^{J_{y_0}}(B_\cdot)\right].
\end{align*}
For $\mu\in \R$, $\sigma^2>0$ write
\begin{equation}
g_{\mu,\sigma^2}(x) := \frac{1}{\sqrt{2\pi \sigma^2}}e^{-\frac{(x-\mu)^2}{2\sigma^2}} \label{Def:Gausz}
\end{equation}
for the Gaussian density function.
\end{definition}

Define 
\begin{align}
1\vee\sup_{s\in [e,t]}\frac{|W_s|}{\sqrt{s}\cdot \sqrt{\log(s)}} &=: \mathcal{C}_{\log}(W,t) =: \ClogW ,\label{Ass:WBMlog}\\
1\vee\sup_{t\in \R}\sup_{s\in [0,t]} \frac{\max\{|h_t(s)|, |h_t(t-s)|\}}{\sqrt{1+s}} &=: C_1(h) := \Coh , \label{Ass:WachsH}
\end{align}
and assume that both quantities are finite. Furthermore, set
\begin{align}
\ExpPos := \gamma(W,h,t) &:= \inf\left\{\gamma\ge 0 : \forall_{2\le s\le t/2}\,  \mathbb{P}\left[\mathcal{B}_{[0,s], h_t(\cdot)-W_\cdot}^{-1,[-4(\Coh-1) \sqrt{s}, -2\Coh\sqrt{s}]}(B_\cdot)\right] \ge s^{-\gamma}\right\}\label{Ass:QuenchedWall}.
\end{align}
We assume that $\gamma>0$.  Furthermore, define
\begin{align}
e^{-128-16\Coh ^2-80\ClogW ^2-134\Coh -96\ClogW -32 \Coh \ClogW }&=: C_2(W,h,t) =:\ConstWnW. \label{Def:WeirdC1}
\end{align}
The constants in the definition of $\ConstWnW$ aren't important, what matters is, that $\ConstWnW$ is a continuous function of $\ClogW $ and $\Coh $.

Finally, for $\lambda>0$ arbitrary but fixed, set
\begin{align}
c_1 &:= c_1(\lambda) := \sqrt{8(\lambda+1)} \label{eq:c3},\\
\Cnull &:= C_3(W,h,\lambda,t) := \max\{48, 128\ClogW ^2, 64\sqrt{3}\Coh^2 ,\left(32\ExpPos+128\right)^4, \sqrt{5}c_1 \}. \label{Ass:AssC0W}
\end{align}

\subsection{Summary of Results} \label{Sec:ResMoveStart}
\begin{theorem}\label{Theo:NebOS}
Let $-2t^{\frac{1}{4}}\le y\le x\le -e$, $y_0\le 0$ and $t\ge 64$. We have that
\begin{align}
\frac{p_{t,h}^{(y_0)}(x,y^2)}{p_{t,h}^{(y_0)}(x)}\in \left[1,\ConstWnW^{-1}|y|^{2\ExpPos}\right] \label{eq: OS1}.
\end{align}
\end{theorem}

\begin{theorem}\label{Theo:MainThLB}
For $t>0$, $y \in [-2\log(t)^2, -e]$, $x\in [y,0]$,  $y_0\le 0$, we have that
\begin{align*}
\frac{p_{t,h}^{(y_0)}(y,y^2)}{p_{t,h}^{(y_0)}(x,y^2)}\ge ce^{-2\Coh }e^{-\ClogW \sqrt{\log(y^2)}}.
\end{align*}
\end{theorem}

Theorem \ref{Theo:MainThLB} can be proved by comparing the density of a Brownian motion at time $y^2$ when starting at $y$ to the same starting at $x$. We point out that we will only apply Theorem \ref{Theo:MainThLB} for positive $h$. In this case, the term $e^{-2\Coh }$ in Theorem \ref{Theo:MainThLB} is not necessary.
\begin{theorem}\label{Theo: MainOS}
For $ y \in[-2\log(t)^2, -e]$, $x\in [y,-e]$, $y_0\le 0$ and $\lambda>0$ there is a $t_0(\lambda)>0$ such that for $t\ge t_0(\lambda)$,
\begin{align}
\frac{p_{t,h}^{(y_0)}(y,y^2)}{p_{t,h}^{(y_0)}(x,y^2)}&\le \ConstWnW^{-1}4^{\ExpPos+2}e^{4\Coh }|y|^{2\ExpPos+3\Cnull }+Ce^{-\lambda \log(t)}p_{t,h}^{(y_0)}(x)^{-1}. \label{eq: OS3}
\end{align}
\end{theorem}

The proof of Theorem \ref{Theo: MainOS} takes up most of the section, an overview is given in Section \ref{subsec:DefOutl}. We also want to point out that the interval in which the Brownian motion is at time $t$ being of length 1 (or even $O(1)$) is not necessary for Theorem \ref{Theo: MainOS} to hold.

For these theorems to be useful, it is necessary to control the constants defined in Section \ref{Subsec:DefSum}, for the specific choices of $W$ and $h$ relevant in the first half of the paper, we will do this in Section \ref{Sec:Constants}.
\subsection{Proof of Theorem \ref{Theo:NebOS}} \label{Sec:NebOS}

We proceed with the upper bound in \eqref{eq: OS1}. We recall Definition \ref{Def:Barr}.
In this section, for $0\le t_1<t_2<\infty$ we set
\[
A_{t_1,t_2} := \mathcal{B}_{[t_1,t_2],h_t(\cdot)-W_\cdot}(B_\cdot).
\]

We first condition on the endpoint. For $x,y_0\le 0$, we have that  $B_t-W_t\sim \mathcal{N}(x-W_t, t)$ under $\mathbb{P}_x$ and thus, see \eqref{Def:Gausz},
\begin{align}
p_{t,h}^{(y_0)}(x) &= \int_{y_0-1}^{y_0} g_{x-W_t,t}(z)\mathbb{P}_{x}\left[ A_{0,t} | B_t- W_t = z\right]\;\mathrm{d}z. \label{eq: BBridgeSetup}
\end{align}

Under $\mathbb{P}_{x}\left[\;\cdot\;|B_t- W_t = z\right]$ the process $(B_s)_{s\in [0,t]}$ is a Brownian bridge from $x$ to $z+W_t$. We use that to prove that the marginals of $(B_s)_{s\in [0,t]}$ are associated w.r.t.\@ that measure.
\begin{lemma}\label{Lem: Association}
Let $(X_s)_{s\in [0,t]}$ be a Brownian bridge from $x$ to $y$, $x,y\in\R$. Let $t_1<\dots< t_n\in [0,t]$  and $A,B$ be decreasing in $(X_{t_1},\dots, X_{t_n})$. Then
\begin{align*}
\mathbb{P}[A\cap B]\ge \mathbb{P}[A]\mathbb{P}[B].
\end{align*}
\end{lemma}

\begin{proof}
Let $(Y_t)_{t\ge0}$ be a Brownian motion starting at $x$, we have that $(X_s)_{s\in[0,t]} \stackrel{d}{=}\left( Y_t-\frac{s}{t}(Y_t-y)\right)_{s\in[0,t]}$, in particular $(X_{t_1},\dots, X_{t_n})$ are normally distributed and a straightforward calculation yields that for $j\le k$ 
\begin{align}
\Cov[X_{t_j},X_{t_k}] &= \Cov\left[Y_{t_j}-\frac{t_j}{t}Y_{t}-\frac{t_j}{t}y, Y_{t_k}-\frac{t_k}{t}Y_{t}-\frac{t_k}{t}y\right]= t_j\left(\frac{t-t_k}{t}\right)\ge 0. \label{eq:Covpos}
\end{align}
 By \cite{PaperAssoc}, equation \eqref{eq:Covpos} implies that $(X_{t_1},\dots, X_{t_n})$ are associated. Since $\mathbf{1}_A$, $\mathbf{1}_B$ are decreasing in $(X_{t_1},\dots, X_{t_n})$, this gives that $\Cov[\mathbf{1}_A,\mathbf{1}_B]\ge 0$.
\end{proof}
Lemma \ref{Lem: Association} and \eqref{eq: BBridgeSetup} imply that for $t\ge 0$,  $0\ge x\ge y$ and $y_0\le 0$,
\begin{equation}
p_{t,h}^{(y_0)}(x)\ge \int_{y_0-1}^{y_0} g_{x-W_t,t}(z)\mathbb{P}_x\left[A_{0,y^2}|B_t-W_t = z\right]\mathbb{P}_x\left[A_{y^2,t}|B_t-W_t = z\right]\;\mathrm{d}z \label{eq:Decomp}
\end{equation}

We want to pull the factor $\mathbb{P}_x\left[A_{0,y^2}|B_t-W_t = z\right]$ in \eqref{eq:Decomp} out of the integral and for this need to bound it uniformly for $z\in [y_0-1,y_0]$. This is done in Lemma \ref{Lem: FKGHilf2} whose proof is postponed. The main idea is, that conditioning on $B_t-W_t = z$ does barely change the distribution of $(B_s)_{s\le y^2}$ for $y$ and $z$ small enough. 
\begin{lemma}\label{Lem: FKGHilf2}
For $t\ge 64$, $0\ge x \ge y$,  $-2t^{\frac{1}{4}} \le y\le -e$, $z\le 0$,  we have that
\begin{align*}
\mathbb{P}_{x}\left[A_{0,y^2} |  B_t-W_t = z\right]&\ge \ConstWnW\mathbb{P}_{x}\left[A_{0,y^2}, B_{y^2}-W_{y^2}\in [4\Coh y,2\Coh y]\right].
\end{align*}
\end{lemma}
The bound $t\ge 64$ is of no further significance and is just a more concrete way to say, that we need $t$ big enough, but importantly not depending on $W$.

Equation \eqref{eq:Decomp} and Lemma \ref{Lem: FKGHilf2} directly imply the following corollary.
\begin{corollary} \label{Cor: Prep}
For $t\ge 64$, $0\ge x\ge y\ge -2t^\frac{1}{4}$, $y_0\le 0$, we have that
\[
p_{t,h}^{(y_0)}(x)\ge \ConstWnW\mathbb{P}_{x}\left[A_{0,y^2}, B_{y^2}-W_{y^2}\in [4\Coh y,2\Coh y]\right]p_{t,h}^{(y_0)}(x,y^2).
\]
\end{corollary}
We can control $\mathbb{P}_{x}\left[A_{0,y^2}, B_{y^2}-W_{y^2}\in [4\Coh y,2\Coh y]\right]$ for $y\le x\le -1$ by applying \eqref{Ass:QuenchedWall}. This is done in the next lemma.
\begin{lemma} \label{Lem: Concl}
We have that for $-2t^{\frac{1}{4}}\le y\le x\le -1$ and $t\ge 64$,
\[
\mathbb{P}_{x}\left[A_{0,y^2}, B_{y^2}-W_{y^2}\in [4\Coh y,2\Coh y]\right] \ge |y|^{-2\ExpPos}.
\] 
\end{lemma}
\begin{proof}
We have that
\begin{align*}
&\mathbb{P}_{x}\left[A_{0,y^2}, B_{y^2}-W_{y^2}\in [4\Coh y,2\Coh y]\right]\\
&\stackrel{y\le x\le -1}{\ge} \mathbb{P}_{x}\left[\forall_{s\le y^2}\, B_s+h_t(s)-(1+x)\le W_s, B_{y^2}-(1+x)-W_{y^2}\in [(4\Coh -1)y,2\Coh y]\right]\\
&= \mathbb{P}_{-1}\left[A_{0,y^2}, B_{y^2}-W_{y^2}\in [(4\Coh -1)y,2\Coh y]\right]\stackrel{\eqref{Ass:QuenchedWall}}{\ge}|y|^{-2\ExpPos},
\end{align*}
where we note that the application of \eqref{Ass:QuenchedWall} uses that $y^2\le 4t^{1/2}\le t/2$ for $t\ge 64$.
\end{proof}
\begin{proof}[Proof of Theorem \ref{Theo:NebOS} assuming Lemma \ref{Lem: FKGHilf2}]
The lower bound in \eqref{eq: OS1} is immediate by monotonicity.

Combining Corollary \ref{Cor: Prep} and Lemma \ref{Lem: Concl} yields the upper bound in \eqref{eq: OS1}.
\end{proof}

\begin{proof}[Proof of Lemma \ref{Lem: FKGHilf2}]
First, we note that the event $A_{0,y^2}$ is decreasing in $B_\cdot$. In particular, it is monotonous in the starting and endpoint of a Brownian bridge, such that for all $z\le 0$ we have
\begin{equation}
\mathbb{P}_x[A_{0,y^2}|B_t-W_t = z] \ge \mathbb{P}_x[A_{0,y^2}| B_t-W_t = 0]. \label{eq:MonBBr}
\end{equation} 
Next, we define
\begin{alignat*}{2}
\mu_1 &= x+y^2/t(W_t-x)-W_{y^2},&\quad \mu_2 &= x-W_{y^2},\\
\sigma_1^2 &= y^2\frac{t-y^2}{t},&\quad \sigma_2^2 &= y^2,
\end{alignat*}
and recall \eqref{Def:Gausz}. The relevance of these is that   the random variable $B_{y^2}-W_{y^2}$ has the density function $g_{\mu_2,\sigma_2^2}$  under $\mathbb{P}_x$ and the density function $g_{\mu_1,\sigma_1^2}$ under $\mathbb{P}_x\left[\;\cdot\;|B_t-W_t = 0\right]$. The latter implies that
\begin{align}
\mathbb{P}_{x}\left[ A_{0,y^2} |B_t- W_t = 0\right] &= \int_{-\infty}^{-h_t(y^2)} g_{\mu_1,\sigma_1^2}(r)\mathbb{P}_{x}\left[A_{0,y^2}| B_t- W_t = 0, B_{y^2}-W_{y^2} = r\right]\notag\\
&= \int_{-\infty}^{-h_t(y^2)} g_{\mu_1,\sigma_1^2}(r)\mathbb{P}_{x}\left[A_{0,y^2}|  B_{y^2}-W_{y^2} = r\right], \notag\\
&\stackrel{\eqref{Ass:WachsH}}{\ge} \int_{4\Coh y}^{2\Coh y} \frac{g_{\mu_1,\sigma_1^2}(r)}{g_{\mu_2,\sigma_2^2}(r)}g_{\mu_2,\sigma_2^2}(r)\mathbb{P}_x\left[A_{0,y^2}|B_{y^2}-W_{y^2} = r\right]\mathrm{d}r.\label{eq: MidEnd}
\end{align}

 We have that
\begin{align}
\frac{g_{\mu_1,\sigma_1^2}(r)}{g_{\mu_2,\sigma_2^2}(r)} &= \sqrt{\frac{t}{t-y^2}} e^{\frac{(r-\mu_2)^2}{\sigma_2^2}-\frac{(r-\mu_1)^2}{\sigma_1^2}}\ge e^{\frac{(r-\mu_2)^2}{\sigma_2^2}-\frac{(r-\mu_1)^2}{\sigma_1^2}}.\label{eq: DensFKG}
\end{align}
Direct calculation using that $-2t^{1/4}\le y\le x\le -e$ and \eqref{Ass:WBMlog} yields that for $r\in [4\Coh y,2\Coh y]$
\begin{align}
e^{\frac{(r-\mu_2)^2}{\sigma_2^2}-\frac{(r-\mu_1)^2}{\sigma_1^2}}&\ge \ConstWnW \label{eq: Ineq1}
\end{align}

Plugging \eqref{eq: DensFKG} and the inequality \eqref{eq: Ineq1} into \eqref{eq: MidEnd} yields that 
\begin{align}
\mathbb{P}_{x}\left[A_{0,y^2} | B_t- W_t = 0\right] &\ge  \ConstWnW\int_{4\Coh y}^{2\Coh y}g_{\mu_2,\sigma_2^2}(r)\mathbb{P}_x\left[A_{0,y^2}|B_{y^2}-W_{y^2} = r\right]\mathrm{d}r \notag\\
&= \ConstWnW\mathbb{P}_{x}\left[A_{0,y^2}, B_{y^2}- W_{y^2}\in [4\Coh y,2\Coh y]\right], \label{eq:BBrBoundLS}
\end{align}
where the last step uses that under  $\mathbb{P}_x$ the random variable $B_{y^2}-W_{y^2}$ has the density $g_{\mu_2,\sigma_2^2}(r)$. Combining \eqref{eq:MonBBr} with \eqref{eq:BBrBoundLS} yields the claim of the lemma.
\end{proof}
\subsection{Proof of Theorem \ref{Theo:MainThLB} } \label{Sec:MainThLB}

For $z\le 0$ and $r\in \R$, define
\begin{align*}
\mu_{z,y^2,h} &:= z+h_t(y^2)-W_{y^2}, \\
\Delta_{y^2,W,h}(s) &:= h_t(s+y^2)-h_t(y^2)-(W_{s+y^2}-W_{y^2}),\\
p_{t;y^2,h}^{(y_0)}(r) &:= \mathbb{P}_{r}\Big[\mathcal{B}_{[0,t-y^2], \Delta_{y^2,W,h}}^{J_{y_0}} \Big].
\end{align*}
We recall \eqref{Def:Gausz}. Using the Markov property we have that
\begin{equation}
p_{t,h}^{(y_0)}(z,y^2) = \int_{-\infty}^0 g_{\mu_{z,y^2,h},y^2}(r)p_{t;y^2,h}^{(y_0)}(r)\;\mathrm{d}r. \label{eq: SeedRec}
\end{equation}

We finish the proof by taking $z = y$ in \eqref{eq: SeedRec} and bounding $g_{\mu_{y,y^2,h},y^2}(r)g_{\mu_{x,y^2,h},y^2}(r)^{-1}$. A direct computation yields that
\begin{equation}
g_{\mu_{y,y^2,h},y^2}(r)g_{\mu_{x,y^2,h},y^2}(r)^{-1} = \exp\left(\frac{r(y-x)}{y^2}+\frac{x^2-y^2}{2y^2}+\frac{( W_{y^2}-h_t(y^2))(y-x)}{y^2} \right). \label{eq: LBQuotDens}
\end{equation}
Applying \eqref{eq: SeedRec} for $z = y$ yields that
\pushQED{\qed} 
\begin{align*}
p_{t,h}^{(y_0)}(y,y^2) &= \int_{-\infty}^0 g_{\mu_{y,y^2,h},y^2}(r)p_{t;y^2,h}^{(y_0)}(r)\;\mathrm{d}r\\
&\stackrel{\eqref{eq: LBQuotDens}}{=} e^{\frac{x^2-y^2}{2y^2}+\frac{(y-x) W_{y^2}}{y^2}-\frac{h_t(y^2)(y-x)}{y^2}}\int_{-\infty}^0 e^{\frac{r(y-x)}{y^2}}g_{\mu_{y,y^2,h},y^2}(r)p_{t;y^2,h}^{(y_0)}(r)\;\mathrm{d}r\\
&\stackrel{y-x\le 0}{\ge} e^{\frac{x^2-y^2}{2y^2}+\frac{(y-x) W_{y^2}}{y^2}-\frac{h_t(y^2)(y-x)}{y^2}}\int_{-\infty}^0 g_{\mu_{y,y^2,h},y^2}(r)p_{t;y^2,h}^{(y_0)}(r)\;\mathrm{d}r\\
&\stackrel{\eqref{eq: SeedRec}}{=}  e^{\frac{x^2-y^2}{2y^2}+\frac{(y-x) W_{y^2}}{y^2}-\frac{h_t(y^2)(y-x)}{y^2}}p_{t,h}^{(y_0)}(x,y^2)\\
&\stackrel{|y|\ge e, \eqref{Ass:WBMlog}, \eqref{Ass:WachsH}}{\ge} ce^{-2\Coh }e^{-\ClogW \sqrt{\log(y^2)}}p_{t,h}^{(y_0)}(x,y^2). \qedhere
\end{align*}
\popQED
\subsection{Proof of Theorem \ref{Theo: MainOS}} \label{Sec:MainOS}
\subsubsection{Definitions and Outline}\label{subsec:DefOutl}

As in Section \ref{Sec:MainThLB} we will use \eqref{eq: SeedRec}, but the factor $e^{\frac{r(1+y)}{y^2}}$ is not a-priori controllable. To fix this we will split the domain of integration into $[-L,0]$, $[-c_1(t), -L]$, $(-\infty,-c_1(t)]$. On the first and the last interval, we will be able to control the integrand. On the middle region, we will use an analogue to \eqref{eq: SeedRec}, which is stated in \eqref{eq: Iter}. This will give us a double integral over $[-c_1(t),-L]\times (-\infty,0]$, we again split up the domain of integration into three parts and iterate the process. The $L$ in the $l$-th iteration will depend on the variable of integration of the $(l-1)$-th iteration in a way such that staying in the middle region for the long time is unlikely.

\begin{definition}\label{Def:SetupMovey}
For $t\ge0$, $u,v\in[0,t]$, set $\Delta_{u,v}(W):= W_u-W_v$, $\Delta_{u,v}(h_t):= h_t(u)-h_t(v)$ and $\Delta_{u,v} := \Delta_{u,v}(h_t)-\Delta_{u,v}(W)$.

 For $\sigma, z, z'\ge 0$ with $\sigma\le t$, $\sigma+z\le t$ and $r\in\R$, define
\begin{align}
p_{t;\sigma,h}^{(y_0)}(r) &:= \mathbb{P}_{r}\left[\mathcal{B}_{[0,t-\sigma],\Delta_{\cdot+\sigma,\sigma}}^{J_{y_0}}\right],\label{eqref:Fig21}\\
p_{t;\sigma,h}^{(y_0)}(r,z) &:= \mathbb{P}_{r}\left[\mathcal{B}_{[z,t-\sigma],\Delta_{\cdot+\sigma,\sigma}}^{J_{y_0}}\right], \label{eqref:Fig22}\\
\mu({z,\sigma,z',h}) &:= z+\Delta_{\sigma+z',\sigma}. \label{Def:mu}
\end{align}
If $\sigma=0$ it is suppressed from notation.
\end{definition}

\begin{figure}[htbp]
\centering
\includegraphics{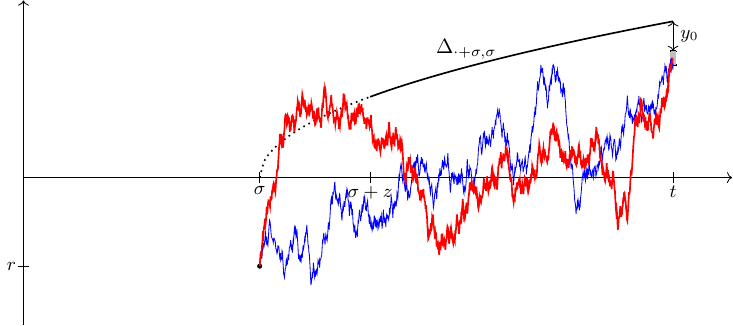}
\caption{The events in Definition \ref{Def:SetupMovey}. The event in \eqref{eqref:Fig21} corresponds to a Brownian motion starting at time $\sigma$ and location $r$ staying below the black line from time $\sigma$ to time $t$ and ending up in the gray interval, the thin blue curve is a sample in this event. The event in \eqref{eqref:Fig22} is the same, but ignoring the black barrier until time $\sigma+z$, the thick red curve serves as a sample.}
\end{figure}
We recall \eqref{Def:Gausz}. A direct consequence of the Markov property for Brownian motion is that  for $z\in \R$ and $z',\sigma\ge 0$ with $\sigma+z'\le t$ we have that
\begin{equation}
p_{t;\sigma,h}^{(y_0)}(z,z') =\int_{-\infty}^0 g_{\mu(z,\sigma,z',h),z'}(r)p_{t;\sigma+z',h}^{(y_0)}(r)\;\mathrm{d}r.\label{eq: Iter}
\end{equation}

\begin{definition}\label{Def:Lly}
Fix $\lambda>0$ and $k = \left\lceil (2\lambda\log(t))^{1/3}+1\right\rceil$. Set $z_0 :=y$. For $l\in\N$ and  $z_1,\dots, z_l\in \R$ set $\sigma_l(z_0,\dots, z_l) := \sum_{j=0}^l z_j^2$ and for $j\le l$ set $L_j(z_{j-1}) := \Cnull z_{j-1}\log(|z_{j-1}|)$.

Introduce the shorthand $\mu_j := \mu(z_j, \sigma_{j-1}, z_j^2,h)$, using the notation introduced in \eqref{Def:mu}.

The dependence of $L_j$ and $\sigma_j$ on $z_0,\dots, z_{j}$ will be omitted in the following. Let  $L^{(l)}(y) := L_l(L^{(l-1)}(y))$, $L^{(0)}(y) = y$. Finally, set $c_1(t) := c_1t^{1/5}(\log(t))^{1/2}$.

\end{definition}

\begin{figure}
\centering
\subfigure[The black line is the start of a trajectory in \eqref{eq: Main1}, the combination of the black line until $y^2$ and the thin green line is the start of a trajectory in \eqref{eq: Main4} and the combination of the black line until $y^2$ and the thick blue line is the start of a trajectory in \eqref{eq: Main2}]
{\includegraphics[scale = 0.8]{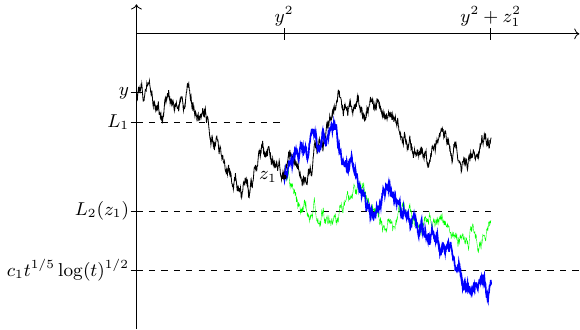}}
\hskip0.5cm
\subfigure[The thick blue line is the start of a trajectory in \eqref{eq: Main3}, the thin red one  of a trajectory in \eqref{eq: Main5}.]
{\includegraphics[scale = 0.8]{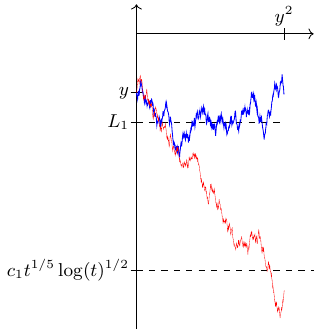}}
\caption{Sketching the summands in \eqref{eq: Main1} to \eqref{eq: Main5} for $k=2$ and $h = W = 0$.}

\label{fig:BigSum}
\end{figure}

We claim that under the assumption that $kc_1(t)\le t$, which will be proved to hold for big enough $t$ in Lemma \ref{Lem: Choicek},
\begin{align}
&p_{t,h}^{(y_0)}(y,y^2)\notag \\
&\le \sum_{l=2}^k \int_{-c_1(t)}^{L_1}\dots \int_{-c_1(t)}^{L_{l-1}} \int_{L_l}^{0} \prod_{i=1}^l \left(g_{\mu_{i-1},z_{i-1}^2}(z_i)\right)p_{t;\sigma_{l-1},h}^{(y_0)}(z_l) \;\mathrm{d}z_l\dots\mathrm{d}z_1 \label{eq: Main1}\\
&+\sum_{l=2}^{k-1}\int_{-c_1(t)}^{L_1}\dots\int_{-c_1(t)}^{L_{l-1}}\int^{L_l\wedge (-c_1(t))}_{-\infty}  \prod_{i=1}^l \left(g_{\mu_{i-1},z_{i-1}^2}(z_i)\right)p_{t;\sigma_{l-1},h}^{(y_0)}(z_l)\;\mathrm{d}z_l\dots\mathrm{d}z_1 \label{eq: Main2}\\
&+\int_{-c_1(t)}^{L_1}\dots \int^{L_{k-1}}_{-c_1(t)}\int_{-\infty}^{L_k} \prod_{i=1}^k \left(g_{\mu_{i-1},z_{i-1}^2}(z_i)\right)\;\mathrm{d}z_k\dots\mathrm{d}z_1 \label{eq: Main4}\\
&+\int_{L_1}^{0} g_{\mu(y,0, y^2,h),y^2}(z_1)p_{t;y^2,h}^{(y_0)}(z_1)\;\mathrm{d}z_1 \label{eq: Main3}\\
&+\int_{-\infty}^{L_1\wedge (-c_1(t))} g_{\mu(y,0, y^2,h),y^2}(z_1)p_{t;y^2,h}^{(y_0)}(z_1)\;\mathrm{d}z_1. \label{eq: Main5}
\end{align}
where all integrals for which the lower limit is bigger than the upper limit are defined to be 0. The previous display follows by repeatedly applying \eqref{eq: Iter} and noticing that $p_{t;\sigma_{l-1},h}^{(y_0)}(z_l)\le p_{t;\sigma_{l-1},h}^{(y_0)}(z_l, z_l^2)$  and that the condition $kc_1(t)\le t$ guarantees that $\sigma_{l-1}+z_l^2\le t$ on the region of integration considered. Example trajectories for the summands \eqref{eq: Main1}--\eqref{eq: Main5} are displayed in Figure \ref{fig:BigSum}.

\begin{remark}
In terms of the iteration described at the beginning of this section the summands in \eqref{eq: Main1} correspond to the probability to stay in the ``middle'' region for the first $l-1$ steps and then end up in $[L_l,0]$. 

The summands in \eqref{eq: Main2} correspond to the probability to stay in the middle region for the first $l-1$ steps and then end up very low in the $l$-th step, this can be controlled, since in the $l$-th step a big jump is made.

The term in \eqref{eq: Main4} corresponds to the probability that a particle stays in the middle region for so long, that even though it never made a big jump, the particle has ended up very low. Here in every step, we pick up a factor and will need to prove, that the product of these factors is small. For this, the choice of $k$, i.e.\@ $k$ being large, will be relevant.

The term in \eqref{eq: Main3} corresponds to landing in $[0,L_1]$ in the first step. Here we can directly control the integrand without needing to iterate.

The term in \eqref{eq: Main5} corresponds to making a big jump in the first step, which has low probability.
\end{remark}

We proceed by proving properties of the $k$ we have chosen.

\begin{lemma}\label{Lem: Choicek}
For $k = \left\lceil (2\lambda\log(t))^{1/3}+1\right\rceil$ and $t\ge (3c_1\lambda)^\frac{10}{3}$, we have that
\begin{equation}
\frac{1}{8}\sum_{l=1}^k \frac{L^{(l)}(y)^2}{L^{(l-1)}(y)^2}\ge 2\lambda\log(t). \label{eq: Bedk}
\end{equation}
 as well as $kc_1(t)<t$.
\end{lemma}
\begin{proof}
It is immediate from the definitions that $L^{(l)}\ge \Cnull ^l\cdot |y|$, where we use that $\log(|y|)\ge 1$. Thus we have that
\begin{align*}
\frac{L^{(l)}(y)^2}{L^{(l-1)}(y)^2} &= \Cnull ^2\log(|L_{l-1}(y)|)^2\ge \Cnull ^2(l-1)^2\log(\Cnull )\ge \Cnull (l-1)^2
\end{align*}
where in the last step we used that $\Cnull \ge e$. Thus
\[
\frac{1}{8}\sum_{l=1}^k \frac{L^{(l)}(y)^2}{L^{(l-1)}(y)^2}\ge \frac{1}{8}\Cnull \sum_{l=1}^{k-1}l^2\ge (k-1)^3\ge 2\lambda\log(t)
\]
by the choice of $k$. Furthermore, we have that 
\[kc_1(t)\le (3\lambda\log(t))^{1/3}c_1(t)\le t,\]
since we assume that $t\ge (3c_1\lambda)^\frac{10}{3}$.
\end{proof}
In the following, we will consider the summands \eqref{eq: Main1}--\eqref{eq: Main5} in the decomposition one by one and provide bounds for them. 
\subsubsection{Bounds for \eqref{eq: Main1} to \eqref{eq: Main5}, proof of theorem \ref{Theo: MainOS} assuming these bounds}

\begin{lemma} \label{Lem: FirstSum}
For $t\ge (2\sqrt{\lambda+1})^{40}$, $ y\in [-2\log(t)^{2}, -e]$ and $y_0\le0$ we have
\begin{align*}
&\sum_{l=2}^k \int_{-c_1(t)}^{L_1}\dots \int_{-c_1(t)}^{L_{l-1}} \int_{L_l}^{0} \prod_{i=1}^l \left(g_{\mu_{i-1},z_{i-1}^2}(z_i)\right)p_{t;\sigma_{l-1},h}^{(y_0)}(z_l) \;\mathrm{d}z_l\dots\mathrm{d}z_1 \\
&\le \ConstWnW^{-1}4^{\ExpPos+2}y^{2\ExpPos} \inf_{x\in[y,-e]} p_{t,h}^{(y_0)}(x,y^2).
\end{align*}
\end{lemma}
The main step in the proof of Lemma \ref{Lem: FirstSum} will be to bound $\int_{L_l}^0 g_{\mu_{l-1},z_{l-1}^2}(z_l)\;\mathrm{d}z_l$ by replacing $g_{\mu_{l-1},z_{l-1}^2}(z_l)$ by $g_{\mu(x,0,\sigma_{l-1},h),\sigma_{l-1}}(z_l)$, which can be done since we have $z_l\in [L_l,0]$.

\begin{lemma}\label{Lem: SecondSum}
For $t>0$, $y \in[-2\log(t)^{2},-e]$, $y_0\le0$, we have that
\begin{align*}
&\sum_{l=2}^{k-1}\int_{-c_1(t)}^{L_l}\dots\int_{-c_1(t)}^{L_{l-1}}\int^{L_l\wedge (-c_1(t))}_{-\infty}  \prod_{i=1}^l \left(g_{\mu_{i-1},z_{i-1}^2}(z_i)\right)p_{t;\sigma_{l-1},h}^{(y_0)}(z_l)\;\mathrm{d}z_l\dots\mathrm{d}z_1\le 8ke^{-\frac{c_1^2\log(t)}{8}}.
\end{align*}
\end{lemma}
For proving Lemma \ref{Lem: SecondSum} it suffices to bound $p_{t;\sigma_{l-1},h}^{(y_0)}(z_l)$ by 1 and to use, that $|z_l|$ is big on the region of integration of the last integral, in particular the bound only involves the last two integrals, for the other integrals all we use is that the $g_{\mu_{i-1},z_{i-1}^2}$ are density functions.
\begin{lemma}\label{Lem: ThirdTerm}
For $t\ge (3c_1\lambda)^\frac{10}{3}$, $ y\in[-2\log(t)^{2},-e]$, $y_0\le 0$, we have that
\[
\int_{-c_1(t)}^{L_1}\dots \int^{L_{k-1}}_{-c_1(t)}\int_{-\infty}^{L_k} \prod_{i=1}^k \left(g_{\mu_{i-1},z_{i-1}^2}(z_i)\right)\;\mathrm{d}z_k\dots\mathrm{d}z_1 \le e^{-2\lambda\log(t)}.
\]
\end{lemma}
To prove Lemma \ref{Lem: ThirdTerm} we can again use $p_{t,h}^{(y_0)}(z_l)\le 1$, but here the last two integrals alone aren't enough. Instead we prove by induction that in the $l$-th integral we pick up a factor $e^{-\frac{1}{8}\frac{L^{(l)}(y)^2}{L^{(l-1)}(y)^2}}$, \eqref{eq: Bedk} then allows us to conclude.

\begin{lemma}\label{Lem: StartBound}
For $t>0$, $ y\in[ -2\log(t)^{2}, -e]$, $x\in [y,-e]$, $y_0\le 0$, we have that
\[
\int_{L_1}^{0} g_{\mu(y,0, y^2,h),y^2}(z_1)p_{t;y^2,h}^{(y_0)}(z_1)\;\mathrm{d}z_1 \le Ce^{4\Coh }e^{3\Cnull \log(|y|)}\inf_{x\in[y,-e]} p_{t,h}^{(y_0)}(x,y^2).
\]
\end{lemma}
The proof of Lemma \ref{Lem: StartBound} closely mirrors the proof of Theorem \ref{Theo:MainThLB}.

\begin{lemma} \label{Lem: LastSum}
There is a $t_0(\lambda)>0$ such that for $t\ge t_0(\lambda)$, $y\in [-2\log(t)^2,-e]$, $y_0\le 0$,
\[
\int_{-\infty}^{L_1\wedge (-c_1(t))} g_{\mu(y,0, y^2,h),y^2}(z_1)p_{t;y^2,h}^{(y_0)}(z_1)\;\mathrm{d}z_1\le 8e^{-\frac{c_1^2t^{2/5}\log(t)}{8y^2}}.
\]
\end{lemma}
Lemma \ref{Lem: LastSum} is quickly proved by applying Gaussian tail estimates.

\begin{proof}[Proof of Theorem \ref{Theo: MainOS} assuming Lemma \ref{Lem: FirstSum} to Lemma \ref{Lem: LastSum}] 
By using Lemmata \ref{Lem: FirstSum}--\ref{Lem: LastSum} and \eqref{eq: Main1}--\eqref{eq: Main5} and $y\in [-2\log(t)^2,-e]$, we conclude that
\begin{align}
\frac{p_{t,h}^{(y_0)}(y,y^2)}{p_{t,h}^{(y_0)}(x,y^2)}&\le \ConstWnW^{-1}4^{\ExpPos+2}e^{4\Coh }|y|^{2\ExpPos+3\Cnull }\notag\\
&\qquad+Cp_{t,h}^{(y_0)}(x,y^2)^{-1}\cdot \left(ke^{-\frac{c_1^2\log(t)}{8}}+e^{-\frac{c_1^2}{32}t^{2/5}\log(t)^{1-4}}+e^{-2\lambda\log(t)}\right).\label{eq:Pf53EinzEq}
\end{align}
We recall our choice $c_1 =\sqrt{8(\lambda+1)}$ in \eqref{eq:c3}. Furthermore, we note that by definition $k\le t$ for $t\ge ((2\lambda)^{1/3}+2)^{3/2}$, that for $t$ big enough $\frac{c_1^2}{4}t^{2/5}\log(t)^{1-4}\ge \lambda\log(t)$ and that by \eqref{eq: OS1} $p_{t,h}^{(y_0)}(x,y^2)^{-1}\le p_{t,h}^{(y_0)}(x)^{-1}$. Combining these observations with \eqref{eq:Pf53EinzEq} finishes the proof.
\end{proof}
\subsubsection{Details}

\paragraph{Preparation}
We recall \eqref{Def:Gausz} and Definition \ref{Def:SetupMovey}. The aim of this paragraph is to prove the following lemma, which allows replacing $g_{\mu_{l-1},z_{l-1}^2}(z_l)$ by $2g_{0,4 z_{l-1}^2}(z_l)$.
\begin{lemma}\label{Lem:Removemu}
Fix $l\in\N$. Assume that for all $r\le l-1$ we have that $z_r\le \Cnull z_{r-1}\log(|z_{r-1}|)$. Then we have that for all $r \le l-1$,
\begin{align*}
g_{\mu_{r},z_{r}^2}(z_{r+1}) &\le 2g_{0, 4 z_{r}^2}(z_{r+1}).
\end{align*}
\end{lemma}
To prove Lemma \ref{Lem:Removemu} we will compare $\sigma_r$ to $z_r^2$, $z_{r+1}^2$ and establish bounds on $\Delta_{\sigma_{i-1},\sigma_{i-2}}(f)$, $f\in \{W,h_t\}$, these are stated in the next three lemmata, whose proofs are postponed since they are mostly calculation.

\begin{lemma}\label{Lem: IneqSig}
Fix $l\in\N$. Assume that for all $r\le l-1$ we have that $z_r\le \Cnull z_{r-1}\log(|z_{r-1}|)$. Then we have that for all $r \le l-1$
\begin{align}
\sigma_r &\le 2z_r^2, \label{eq:sigma1}\\
\sigma_{r-1} &\le 2\frac{z_r^2}{\log(|z_{r-1}|)^2}. \label{eq:sigma2}
\end{align}
\end{lemma}

\begin{lemma}\label{Lem: WExp}
Fix $l\in\N$. Assume that for all $r\le l-1$ we have that $z_r\le \Cnull z_{r-1}\log(|z_{r-1}|)$. Then we have that for all $r \le l-1$
\begin{align*}
 |\Delta_{\sigma_{r-1},\sigma_{r-2}}(W)|\le \frac{1}{8}|z_r-z_{r-1}|
\end{align*}
\end{lemma}

\begin{lemma}\label{Lem: hExp}
Fix $l\in\N$. Assume that for all $r\le l-1$ we have that $z_r\le \Cnull z_{r-1}\log(|z_{r-1}|)$. Then we have that for all $r \le l-1$
\[
|\Delta_{\sigma_{r-1},\sigma_{r-2}}(h_t)|\le \frac{1}{8}|z_r-z_{r-1}|.
\]
\end{lemma}

\begin{proof}[Proof of Lemma \ref{Lem:Removemu} assuming Lemmata \ref{Lem: WExp}, \ref{Lem: hExp} ]
Follows by plugging Lemma \ref{Lem: WExp}, Lemma \ref{Lem: hExp} into the definition of $g_{\mu_{r-1}, z_{r-1}^2}(z_r)$ and noticing that $z_r\le 4z_{r-1}$ by \eqref{Ass:AssC0W}.
\end{proof}

\begin{proof}[Proof of Lemma \ref{Lem: IneqSig}]
Equation \eqref{eq:sigma2} is implied by \eqref{eq:sigma1} using that $z_{r+1}^2\ge \log(|z_r|)^2z_r^2$. Thus it remains to prove \eqref{eq:sigma1}. We have that $z_r\le \Cnull z_{r-1}\log(|z_{r-1}|)$ and since $z_{r-1}\le y\le -e$ we have that $\log(|z_{r-1}|)\ge 1$ and thus 
\[
z_r\le \Cnull z_{r-1}<0.
\]
Squaring and iterating this bound gives that for $k\in\{0,\dots, r\}$ we have $z_{r-k}^2\le \frac{1}{\Cnull ^{2k}}z_r^2$ which yields that
\[
\sigma_r = \sum_{k=0}^r z_k^2 \le \sum_{k=0}^r \frac{1}{\Cnull ^{2k}}z_r^{2}\le \frac{1}{1-1/\Cnull ^2}z_r^2\le 2z_r^2
\]
where the second to last step uses that $\Cnull ^2>1$ and the last step uses that $1-1/\Cnull ^2\ge 1/2$ since $\Cnull ^2\ge 2$.
\end{proof}

\begin{proof}[Proof of Lemma \ref{Lem: WExp}]
Using the triangle inequality and $|z_{i-1}|\le \frac{1}{\Cnull } |z_{i}|$ shows that it suffices to prove that for $k\in\{1,2\}$
\begin{equation*}
|W_{\sigma_{i-k}}|\le \frac{\Cnull -1}{16\Cnull }|z_{i}|.
\end{equation*}
By \eqref{Ass:AssC0W} we have that $(\Cnull -1)/\Cnull \ge \frac{1}{2}$ and thus it suffices to prove that
\begin{equation}
\frac{|W_{\sigma_{i-k}}|}{|z_i|}\le \frac{1}{32} \label{eq:DWterm}
\end{equation}
 for $k\in\{1,2\}$. We note that for $i = 1$ and $k = 2$ \eqref{eq:DWterm} is trivial, since $W_{\sigma_{-1}} = W_0 = 0$. For the other pairs $i,k$ by Lemma \ref{Lem: IneqSig} we have that
\[
0<\sigma_{i-k}\le 2z_{i-1}^{2},
\] which implies that
\begin{align*}
\frac{|W_{\sigma_{i-k}}|}{|z_{i}|} = \frac{|W_{\sigma_{i-k}}|}{\sigma_{i-k}^{1/2}\log(\sigma_{i-k})}\frac{\sigma_{i-k}^{1/2}\log(\sigma_{i-k})}{|z_i|}\le \ClogW \frac{2|z_{i-1}|\log(2z_{i-1}^2)}{\Cnull |z_{i-1}|\log(|z_{i-1}|)}\le \frac{1}{32},
\end{align*}
where the last step uses \eqref{Ass:AssC0W}. This yields \eqref{eq:DWterm} and finishes the proof.
\end{proof}

\begin{proof}[Proof of Lemma \ref{Lem: hExp}]
As in Lemma \ref{Lem: WExp} it suffices to prove that
\[
\frac{|h_t(\sigma_{i-k})|}{|z_i|}\le \frac{1}{32}
\]
for $k\in\{1,2\}$. 
By \eqref{Ass:WachsH} we have that for $k\in\{1,2\}$
\begin{align*}
|h_t(\sigma_{i-k})| &\le \Coh |\sqrt{1+\sigma_{i-k}}|\stackrel{L. \ref{Lem: IneqSig}}{\le} \Coh \sqrt{1+2z_{i-k}^2}\stackrel{|z_{i-k}|\ge 1}\le \Coh \sqrt{3}|z_{i-k}|\le 2\Coh \frac{\sqrt{3}}{\Cnull }|z_i|,
\end{align*}
where the last step uses our assumption that $z_3\le C_3z_{r-1}\log(|z_{r-1}|)\le 0$ and the fact that $\log(|z_{r-1}|)\ge 1$, since $z_{r-1}\le y\le -e$.
Using \eqref{Ass:AssC0W} this yields the claim.
\end{proof}

\paragraph{Proof of Lemma \ref{Lem: FirstSum}}
The proof of this lemma is split into two steps. First, we bound the last two integrals in \eqref{eq: Main1}. For this, we control $g_{\mu_{l-1}, z_{l-1}^2}(z_l)g_{\mu(x,0,\sigma_{l-1},h),\sigma_{l-1}}(z_l)^{-1}$ for $z_l\in [L_l,0]$, giving the following lemma.

\begin{lemma}\label{Lem: Last2}
 Let $t\ge (2\sqrt{\lambda+1})^{40}$, $y\in [-2\log(t)^{2}, -e]$, $y_0\le 0$ and $l\ge 2$. Assume that for all  $r\le l-1$ we have that $z_r\le \Cnull z_{r-1}\log(|z_{r-1}|)$. Then
\begin{align*}
&\int_{-c_1(t)}^{L_{l-1}} g_{\mu_{l-2},z_{l-2}^2}(z_{l-1})\int_{L_l}^{0} g_{\mu_{l-1}, z_{l-1}^2}(z_l)p_{t;\sigma_{l-1},h}^{(y_0)}(z_l)\;\mathrm{d}z_l\mathrm{d}z_{l-1} \\
&\le \ConstWnW^{-1}4^{\ExpPos+2}|z_{l-2}|^{2\ExpPos}e^{-\frac{\Cnull ^2}{16}} \inf_{x\in[y,-e]}p_{t,h}^{(y_0)}(x,y^2)
\end{align*}
\end{lemma}

We note that the condition $t\ge (2\sqrt{\lambda+1})^{40}$ is not optimal, and should be thought of as $t$ big enough in a way, which does not depend on $W$.

After this, we bound the other integrals in \eqref{eq: Main1} by induction giving the following lemma.

\begin{lemma}\label{Lem: Ind}
For $y\le \xi_0$ arbitrary and all $l\in\N$, we have that 
\[
\int_{-\infty}^{L_1}\dots \int_{-\infty}^{L_l}  z_l^{2\ExpPos}\prod_{i=1}^l g_{\mu_{i-1}, z_{i-1}^2}(z_i)\;\mathrm{d}z_l\dots\mathrm{d}z_1\le \sqrt{8}e^{-\frac{l \Cnull ^2}{16}}|y|^{2\ExpPos}.
\]
\end{lemma}

\begin{proof}[Proof of Lemma \ref{Lem: FirstSum} assuming Lemmata \ref{Lem: Last2} and \ref{Lem: Ind}]
Combining Lemma \ref{Lem: Last2} and Lemma \ref{Lem: Ind} gives that
\begin{align*}
&\int_{-c_1(t)}^{L_1}\dots \int_{-c_1(t)}^{L_{l-1}} \int_{L_l}^{0} \prod_{i=1}^l \left(g_{\mu_{i-1},z_{i-1}^2}(z_i)\right)p_{t;\sigma_{l-1},h}^{(y_0)}(z_l) \;\mathrm{d}z_l\dots\mathrm{d}z_1\\
&\le  \ConstWnW^{-1}4^{\ExpPos+2}|y|^{2\ExpPos}e^{-\frac{(l-1)\Cnull ^2}{16}} \inf_{x\in[y,-e]} p_{t,h}^{(y_0)}(x,y^2).
\end{align*}
Summing this inequality over $l\ge 2$ directly yields the claim of Lemma \ref{Lem: FirstSum} since by \eqref{Ass:AssC0W} we have $1/\left(1-e^{-\frac{\Cnull ^2}{16}}\right)\le 2$.
\end{proof}
\subparagraph{Proof of Lemma \ref{Lem: Last2}}
\begin{lemma}\label{Lem: CrudeDensQuot}
We have that for $y\in[-2\log(t)^2,-e]$, $x\in [y,-e]$, $l\ge 2$,  $z_l\in [L_l,0]$, $z_{l-1}\le L_{l-1}$,
\begin{align*}
&\frac{g_{\mu_{l-1}, z_{l-1}^2}(z_l)}{ g_{\mu(x,0,\sigma_{l-1},h),\sigma_{l-1}}(z_l)}\le \sqrt{2}e^{\Cnull ^{3/2}\left(\log(|z_{l-2}|)\right)^{3/2}}.
\end{align*}
\end{lemma}
\begin{proof}
Using the definition of $g$ in \eqref{Def:Gausz}, dropping negative terms, using that $z_l\ge L_l$ and simplifying gives the bound
\begin{align}
&g_{\mu_{l-1}, z_{l-1}^2}(z_l)g_{\mu(x,0,\sigma_{l-1},h),\sigma_{l-1}}(z_l)^{-1}\notag\\
&\hskip4cm\le\frac{\sqrt{\sigma_{l-1}}}{|z_{l-1}|}e^{\frac{L_l}{z_{l-1}}+\frac{L_l}{z_{l-1}}\left(\frac{\sigma_{l-2}(|W_{\sigma_{l-1}}|+|h_t(\sigma_{l-1})|)}{\sigma_{l-1}|z_{l-1}|}+\frac{|W_{\sigma_{l-2}}|+|h_t(\sigma_{l-2})|}{|z_{l-1}|} \right)}\notag\\
&\hskip4cm\qquad\cdot e^{\frac{x^2+h_t(\sigma_{l-1})^2+ W_{\sigma_{l-1}}^2+2|x|\cdot |h_t(\sigma_{l-1})|+2|x|\cdot |W_{\sigma_{l-1}}|+2 |h_t(\sigma_{l-1})|\cdot |W_{\sigma_{l-1}}|}{2\sigma_{l-1}}}. \label{eq:FUBgmu}
\end{align}
This upper bound can be simplified further by plugging in the definition of $L_l$ and bounding the terms involving $W$ and $h$. Using  Assumptions \eqref{Ass:WachsH}, \eqref{Ass:WBMlog} and Lemma \ref{Lem: IneqSig} direct calculation gives that
\begin{align*}
\frac{(|W_{\sigma_{l-1}}|+|h_t(\sigma_{l-1})|)\cdot\sigma_{l-2}}{\sigma_{l-1}z_{l-1}} &\le 8\ClogW \left(\log(|z_{l-1}|)\right)^{1/2}+8\Coh,\\
\frac{|W_{\sigma_{l-2}}|+|h_t(\sigma_{l-2})|}{z_{l-1}} &\le 8\ClogW+8\Coh \\
\frac{x^2+W_{\sigma_{l-1}}^2+h_t(\sigma_{l-1})^2}{\sigma_{l-1}} &\le 1+\ClogW ^2\log(2z_{l-1}^2)+\Coh^2,\\
\frac{\left|xW_{\sigma_{l-1}}\right|+\left|x h_t(\sigma_{l-1})\right|}{\sigma_{l-1}}&\le \sqrt{2}\ClogW+\sqrt{2}\Coh ,\\
\frac{|h_t(\sigma_{l-1})|\cdot |W_{\sigma_{l-1}}|}{\sigma_{l-1}} &\le  \Coh \ClogW \left(\log(2z_{l-1}^2))\right)^{1/2}.
\end{align*}
Plugging these and the definition of $L_l$ in Definition \ref{Def:Lly} into \eqref{eq:FUBgmu} and simplifying yields the claim of the lemma, we note that while simplifying we also used Assumption \eqref{Ass:AssC0W} as well as that $\sqrt{\sigma_{l-1}}/|z_{l-1}| \le \sqrt{2}$ by Lemma \ref{Lem: IneqSig}.
\end{proof}

Lemma \ref{Lem: CrudeDensQuot} and \eqref{eq: OS1} allow us to bound the last integral in \eqref{eq: Main1}.
\begin{lemma}\label{Lem: BoundLastInt}
For $t\ge (2\sqrt{\lambda+1})^{40}$, $ y\in [-2\log(t)^{2}, -e]$, $y_0\le 0$, $l\ge 2$ and $z_j\in [-c_1(t),0]$ we have that for all $j\le l-1$,
\[
\int_{L_l}^{0} g_{\mu_{l-1}, z_{l-1}^2}(z_l)p_{t;\sigma_{l-1},h}^{(y_0)}(z_l)\;\mathrm{d}z_l\le \ConstWnW^{-1}\sqrt{2}2^{2\ExpPos}e^{\Cnull ^{3/2}\left(\log(|z_{l-1}|)\right)^{3/2}}|z_{l-1}|^{2\ExpPos} \inf_{x\in[y,-e]}p_{t,h}^{(y_0)}(x,y^2).
\]
\end{lemma}
\begin{proof}
Fix $x\in [y,-e]$ arbitrary. We have that 
\begin{align}
&\int_{L_l}^{0} g_{\mu_{l-1}, z_{l-1}^2}(z_l)p_{t;\sigma_{l-1},h}^{(y_0)}(z_l)\;\mathrm{d}z_l\notag \\
&\stackrel{L.\ref{Lem: CrudeDensQuot}}{\le}\sqrt{2}e^{\Cnull ^{3/2}\left(\log(|z_{l-1}|)\right)^{3/2}}\int_{-\infty}^{0}g_{\mu(x,0,\sigma_{l-1},h),\sigma_{l-1}}(z_l)p_{t;\sigma_{l-1},h}^{(y_0)}(z_l)\;\mathrm{d}z_l \notag \\
&\stackrel{\eqref{eq: Iter}}{=} \sqrt{2}e^{\Cnull ^{3/2}\left(\log(|z_{l-1}|)\right)^{3/2}}p_{t,h}^{(y_0)}(x,\sigma_{l-1}). \label{eq:7211}
\end{align}

 Since $z_j\in [-c_1(t),0]$ for all $j\le l-1$ we have that
\[
\sqrt{\sigma_{l-1}}\le \sqrt{2}|z_{l-1}|\le \sqrt{2}c_1(t) \le 2t^{1/4},
\] where the last inequality uses that $t\ge (2\sqrt{\lambda+1})^{40}$ and $c_1(t) = \sqrt{8(\lambda+1)}t^{1/5}(\log(t))^{1/2}$. Thus, we can apply \eqref{eq: OS1} to conclude
\begin{align}
p_{t,h}^{(y_0)}(x,\sigma_{l-1}&)\stackrel{\eqref{eq: OS1}}{\le } \ConstWnW^{-1}p_{t,h}^{(y_0)}(x)\sigma_{l-1}^{\ExpPos}\stackrel{L. \ref{Lem: IneqSig}}{\le} \ConstWnW^{-1}p_{t,h}^{(y_0)}(x,y^2)|z_{l-1}|^{2\ExpPos} 2^{2\ExpPos}. \label{eq:7212}
\end{align}
Together \eqref{eq:7211}, \eqref{eq:7212} and the fact that $x\in [y,-e]$ was arbitrary, finish the proof.
\end{proof}
A direct calculation gives that for  $s\ge \Cnull \log(|z_{l-2}|)$ 
\begin{equation}
s^{2\ExpPos}e^{\Cnull ^{3/2}\log(s|z_{l-2}|)^{3/2}} \le e^{\frac{s^2}{16}}. \label{eq:Fu}
\end{equation}

With this, we have all tools needed for the proof of Lemma \ref{Lem: Last2}.
\begin{proof}[Proof of Lemma \ref{Lem: Last2}]
Set
\[
I(z_{l-1}):= \int_{-c_1(t)}^{L_{l-1}} g_{\mu_{l-2},z_{l-2}^2}(z_{l-1})\int_{L_l}^{0} g_{\mu_{l-1}, z_{l-1}^2}(z_l)p_{t;\sigma_{l-1},h}^{(y_0)}(z_l)\;\mathrm{d}z_l\mathrm{d}z_{l-1}.
\] Using Lemma \ref{Lem: BoundLastInt} and Lemma \ref{Lem:Removemu} yields that
\begin{align*}
&\frac{I(z_{l-1})}{\sqrt{2}2^{2\ExpPos}\ConstWnW^{-1} p_{t,h}^{(y_0)}(x,y^2)} \\
&\stackrel{L. \ref{Lem: BoundLastInt}}{\le} \int_{-\infty}^{L_{l-1}}g_{\mu_{l-2},z_{l-2}^2}(z_{l-1})e^{\Cnull ^{3/2}\left(\log(|z_{l-1}|)\right)^{3/2}}|z_{l-1}|^{2\ExpPos}\;\mathrm{d}z_{l-1}\\
&\stackrel{L. \ref{Lem:Removemu}}{\le} \int_{-\infty}^{\Cnull z_{l-2}\log(|z_{l-2}|)} |z_{l-1}|^{2\ExpPos}2g_{0,4z_{l-2}^2}(z_{l-1})e^{\Cnull ^{3/2}\left(\log(|z_{l-1}|)\right)^{3/2}}\;\mathrm{d}z_{l-1}\\
&\stackrel{ s= \frac{z_{l-1}}{z_{l-2}}, \eqref{eq:Fu}}{\le} |z_{l-2}|^{2\ExpPos} \sqrt{8}\int_{\Cnull \log(|z_{l-2}|)}^\infty g_{0,8}(s)\;\mathrm{d}s\le 8|z_{l-2}|^{2\ExpPos}  e^{-\frac{\Cnull ^{2}\log(|z_{l-2}|)^2}{16}}\\
&\le 8|z_{l-2}|^{2\ExpPos}  e^{-\frac{\Cnull ^2}{16}},
\end{align*} 
where the second to last step used Gaussian tail estimates\footnote{See for example Theorem 1.2.6 in  \cite[p. 13]{Durrett}.} and the last that by our assumptions $|z_{l-2}|\ge |y|\ge e$.  Rearranging yields the claim.
\end{proof}
\subparagraph{Proof of Lemma \ref{Lem: Ind}}

We prove this by induction. Induction basis ($l = 1$): By applying Lemma \ref{Lem:Removemu} and \eqref{eq:Fu} we get analogously to the calculation in the proof of Lemma \ref{Lem: Last2}
\begin{align*}
\int_{-\infty}^{L_1} |z_1|^{2\ExpPos} g_{\mu_0,z_0^2}(z_1)\;\mathrm{d}z_1 &\le |z_0|^{2\ExpPos}\int_{\Cnull \log(|z_0|)}^{\infty}\sqrt{8}g_{0,8}(s)\;\mathrm{d}s\le \sqrt{8}|z_0|^{2\ExpPos} e^{-\frac{\Cnull ^2}{16}}.
\end{align*}

Induction step ($l\to l+1$): We have that 
\begin{align*}
\int_{-\infty}^{L_1}\dots \int_{-\infty}^{L_{l+1}} |z_{l+1}|^{2\ExpPos}\prod_{i=1}^l g_{\mu_{i-1}, z_{i-1}^2}(z_i)\;\mathrm{d}z_{l+1}\dots\mathrm{d}z_1&\stackrel{\text{IH}}{\le} \int_{-\infty}^{L_1} |z_1|^{2\ExpPos}e^{-\frac{l\Cnull ^2}{16}} g_{\mu_0,z_0^2}(z_1)\;\mathrm{d}z_1\\
&\stackrel{\text{IB}}{\le} e^{-\frac{l\Cnull ^2}{16}} |y|^{2\ExpPos}e^{-\frac{\Cnull ^2}{16}}= |y|^{2\ExpPos}e^{-\frac{(l+1)\Cnull ^2}{16}}.\qed
\end{align*}

\paragraph{Proof of Lemma \ref{Lem: SecondSum}}

Let $l\in \{2,\dots,k-1\}$. To keep notation slightly lighter, we will as in the proof of Lemma \ref{Lem: FirstSum} handle the last two integrals in \eqref{eq: Main2} first. Since $p_{t;\sigma_{l-1},h}^{(y_0)}(z_l)\le 1$ we have that
\begin{align}
&\int_{-c_1(t)}^{L_{l-1}}\int_{-\infty}^{-c_1 t^{1/5}\log(t)^{1/2}\wedge L_l}  \prod_{i=l-1}^l \left(g_{\mu_{i-1},z_{i-1}^2}(z_i)\right)p_{t;\sigma_{l-1},h}^{(y_0)}(z_l)\;\mathrm{d}z_l\mathrm{d}z_{l-1} \notag\\
&\le \int_{-c_1(t)}^{L_{l-1}}\int_{-\infty}^{-c_1 t^{1/5}\log(t)^{1/2}\wedge L_l}  \prod_{i=l-1}^l \left(g_{\mu_{i-1},z_{i-1}^2}(z_i)\right)\;\mathrm{d}z_l\mathrm{d}z_{l-1}\notag\\
&\stackrel{L. \ref{Lem:Removemu}}{\le} \int_{-c_1(t)}^{L_{l-1}}\int_{-\infty}^{-c_1 t^{1/5}\log(t)^{1/2}\wedge L_l} \prod_{i=l-1}^l 2g_{0,4z_{i-1}^2}(z_i)\;\mathrm{d}z_l\mathrm{d}z_{l-1}\notag\\
&\le 4\int_{-c_1(t)}^{L_{l-1}} g_{0,4z_{l-2}^2}(z_{l-1})e^{-\frac{(-c_1(t)\wedge L_l)^2}{8 z_{l-1}^2}}\;\mathrm{d}z_{l-1} \label{eq:FBL511}
\end{align}
On $[-t^{1/5}, L_{l-1}]$ we have that 
\[
e^{-\frac{(-c_1(t)\wedge L_l)^2}{8 z_{l-1}^2}} \le e^{-\frac{c_1^2\log(t)}{8}}
\] and on $[-c_1(t),-t^{1/5}]$ we have that
\begin{align*}
e^{-\frac{(-c_1(t)\wedge L_l)^2}{8 z_{l-1}^2}} &\le e^{-\frac{\Cnull ^2\log(|z_{l-1}|)^2}{8}}\le e^{-\frac{c_1^2\log(t)}{8}},
\end{align*}
where the last step used that by \eqref{Ass:AssC0W} $\Cnull ^2\ge 5c_1^2$.

Plugging both of these into \eqref{eq:FBL511} and using that $g_{\mu_{l-2},z_{l-2}^2}$ is a density function gives that 

\begin{align}
\int_{-c_1(t)}^{L_{l-1}}\int_{-\infty}^{-c_1 t^{1/5}\log(t)^{1/2}\wedge L_l}  \prod_{i=l-1}^l \left(g_{\mu_{i-1},z_{i-1}^2}(z_i)\right)p_{t;\sigma_{l-1},h}^{(y_0)}(z_l)\;\mathrm{d}z_l\mathrm{d}z_{l-1}&\le 8e^{-\frac{c_1^2\log(t)}{8}} \label{eq:UB515}
\end{align}

 Since for $i<l-1$ all $g_{\mu_{i-1}, z_{i-1}^2}$ are density functions, \eqref{eq:UB515} implies that 
\begin{align}
&\int_{-c_1(t)}^{L_l}\dots\int_{-c_1(t)}^{L_{l-1}}\int^{L_l\wedge (-c_1(t))}_{-\infty}  \prod_{i=1}^l \left(g_{\mu_{i-1},z_{i-1}^2}(z_i)\right)p_{t;\sigma_{l-1},h}^{(y_0)}(z_l)\;\mathrm{d}z_l\dots\mathrm{d}z_1 \notag \\
&\le 8e^{-\frac{c_1^2\log(t)}{8}}. \label{eq:UB516}
\end{align}

Finally, summing \eqref{eq:UB516} for $l\in \{2,\dots, k-1\}$ yields the statement of Lemma \ref{Lem: SecondSum}. \qed
\paragraph{Proof of Lemma \ref{Lem: ThirdTerm} }
One preparatory lemma is needed for this. 
\begin{lemma}\label{Lem: LastHilf}
For $m\in\N$ arbitrary, $y\le \xi_0$
\begin{align*}
\int_{-\infty}^{L_1}\dots \int_{-\infty}^{L_{m-1}} \prod_{i=1}^{m-1}\frac{1}{\sqrt{2\pi}|z_{i-1}|}e^{-\frac{z_i^2}{8z_{i-1}^2}}e^{-\frac{L_m^2}{8z_{m-1}^2}}\;\mathrm{d}z_{m-1}\dots\mathrm{d}z_1\le e^{-\frac{1}{8}\sum_{l=1}^m \frac{L^{(l)}(y)^2}{L^{(l-1)}(y)^2}}.
\end{align*}
\end{lemma}
\begin{proof}
We prove this using induction. Induction basis ($m=1$): Here no integrals remain and the statement trivially holds by definition of $L^{(1)}(y)$ and $L^{(0)}(y)$.

Induction step ($m\to m+1$): We have that
\begin{align*}
&\int_{-\infty}^{L_1}\dots \int_{-\infty}^{L_m} \prod_{i=1}^m \frac{1}{\sqrt{2\pi}|z_{i-1}|}e^{-\frac{z_i^2}{8z_{i-1}^2}-\frac{L_{m+1}^2}{8z_m^2}}\;\mathrm{d}z_m\dots\mathrm{d}z_1\\
&\stackrel{\text{IH}}{\le}\int_{-\infty}^{L_1} \frac{1}{\sqrt{2\pi}|y|}e^{-\frac{z_1^2}{8y^2}}e^{-\frac{1}{8}\sum_{l=1}^m \frac{L^{(l)}(z_1)^2}{L^{(l-1)}(z_1)^2}}\;\mathrm{d}z_1
\end{align*}
since  $\frac{L^{(l)}(x)^2}{L^{(l-1)}(x)^2} = \Cnull ^2\log(|L^{(l-1)}(x)|)^2$ is decreasing in $x$ for $x\le -1$ and $z_1\le y$, which can be seen directly from the definition in \ref{Def:Lly}, we have that 
\begin{align*}
&\int_{-\infty}^{L_1}\dots \int_{-\infty}^{L_m} \prod_{i=1}^m \frac{1}{\sqrt{2\pi}|z_{i-1}|}e^{-\frac{z_i^2}{8z_{i-1}^2}-\frac{L_{m+1}^2}{8z_m^2}}\;\mathrm{d}z_m\dots\mathrm{d}z_1\\
&\le e^{-\sum_{l=1}^{m} \frac{L^{(l)}(L_1(y))^2}{L^{(l-1)}(L_1(y))^2}}\int_{-\infty}^{L_1} \frac{1}{\sqrt{2\pi}|y|}e^{-\frac{z_1^2}{8y^2}}\;\mathrm{d}z_1\le e^{-\sum_{l=2}^{m+1} \frac{L^{(l)}(y)^2}{L^{(l-1)}(y)^2}} e^{-\frac{L_1(y)^2}{8y^2}}.\qedhere
\end{align*}
\end{proof}

\begin{proof}[Proof of Lemma \ref{Lem: ThirdTerm}]
Repeating the calculation in Lemma \ref{Lem: SecondSum} we get that
\begin{align*}
&\int_{-c_1(t)}^{L_1}\dots \int^{L_{k-1}}_{-c_1(t)}\int^{L_k}_{-\infty} \prod_{i=1}^k \left(g_{\mu_{i-1},z_{i-1}^2}(z_i)\right)\;\mathrm{d}z_k\dots\mathrm{d}z_1 \\
&\le\int_{-\infty}^{L_1}\dots \int_{-\infty}^{L_{k-1}} \prod_{i=1}^{k-1}\frac{1}{\sqrt{2\pi}|z_{i-1}|}e^{-\frac{z_i^2}{8z_{i-1}^2}}e^{-\frac{L_k^2}{8z_{k-1}^2}}\;\mathrm{d}z_{k-1}\dots\mathrm{d}z_1\stackrel{L. \ref{Lem: LastHilf}}{\le} e^{-\frac{1}{8}\sum_{l=1}^k \frac{L^{(l)}(y)^2}{L^{(l-1)}(y)^2}}\\
&\stackrel{\eqref{eq: Bedk}}{\le} e^{-2\lambda \log(t)}.\qedhere
\end{align*}
\end{proof}
\paragraph{Proof of Lemma \ref{Lem: StartBound} }
\pushQED{\qed} 
We have that
\begin{align*}
&\int_{L_1}^{0}g_{\mu(y,0,y^2,h),y^2}(z_1)p_{t;y^2,h}^{(y_0)}(z_1)\;\mathrm{d}z_1 \\
&\stackrel{\eqref{eq: LBQuotDens}}{=} e^{\frac{x^2-y^2}{2y^2}+\frac{(y-x) W_{y^2}}{y^2}-\frac{h_t(y^2)(y-x)}{y^2}}\int_{L_1}^{0}e^{\frac{z_1(y-x)}{y^2}}g_{\mu(x,0,y^2,h),y^2}(z_1)p_{t;y^2,h}^{(y_0)}(z_1)\;\mathrm{d}z_1\\
&\le e^{\frac{x^2-y^2}{2y^2}+\frac{(y-x) W_{y^2}}{y^2}-\frac{h_t(y^2)(y-x)}{y^2}} e^{\frac{\Cnull  y\log(|y|)(y-x)}{y^2}}\int_{-\infty}^0 g_{\mu(x,0,y^2,h),y^2}(z_1)p_{t;y^2,h}^{(y_0)}(z_1)\;\mathrm{d}z_1\\
&\stackrel{\eqref{eq: SeedRec}, \eqref{Ass:WBMlog}, \eqref{Ass:WachsH}}{\le}  Ce^{4\Coh }e^{2\ClogW \sqrt{\log(y^2)}}e^{2\Cnull \log(|y|)}p_{t,h}^{(y_0)}(x,y^2)\\
&\stackrel{\eqref{Ass:AssC0W}}{\le} Ce^{4\Coh }e^{3\Cnull \log(|y|)}p_{t,h}^{(y_0)}(x,y^2).\qedhere
\end{align*}
\popQED
\paragraph{Proof of Lemma \ref{Lem: LastSum}}
\pushQED{\qed} 
By definition $L_1\wedge (-c_1(t)) \le -c_1(t)$ and thus
\begin{align*}
\int_{-\infty}^{L_1\wedge (-c_1(t))} g_{\mu(y,0, y^2,h),y^2}(z_1)p_{t;y^2,h}^{(y_0)}(z_1)\;\mathrm{d}z_1 &\stackrel{L. \ref{Lem: WExp}, \ref{Lem: hExp}}{\le} \int_{-\infty}^{-c_1(t)} 2g_{0,4y^2}(z_1)\;\mathrm{d}z_1\\
&\le 8e^{-\frac{c_1^2t^{2/5}\log(t)}{8y^2}}.\qedhere
\end{align*}
\popQED

\section{Getting rid of $h_t$ -- Preparation for Section \ref{Sec:BarrComps}}\label{Sec:Moveh}
In all of this section we assume that $t\ge 0$, $y\in [-\log(t)^{2},0]$, $ y_0\in [-\log(t)^{2},0]$ are fixed.

 Furthermore, let $(X_s)_{s\ge0}$ be the universal object in what follows and assume that under $\mathbb{P}_y$ it is distributed as a Brownian motion starting at $y$. In this section, we will analyze how adding a shift $h_t$ to a barrier changes the probability that $X_s$ stays below the barrier. In Section \ref{Sec:LinShift} the additional shift will be linear, in Sections \ref{Sec:htUBRT} and \ref{Sec:htLBRT} we add a positive/negative banana. The main tool in this section is the Girsanov theorem, we define measures with regard to which $(X_s+h_t(s))_{s\in [0,t]}$ is a Brownian motion and control the corresponding Radon-Nikodym derivatives.

\subsection{Linear Shifts} \label{Sec:LinShift}

We recall Definition \ref{Def:Barr}.
\begin{definition}
Consider $(c_t)_{t\ge0}$ in $\R$ and set
\begin{align*}
A_{t,h}^{(y_0)}&:= \mathcal{B}_{[0,t],h_t(\cdot)-W_\cdot}^{J_{y_0}}(X_\cdot),\\
A_{t,h,c}^{(y_0)} &:=  \mathcal{B}_{[0,t],h_t(\cdot)-W_\cdot}^{J_{y_0}}(X_s+sc_t).
\end{align*}
\end{definition}
In this section, we prove the following theorem.
\begin{theorem} \label{Theo: LinTerm}
We have that for $t\ge 0$ and $y,y_0\le 0$,
\begin{align*}
\frac{\mathbb{P}_y[A_{t,h,c}^{(y_0)}]}{\mathbb{P}_y[A_{t,h}^{(y_0)}]} &\le  e^{-\frac{tc_t^2}{2}}e^{-yc_t-h_t(t)c_t+|y_0-1|\cdot|c_t|}e^{ W_t c_t},\\
\frac{\mathbb{P}_y[A_{t,h,c}^{(y_0)}]}{\mathbb{P}_y[A_{t,h}^{(y_0)}]} &\ge e^{-\frac{tc_t^2}{2}}e^{-yc_t-h_t(t)c_t-|y_0-1|\cdot |c_t|}e^{ W_t c_t}.
\end{align*}
\end{theorem}
\begin{proof}
Set
\begin{align*}
Z_t&:= \frac{\mathrm{d}\mathbb{Q}_y}{\mathrm{d}\mathbb{P}_y}:= \exp\left(-\int_0^t c_t\;\mathrm{d}X_s-\frac{1}{2}\int_0^t c_t^2\;\mathrm{d}s\right)=\exp\left(-c_t(X_t-X_0)-\frac{tc_t^2}{2}\right)
\end{align*}

The Girsanov Theorem then gives
\begin{align}
\mathbb{P}_y\left[A_{t,h,c}^{(y_0)}\right] &= \mathbb{E}_{\mathbb{Q}_y}\left[(Z_t)^{-1}\mathbf{1}_{A_{t,h,c}^{(y_0)}}\right]= e^{\frac{tc_t^2}{2}}\mathbb{E}_{\mathbb{Q}_y}\left[e^{c_tX_t-yc_t} \mathbf{1}_{A_{t,h,c}^{(y_0)}}\right]. \label{eq:LinctSetup}
\end{align}
On $A_{t,h,c}^{(y_0)}$ we have that
\begin{equation}
 W_t-tc_t+y_0-1-h_t(t)\le X_t\le  W_t-tc_t-h_t(t)+y_0. \label{eq:XtonAthc}
\end{equation}
Plugging \eqref{eq:XtonAthc} into \eqref{eq:LinctSetup} yields that
\begin{align*}
\mathbb{P}_y\left[A_{t,h,c}^{(y_0)}\right] &\le e^\frac{tc_t^2}{2}\mathbb{E}_{\mathbb{Q}_y}\left[e^{( W_t -tc_t-h_t(t)+y_0-\mathbf{1}_{c_t<0})c_t-yc_t}\mathbf{1}_{A_{t,h,c}}\right]\\
&\le  e^{-\frac{tc_t^2}{2}}e^{-yc_t-h_t(t)c_t+|y_0-1|\cdot|c_t|}e^{ W_t c_t}\mathbb{P}_y\left[A_{t,h}^{(y_0)}\right],\\
\mathbb{P}_y[A_{t,h,c}^{(y_0)}] &\ge e^{\frac{tc_t^2}{2}}\mathbb{E}_{\mathbb{Q}_y}\left[e^{( W_t-tc_t+y_0-\mathbf{1}_{c_t>0}-h_t(t))c_t-yc_t}\mathbf{1}_{A_{t,h,c}^{(y_0)}}\right]\\
&\ge e^{-\frac{tc_t^2}{2}}e^{-yc_t-h_t(t)c_t-|y_0-1|\cdot |c_t|}e^{ W_t c_t}\mathbb{P}_y\left[A_{t,h}^{(y_0)}\right].
\end{align*}
Rearranging yields the claim.
\end{proof} 
\subsection{General Setup}\label{Sec:GenSet}

Fix $t>0$ and $h_t:[0,t]\to \R$ with $h_t(t) = 0 = h_t(0)$. Set $\chi(t) :=\log(t)^{1.5}$. Define $\widehat{W}$ via $\widehat{W}_s := W_{t-s}-W_t$, $s\in [0,t]$. We recall  Section \ref{Subsec:DefSum}, we use $\widehat{C}$ to denote variants of constants using $\widehat{W}$ instead of $W$, for example we set $\widehat{\mathcal{C}}_{\log} := \mathcal{C}_{\log}(\widehat{W},t)$.

 From now on we will need two additional assumptions. Namely, we assume that 
\begin{equation}
\log(t)^{-0.5}\max\{\ClogW ,\ClogWhat \}\le \frac{1}{3}\label{eq:Limgut}
\end{equation}
and that there are $\ExpStart:=C_4(W,t)<\infty$, $\ExpStarthat := C_4(\widehat{W},t)<\infty$, $\ConstStart := c_2(W,t)>0$, $\ConstStarthat := c_2(\widehat{W},t)>0$ such that for all $r\le\chi(t)$
\begin{align}
\mathbb{P}_{-1}\left[\mathcal{B}_{[0,r],-W_\cdot}^{J_{0}}(X_\cdot)\right]\ge \ConstStart(1+r)^{-\ExpStart},\label{Ass:Assdoof}\\
\mathbb{P}_{-1}\left[\mathcal{B}_{[0,r],-\widehat{W}_\cdot}^{J_{0}}(X_\cdot)\right]\ge \ConstStarthat(1+r)^{-\ExpStarthat}.\label{Ass:Assdoofhat}
\end{align}

We will check in Section \ref{Sec:Constants} that these assumptions hold for the specific $W$ we use in the first half of the paper.

\begin{definition}\label{Def:pttildes}
For $t\ge 0$, $h_t:[0,t]\to \R$, $y\in [-\log(t)^{2},0]$, $ y_0\in [-\log(t)^{2},0]$ define
\begin{align*}
\widetilde{p}_{t,h}^{(y_0)}(y) &:= \mathbb{P}_y\left[\widetilde{A}^{(y_0)}_{t,h}\right]:= \mathbb{P}_y\left[\mathcal{B}_{[0,t],-W_\cdot}^{J_{y_0}}(X_\cdot+h_t(\cdot)+\frac{\cdot}{t}W_t)\right],\\
p_{t,h}^{(y_0)}(y) &:= \mathbb{P}_y\left[A_{t,h}^{(y_0)}\right]:= \mathbb{P}_y\left[\mathcal{B}_{[0,t],-W_\cdot}^{J_{y_0}}(X_\cdot+h_t(\cdot))\right]\\
p_t^{(y_0)}(y) &:= \mathbb{P}_y\left[A_t^{(y_0)}\right]:= \mathbb{P}_y\left[\mathcal{B}_{[0,t],-W_\cdot}^{J_{y_0}}(X_\cdot)\right].
\end{align*}
\end{definition}
Set
\begin{align}
Z_t &:= \frac{\mathrm{d}\mathbb{Q}_y}{\mathrm{d}\mathbb{P}_y} = \exp\left(\int_0^t -h_t'(s)-\frac{ W_t}{t}\;\mathrm{d}X_s-\frac{1}{2}\int_0^t \left(h'_t(s)+\frac{ W_t}{t}\right)^2\mathrm{d}s\right)\notag\\
&=\exp\left(-h_t'(t)X_t+h_t'(0)y-\frac{ W_t (X_t-y)}{t}+\int_0^t h_t''(s)X_s\mathrm{d}s-\frac{1}{2}\int_0^t \left(h'_t(s)+\frac{ W_t}{t}\right)^2\mathrm{d}s\right)\label{Def:Zt-1}
\end{align}
We have that 
\begin{align}
\widetilde{p}_{t,h}^{(y_0)}(y) &= \mathbb{P}_y\left[\widetilde{A}_{t,h}^{(y_0)}\right] = \mathbb{E}_{\mathbb{Q}_y}\left[(Z_t)^{-1}\mathbf{1}_{\widetilde{A}_{t,h}^{(y_0)}}\right]\label{eq:Eins}
\end{align}
and since Girsanov implies that $(X_s+h_t(s)+\frac{s}{t}W_t)_{s\le t}$ is a Brownian motion with regard to $\mathbb{Q}_y$ \eqref{eq:Eins} reduces controlling  $\widetilde{p}_{t,h}^{(y_0)}(y)/p_t^{(y_0)}(y)$ to controlling $(Z_t)^{-1}$ on $\mathbf{1}_{\widetilde{A}_{t,h}^{(y_0)}}$. Furthermore, Section \ref{Theo: LinTerm} allows us to control $\widetilde{p}_{t,h}^{(y_0)}(y)/p_{t,h}^{(y_0)}(y)$ since applying Theorem \ref{Theo: LinTerm} with $c_t = \frac{W_t}{t}$ yields that for $t\ge e$,  $y, y_0\in [- 2\log(t)^2,0]$
\begin{align}
\frac{\widetilde{p}_{t,h}^{(y_0)}(y)}{p_{t,h}^{(y_0)}(y)} &\le e^{\frac{W_t^2}{2t}}e^{-y\frac{W_t}{t}}e^{|y_0-1|\cdot\left|\frac{W_t}{t}\right|}\stackrel{\eqref{Ass:WBMlog}}{\le }e^{\frac{W_t^2}{2t}}e^{10\ClogW }, \label{eq:linUB}\\
\frac{\widetilde{p}_{t,h}^{(y_0)}(y)}{p_{t,h}^{(y_0)}(y)} &\ge e^{\frac{W_t^2}{2t}}e^{-y\frac{W_t}{t}}e^{-|y_0-1|\cdot\left|\frac{W_t}{t}\right|}\stackrel{\eqref{Ass:WBMlog}}{\ge} e^{\frac{W_t^2}{2t}}e^{-10\ClogW }\label{eq:linLB}.
\end{align}
\subsection{Bound Used in the Upper Bound of the Right Tail}\label{Sec:htUBRT}
In this section, we consider $h^{\smallfrown}_t : [0,t]\to \R, s\mapsto -\left(((1+s)\wedge (1+t-s))^{1/6}-1\right)/2$. 

The aim of this section is to derive an upper bound on $p_{t,h^{\smallfrown}}^{(y_0)}(y)/p_t^{(y_0)}(y)$, to formulate this upper bound we need some additional definitions. We recall \eqref{Ass:WBMlog}, \eqref{Ass:QuenchedWall}. \eqref{Def:WeirdC1} and Definition \ref{Def:pttildes}.
\begin{definition}\label{Def:ptast}
 For $\widehat{W} := (W_{t-s}-W_t)_{s\in [0,t]}$ and $y,y_0\in [-\log(t)^2,0]$ set
\begin{align*}
(\ptast)^{-1} &:= p^\ast_t(W,t,h^{\smallfrown})^{-1} := e^{5\ClogW +\ClogW ^2}\ConstWnW^{-1}\ConstStart^{-1}(1+\chi(t))^{\ExpStart} \log(t)^{4\ExpPos}p_{t,h^{\smallfrown}}^{(y_0)}(y)^{-1}, \\
(\ptasthat)^{-1} &:=\widehat{p}^\ast_t(W,t,h^{\smallfrown})^{-1} := e^{5\ClogWhat +\ClogWhat ^2}\ConstWnWhat^{-1}\ConstStarthat^{-1}(1+\chi(t))^{\ExpStarthat} \log(t)^{4\ExpPoshat }\widehat{p}_{t,h^{\smallfrown}}^{(y_0)}(y)^{-1},
\end{align*}
where $\widehat{p}_t^{(y_0)}(y)$ uses $\widehat{W}$ instead of $W$. 

Choose $\kappa := \kappa(W,\widehat{W},t,h^{\smallfrown},y,y_0)>0$ such that
\begin{equation}
\forall_{s\le \chi(t)}\qquad -\kappa(1+s)^\frac{2}{3}-h^{\smallfrown}_t(s)\le -2(\max\{\ClogW ,\ClogWhat \})(1+s)^{\frac{2}{3}}\label{Def:Kappa}
\end{equation}
Furthermore, fix $1/8\ge\delta>0$ and choose $\eta := \eta(W,\widehat{W},t,h^{\smallfrown},\delta,y,y_0)<0$ such that
\begin{align}
&\sum_{l=1}^{\infty} 6e^{-\frac{\kappa^2 l^{\frac{4}{3}-1}+\eta^2l^{-1}}{2}}e^{15\ClogW+\ClogW^2}\ConstStart^{-1}(l+1)^{\ExpStart}\ConstWnW^{-2}e^{4\Coh }4^{\ExpPos+2}\left(|y|+3\ClogW(1+l)^\frac{2}{3}+|\eta|\right)^{4\ExpPos+3\Cnull}\notag\\
&\le \delta,\label{eq:eta0}\\
&\sum_{l=1}^{\infty} 6e^{-\frac{\kappa^2 l^{\frac{4}{3}-1}+\eta^2l^{-1}}{2}}e^{15\ClogWhat+\ClogWhat^2}\ConstStarthat^{-1}(l+1)^{\ExpStarthat}\ConstWnWhat^{-2}e^{4\Coh }4^{\ExpPoshat+2}\left(|y_0|+3\ClogWhat(1+l)^\frac{2}{3}+|\eta|\right)^{4\ExpPoshat+3\Cnullhat}\notag\\
&\le \delta.\label{eq:eta1}
\end{align}

\end{definition}
In this section, we will prove the following theorem.
\begin{theorem}\label{Sa:GirsConvH}
There exist $c,C>0$ such that for $\lambda>0$ arbitrary there is a $t_0(\lambda)\ge0$ such that  for $t\ge t_0(\lambda)$ with \eqref{eq:Limgut} and $y,y_0\in [-\log(t)^2,0]$,
\begin{align*}
&\frac{p_{t,h^{\smallfrown}}^{(y_0)}(y)-\frac{\delta}{1-\delta}\left(\sum_{j=0}^2p_{t,h^{\smallfrown}}^{(y_0-1)}(y-j)\right)}{p_t^{(y_0)}(y)}  \\
&\le Ce^{c(\ClogW+\kappa-\eta)}\left(1+t^{-\lambda}p_t^{(y_0)}(y)^{-1}\right)+Ct^{-\lambda}((\ptast)^{-1}+(\ptasthat)^{-1})p_t^{(y_0)}(y)^{-1}.
\end{align*}

\end{theorem}
\begin{remark}\label{Rem:Rem1}
We note that $\eta$ and $\kappa$ can be chosen as continuous functions of $y$, $y_0$, $\Coh$, $\ClogW $, $\ClogWhat $, $\ConstWnW^{-1}$, $\ConstWnWhat^{-1}$, $\ExpPos$, $\ExpPoshat $, $\ExpStart$, $\widehat{\ExpStart}$, $\ConstStart^{-1}$, $\widehat{\ConstStart}^{-1}$, $\delta$. 
\end{remark}
\begin{definition}\label{Def:rho}
To slightly shorten the following definitions shorten $X_s^{h^{\smallfrown}} := X_s+h^{\smallfrown}_t(s)$. Using this define
\begin{align*}
\widetilde{A}_{t,h^{\smallfrown},\text{start}} &:= \{\forall_{s\le \chi(t)}\, X_s^{h^{\smallfrown}}\ge -\kappa(1+s)^\frac{2}{3}+y+\eta\},\\
\widetilde{A}_{t,h^{\smallfrown},\text{end}} &:= \{\forall_{s\ge t-\chi(t)}\, X_s^{h^{\smallfrown}}\ge -\kappa(1+t-s)^\frac{2}{3}+X_t+\eta\},\\
\widetilde{A}_{t,h^{\smallfrown},\chi(t)}^{(y_0)} &:= \widetilde{A}_{t,h^{\smallfrown}}^{(y_0)}\cap\widetilde{A}_{t,h^{\smallfrown},\text{start}}\cap\widetilde{A}_{t,h^{\smallfrown},\text{end}},\\
B_{t,h^{\smallfrown},\chi(t),k} &:= \Bigg\{\min\Big\{\min_{\chi(t)\le s\le t/2} \frac{X_s^{h^{\smallfrown}}+ s/tW_t-y}{(1+s)^\frac{2}{3}},\min_{t/2\le s\le t-\chi(t)} \frac{X_s^{h^{\smallfrown}}- (t-s)/t W_t-y_0}{(1+t-s)^\frac{2}{3}}\Big\}\\
&\hskip9cm\in [-k,-k+1) \Bigg\},\\
\widetilde{p}_{t,h^{\smallfrown},\chi(t)}^{(y_0)}(y) &:= \mathbb{P}_y\left[\widetilde{A}_{t,h^{\smallfrown},\chi(t)}^{(y_0)}\right].
\end{align*}
\end{definition}
Note that
\begin{equation}
\mathbb{Q}_y\left[\widetilde{A}_{t,h^{\smallfrown},\chi(t)}^{(y_0)}\right] \le \mathbb{Q}_y\left[\widetilde{A}_{t,h^{\smallfrown}}^{(y_0)}\right] =  p_t^{(y_0)}(y),\label{eq:Triv}
\end{equation}
using that under $\mathbb{Q}_y$ we have that $(X_s^{h^{\smallfrown}}+\frac{s}{t}W_t)_{s\le t}$ is a Brownian motion starting at $y$.

The proof of Theorem \ref{Sa:GirsConvH} is split into the following two propositions.
\begin{proposition}\label{Prop:Remh1}
There are $c,C>0$ such that for $y,y_0\in [-\log(t)^2,0]$ and $\lambda>0$ there is a $t_0(\lambda)>0$ such that for all $t\ge t_0(\lambda)$,
\begin{align*}
\frac{\widetilde{p}_{t,h^{\smallfrown},\chi(t)}^{(y_0)}(y)}{p_t^{(y_0)}(y)}&\le Ce^{c(\kappa-\eta+\ClogW )}\exp\left(\frac{W_t^2}{2t}\right)\left(1+t^{-\lambda}p_t^{(y_0)}(y)^{-1}\right).
\end{align*}
\end{proposition}

\begin{proposition}\label{Prop:Remh2}
There is a $C>0$ such that for $y,y_0\in [-\log(t)^2,0]$, $\lambda>0$ there is a $t_0(\lambda)>0$ such that for $t\ge t(\lambda)$, for which additionally \eqref{eq:Limgut} holds,
\[
\frac{\widetilde{p}_{t,h^{\smallfrown}}^{(y_0)}(y)-\frac{\delta}{1-\delta}\left(\sum_{j=0}^2\widetilde{p}_{t,h^{\smallfrown}}^{(y_0-1)}(y-j)\right)}{\widetilde{p}_{t,h^{\smallfrown},\chi(t)}^{(y_0)}(y)}\le 2+Ct^{-\lambda}((\ptast)^{-1}+(\ptasthat)^{-1})\widetilde{p}_{t,h^{\smallfrown},\chi(t)}^{(y_0)}(y)^{-1}.
\]
\end{proposition}

\begin{proof}[Proof of Theorem \ref{Sa:GirsConvH} assuming Propositions \ref{Prop:Remh1} and \ref{Prop:Remh2}]
Combine \eqref{eq:linLB} with Proposition \ref{Prop:Remh1} and Proposition \ref{Prop:Remh2}.
\end{proof}

\begin{lemma}\label{Lem:Bound}
There is a $C>0$ such that for $y,y_0\in [-\log(t)^2,0]$ and for all $k\ge 2$,
\begin{align*}
\mathbb{Q}_y\left[\widetilde{A}_{t,h^{\smallfrown},\chi(t)}^{(y_0)}\cap B_{t,h^{\smallfrown},\chi(t),k}\right] &\le Ce^{-\frac{(k-1)^2}{2}\left(\lfloor \chi(t)\rfloor -1\right)^{\frac{4}{3}}}.
\end{align*}
\end{lemma}
\begin{proof}[Proof of Proposition \ref{Prop:Remh1} assuming Lemma \ref{Lem:Bound}]
\end{proof}
We have that
\begin{equation}
\begin{aligned}
\widetilde{p}_{t,h^{\smallfrown},\chi(t)}^{(y_0)}(y) &= \mathbb{P}_y\left[\widetilde{A}_{t,h^{\smallfrown},\chi(t)}^{(y_0)}\right] = \mathbb{E}_{\mathbb{Q}_y}\left[(Z_t)^{-1}\mathbf{1}_{\widetilde{A}_{t,h^{\smallfrown},\chi(t)}^{(y_0)}}\right]\\
&= \sum_{k=1}^\infty \mathbb{E}_{\mathbb{Q}_y}\left[(Z_t)^{-1}\mathbf{1}_{\widetilde{A}_{t,h^{\smallfrown},\chi(t)}^{(y_0)}}\mathbf{1}_{B_{t,h^{\smallfrown},\chi(t),k}}\right] \label{eq:pthchiy}
\end{aligned}
\end{equation}

We can infer by direct calculation from \eqref{Def:Zt-1} and Definition \ref{Def:rho} that there are $c,C>0$ such that on $\widetilde{A}_{t,h^{\smallfrown},\chi(t)}^{(y_0)}\cap B_{t,h^{\smallfrown},\chi(t),k}$
\begin{align}
(Z_t)^{-1} &\le Ce^{c(\ClogW+\kappa-\eta+k)}\exp\left(\frac{W_t^2}{2t}\right). \label{eq:BoundZt}
\end{align}
In this calculation, we use assumption \eqref{Ass:WBMlog} as well as $|y|, |y_0|\le \log(t)^2$.
Plugging \eqref{eq:BoundZt} into \eqref{eq:pthchiy} yields that 
\begin{align*}
\widetilde{p}_{t,h^{\smallfrown},\chi(t)}^{(y_0)}(y)&\le Ce^{c(\kappa-\eta+\ClogW)} \exp\left(\frac{W_t^2}{2t}\right)\sum_{k=1}^\infty e^{ck}\mathbb{Q}_y\left[\widetilde{A}_{t,h^{\smallfrown},\chi(t)}^{(y_0)}\cap B_{t,h^{\smallfrown},\chi(t),k}\right]\\
&\stackrel{L. \ref{Lem:Bound}, \eqref{eq:Triv}}{\le}Ce^{c(\kappa-\eta+\ClogW)} \exp\left(\frac{W_t^2}{2t}\right)\left(p_t^{(y_0)}(y)+\sum_{k=2}^\infty e^{ck}e^{-\frac{(k-1)^2}{2}\left(\lfloor \chi(t)\rfloor -1\right)^{\frac{4}{3}}} \right)\\
&\le C\left(e^{c(\kappa-\eta+\ClogW)} \exp\left(\frac{W_t^2}{2t}\right)\left(p_t^{(y_0)}(y)+e^{-\frac{1}{2}\left(\lfloor \chi(t)\rfloor-1\right)^{\frac{4}{3}}}\right)\right)
\end{align*}
Since $e^{-1/2(\lfloor \chi(t)\rfloor-1)^{\frac{4}{3}}}t^{\lambda}\to 0$ for $t\to \infty$ and all $\lambda>0$, rearranging this yields the claim of Proposition \ref{Prop:Remh1}.

\begin{proof}[Proof of Lemma \ref{Lem:Bound}]
We have that
\begin{align}
\mathbb{Q}_y\left[B_{t,h^{\smallfrown},\chi(t),k}\right] &\le \sum_{l=\lfloor \chi(t)\rfloor}^{\left\lceil\frac{t}{2}\right\rceil} \mathbb{Q}_y\left[\min_{l-1\le s\le l} X_s^{h^{\smallfrown}}+s/t W_t-y< -(k-1)\cdot l^\frac{2}{3}\right]\notag\\
&\hskip-43pt+\sum_{l=  \left\lfloor\frac{t}{2}\right\rfloor}^{\lceil t-\chi(t)\rceil} \mathbb{Q}_y\left[\min_{l\le s\le l+1} X_s^{h^{\smallfrown}}- (t-s)/t W_t-y_0 <- (k-1)(t-l)^\frac{2}{3}, X_t\in[y_0-1,y_0]\right]. \label{lem:1}
\end{align}
First, handle the sum from $\lfloor \chi(t)\rfloor$ to $\left\lceil\frac{t}{2}\right\rceil$ in \eqref{lem:1}. Since $(X_s^{h^{\smallfrown}}+s/tW_t-y)_{s\in[0,t]}$ is a Brownian motion starting at 0 with regard to $\mathbb{Q}_y$ we have that
\begin{align}
\sum_{l=\lfloor \chi(t)\rfloor}^{\left\lceil\frac{t}{2}\right\rceil}\mathbb{Q}_y\left[\min_{l-1\le s\le l} X_s^{h^{\smallfrown}}+s/t W_t-y< -(k-1)\cdot l^\frac{2}{3}\right]&\le  \sum_{l=\lfloor \chi(t)\rfloor}^{\left\lceil\frac{t}{2}\right\rceil}4e^{-\frac{(k-1)^2l^{\frac{4}{3}}}{2}} \notag \\
&\le  Ce^{-\frac{(k-1)^2}{2}(\lfloor\chi(t)\rfloor-1)^{\frac{4}{3}}}.\label{lem:2}
\end{align}
The second sum in \eqref{lem:1} can be handled analogously. Combining \eqref{lem:1} and the last two facts proves the lemma.
\end{proof}
\paragraph{Proof of Proposition \ref{Prop:Remh2}}
We recall Definitions \ref{Def:ptast} and \ref{Def:rho}. We have that
\begin{equation}
\widetilde{p}_{t,h^{\smallfrown},\chi(t)}^{(y_0)}(y)\ge \widetilde{p}_{t,h^{\smallfrown}}^{(y_0)}(y)-\mathbb{P}_y\left[\widetilde{A}_{t,h^{\smallfrown}}^{(y_0)}\cap \widetilde{A}_{t,h^{\smallfrown},\text{start}}^c\right]-\mathbb{P}_y\left[\widetilde{A}_{t,h^{\smallfrown}}^{(y_0)}\cap \widetilde{A}_{t,h^{\smallfrown},\text{end}}^c\right]. \label{eq:chiArg}
\end{equation}
Next, we prove upper bounds for $\mathbb{P}_y\left[\widetilde{A}_{t,h^{\smallfrown}}^{(y_0)}\cap \widetilde{A}_{t,h^{\smallfrown},\text{start}}^c\right]$ and $\mathbb{P}_y\left[\widetilde{A}_{t,h^{\smallfrown}}^{(y_0)}\cap \widetilde{A}_{t,h^{\smallfrown},\text{end}}^c\right]$, we do this in detail for the first term and handle the second term by symmetry. Thus we aim to prove the following two propositions.
\begin{proposition}\label{Prop:Startchi}
There is a $C>0$ such that for $y,y_0\in [-\log(t)^2,0]$, $\lambda>0$ there is a $t_0(\lambda)>0$ such that for $t\ge t_0(\lambda)$, for which \eqref{eq:Limgut} holds,
\[
\mathbb{P}_y\left[\widetilde{A}_{t,h^{\smallfrown}}^{(y_0)}\cap\widetilde{A}_{t,h^{\smallfrown},\text{start}}^c\right] \le \delta\left(\widetilde{p}_{t,h^{\smallfrown}}^{(y_0-1)}(y)+\widetilde{p}_{t,h^{\smallfrown}}^{(y_0)}(y)\right)+Ct^{-\lambda}(\ptast)^{-1}.
\]
\end{proposition}

\begin{proposition}\label{Prop:Endchi}
There is a $C>0$ such that for $y,y_0\in [-\log(t)^2,0]$, $\lambda>0$ there is a $t_0(\lambda)>0$ such that for $t\ge t_0(\lambda)$, for which \eqref{eq:Limgut} holds,
\[\mathbb{P}_y\left[\widetilde{A}^{(y_0)}_{t,h^{\smallfrown}}\cap\widetilde{A}_{t,h^{\smallfrown},\text{end}}^c\right]\le \delta\left(\widetilde{p}_{t,h^{\smallfrown}}^{(y_0-1)}(y-2)+\widetilde{p}_{t,h^{\smallfrown}}^{(y_0-1)}(y-1)\right)+Ct^{-\lambda}(\ptasthat)^{-1}.
\]
\end{proposition}
\begin{remark}
We prove Proposition \ref{Prop:Endchi} by time reversal, thus essentially copying the proof of Proposition \ref{Prop:Startchi}.
\end{remark}
\begin{proof}[Proof of Proposition \ref{Prop:Remh2} assuming Propositions \ref{Prop:Startchi}, \ref{Prop:Endchi}]
Plugging Propositions \ref{Prop:Startchi}, \ref{Prop:Endchi} into equality \eqref{eq:chiArg} immediately yields the claim of Proposition \ref{Prop:Remh2}.
\end{proof}
We need one further definition before proceeding with the proof of Proposition \ref{Prop:Startchi}.
\begin{definition}
For $r\le \chi(t)$ set
\begin{align}
\tau &:= \inf\left\{s\le t : X_s-y\le -2\ClogW (1+s)^\frac{2}{3}+\eta\right\},\notag\\
w_r &:= y-2\ClogW (1+r)^\frac{2}{3}+W_t r/t+\eta,\label{Def:wr}\\
m_r(s) &:= W_{s+r}-h^{\smallfrown}_t(s+r). \label{Def:mr}
\end{align}
\end{definition}

The rough strategy is to bound 
\[
\mathbb{P}_y\left[\widetilde{A}_{t,h^{\smallfrown}}^{(y_0)}\cap \widetilde{A}_{t,h^{\smallfrown},\text{start}}^c\right]\le \sum \mathbb{P}_y\left[\tau\in [l-1,l]\right]\mathbb{P}_{w_r-W_t r/t}\left[\mathcal{B}_{[0,t-r], -m_r(\cdot)+(r+\cdot)/tW_t}^{J_{y_0}}(X_\cdot)\right],
\]
 use that $\tau$ has exponential tails and bound 
\[
\mathbb{P}_{w_r-W_t r/t}\left[\forall_{s\le t-r}\, X_s\le m_r(s), X_{t-r}-W_t\in [y_0-1,y_0]\right]
\] by applying Theorems \ref{Theo: LinTerm}, \ref{Theo:NebOS} and \ref{Theo: MainOS}. We break this up into a few lemmata, which we state next. After this, we will prove Proposition \ref{Prop:Startchi} assuming the lemmata. Finally, we will prove the lemmata. 

\begin{lemma}\label{Lem:chiRem1}
We have that for $t\ge e$ and $y,y_0\in [-\log(t)^2,0]$,
\begin{align*}
\mathbb{P}_y\left[\widetilde{A}_{t,h^{\smallfrown}}^{(y_0)}\cap \widetilde{A}_{t,h^{\smallfrown},\text{start}}^c\right]&\le \sum_{l=1}^{\lceil \chi(t)\rceil}\mathbb{P}_y\left[\tau\in [l-1,l]\right]\cdot\\
&\hskip2cm\cdot\max_{r\in [l-1,l]}\Bigg(e^{\frac{W_t^2}{t}-\frac{W_t^2}{2t}\frac{t-r}{t}+5\ClogW }\mathbb{P}_{w_r}\left[\mathcal{B}_{[0,t-r],-m_r}^{J_{y_0}}(X_\cdot)\right]\Bigg).
\end{align*}
\end{lemma}

\begin{lemma}\label{Lem:chiRemtau}
For $l\in\N$ and $y,y_0\in [-\log(t)^2,0]$, we have that 
\[
\mathbb{P}_y\left[\tau\in [l-1,l]\right]\le 2e^{-\frac{(2\ClogW )^2 l^{\frac{4}{3}-1}+\eta^2l^{-1}}{2}}.
\]
\end{lemma}

\begin{lemma}\label{Lem:RUB}
We have that for $t>0$ and $y,y_0\in [-\log(t)^2,0]$,
\begin{align*}
\mathbb{P}_{w_r}\left[\mathcal{B}_{[0,t-r],-m_r}^{J_{y_0}}(X_\cdot)\right]&\le \left(p_{t,h^{\smallfrown}}^{(y_0-1)}(w_r-W_r)+p_{t,h^{\smallfrown}}^{(y_0)}(w_r-W_r)\right)\mathbb{P}_{-1}\left[\mathcal{B}_{[0,r], -W_\cdot}^{J_{0}}(X_\cdot)\right]^{-1}.
\end{align*}
\end{lemma}

\begin{lemma} \label{Lem:UseStartMove}
We have that for $y,y_0\in [-\log(t)^2,0]$, $\lambda>0$ arbitrary and $t$ for which \eqref{eq:Limgut} holds,
\begin{align*}
p_{t,h^{\smallfrown}}^{(y_0)}(w_r-W_r)&\le p_{t,h^{\smallfrown}}^{(y_0)}(y)\ConstWnW^{-2}e^{4\Coh }4^{\ExpPos+2}\left(|w_r-W_r|\right)^{4\ExpPos+3\Cnull }+C e^{-5\ClogW -\ClogW ^2}t^{-\lambda}(\ptast)^{-1}.
\end{align*}
\end{lemma}

\begin{proof}[Proof of Proposition \ref{Prop:Startchi} assuming Lemmata \ref{Lem:chiRem1}--\ref{Lem:UseStartMove}]
By \eqref{eq:linLB} we have that
\begin{align}
p_{t,h^{\smallfrown}}^{(y_0)}(y)&\le e^{-\frac{W_t^2}{2t}}e^{10\ClogW }\widetilde{p}_{t,h^{\smallfrown}}^{(y_0)}(y).  \label{eq:PP71}
\end{align}
Equation \eqref{eq:PP71}, Assumption \eqref{Ass:Assdoof} and Lemmata \ref{Lem:chiRem1}, \ref{Lem:chiRemtau}, \ref{Lem:RUB}, \ref{Lem:UseStartMove} imply that for $\lambda>0$ arbitrary
\begin{align}
&\mathbb{P}_y\left[\widetilde{A}_{t,h^{\smallfrown}}^{(y_0)}\cap \widetilde{A}_{t,h^{\smallfrown},\text{start}}^c\right] \notag\\
&\le \sum_{l=1}^{\lceil \chi(t)\rceil} 2e^{-\frac{(2\ClogW )^2 l^{\frac{4}{3}-1}+\eta^2l^{-1}}{2}} \max_{r\in [l-1,l]} \Bigg(e^{\frac{W_t^2}{t}-\frac{W_t^2}{2t}\frac{t-r}{t}+5\ClogW}\cdot \notag\\
&\hskip3cm\cdot\Big(\left(p_{t,h^{\smallfrown}}^{(y_0-1)}(w_r-W_r)+p_{t,h^{\smallfrown}}^{(y_0)}(w_r-W_r)\right) \mathbb{P}_{-1}\left[\mathcal{B}_{[0,r], -W_\cdot}^{J_{0}}(X_\cdot)\right]^{-1} \Big)\Bigg)\notag\\
&\le  Ce^{-\ClogW ^2+\frac{W_t^2}{2t}-\frac{W_t^2}{2t}\frac{t-\lceil\chi(t)\rceil}{t}}t^{-\lambda}(\ptast)^{-1}+\left(\widetilde{p}_{t,h}^{(y_0-1)}(y)+\widetilde{p}_{t,h}^{(y_0)}(y)\right)\ConstWnW^{-2}e^{4\Coh }4^{\ExpPos+2}e^{15\ClogW +\ClogW ^2}\notag\\
&\qquad\cdot\left(\sum_{l=1}^{\lceil \chi(t)\rceil} 2e^{-\frac{(2\ClogW )^2 l^{\frac{4}{3}-1}+\eta^2l^{-1}}{2}}\ConstStart^{-1}(l+1)^{\ExpStart}e^{\frac{W_t^2}{2t}-\frac{W_t^2}{2t}\frac{t-l}{t}}\max_{r\in [l-1,l]}|w_r-W_r|^{4\ExpPos+3\Cnull}\right). \label{eq:Yay}
\end{align}
Using $l\le \lceil \chi(t)\rceil$ and the definition of $w_r$ in  \eqref{Def:wr} yields that
\begin{align}
\max_{r\in [l-1,l]}|w_r-W_r|&\le |y|+3\ClogW (1+l)^\frac{2}{3}+|\eta|, \label{eq:PP72}\\
 \frac{W_t^2}{2t}-\frac{W_t^2}{2t}\frac{t-l}{t} &\le \frac{W_t^2\lceil \chi(t)\rceil}{2t^2} \stackrel{\eqref{Ass:WBMlog}}{\le} \ClogW ^2.\label{eq:PP73}
\end{align}
 Plugging \eqref{eq:PP72} and \eqref{eq:PP73} into \eqref{eq:Yay} and remembering the choice of $\eta$ in \eqref{eq:eta0}  finishes the proof of Proposition \ref{Prop:Startchi}. \end{proof}

\begin{proof}[Proof of Lemma \ref{Lem:chiRem1}]
By the choice of $\kappa$ in \eqref{Def:Kappa} we have that $-\kappa(1+s)^\frac{2}{3}-h^{\smallfrown}_t(s)\le -2\ClogW   (1+s)^\frac{2}{3}$ for all $s\le \chi(t)$. Thus we have that
\begin{align*}
&\mathbb{P}_y\left[\widetilde{A}_{t,h^{\smallfrown}}^{(y_0)}\cap\widetilde{A}_{t,h^{\smallfrown},\text{start}}^c\right] \\
&= \mathbb{P}_y\left[\widetilde{A}_{t,h^{\smallfrown}}\cap\left\{\exists_{s\le \chi(t)} : X_s -y\le -\kappa(1+s)^\frac{2}{3}-h^{\smallfrown}_t(s)+\eta\right\}\right]\\
&\le \mathbb{P}_y\left[\widetilde{A}_{t,h^{\smallfrown}}\cap \left\{\exists_{s\le \chi(t)}: X_s-y\le -2\ClogW (1+s)^\frac{2}{3}+\eta\right\}\right]\\
&\le \sum_{l=1}^{\lceil \chi(t)\rceil}\mathbb{P}_y\left[\tau\in [l-1,l]\right]\cdot\max_{r\in [l-1,l]}\mathbb{P}_{w_r-W_t r/t}\left[\mathcal{B}_{[0,t-r],-m_r}^{J_{y_0}}(X_s+W_t(s+r)/t)\right]\\
&=\sum_{l=1}^{\lceil \chi(t)\rceil}\mathbb{P}_y\left[\tau\in [l-1,l]\right]\cdot\max_{r\in [l-1,l]}\mathbb{P}_{w_r}\left[\mathcal{B}_{[0,t-r],-m_r}^{J_{y_0}}(X_s+W_t s/t)\right]\\
&\stackrel{T.\ \ref{Theo: LinTerm}, \eqref{Ass:WBMlog}}{\le} \sum_{l=1}^{\lceil \chi(t)\rceil}\mathbb{P}_y\left[\tau\in [l-1,l]\right]\max_{r\in [l-1,l]}e^{\frac{W_t^2}{t}-\frac{W_t^2}{2t}\frac{t-r}{t}+5\ClogW }\mathbb{P}_{w_r}\left[\mathcal{B}_{[0,t-r],-m_r}^{J_{y_0}}(X_\cdot)\right],
\end{align*}
where in the last step we applied Theorem \ref{Theo: LinTerm} for time $t-r$ and $c_{t-r} = \frac{W_t}{t}$.
\end{proof}

\begin{proof}[Proof of Lemma \ref{Lem:chiRemtau}]
We have that
\begin{align*}
\mathbb{P}_y\left[\tau\in [l-1,l]\right] &\le \mathbb{P}_y\left[\min_{s\in [0,l]} X_s-y \le -2\ClogW l^\frac{2}{3}+\eta\right]\\
&= 2\mathbb{P}_y\left[(X_l-y)/\sqrt{l}\le-2\ClogW  l^{\frac{2}{3}-1/2}+\eta l^{-1/2}\right]\le 2e^{-\frac{(2\ClogW )^2 l^{\frac{4}{3}-1}+\eta^2l^{-1}}{2}}.\qedhere
\end{align*}
\end{proof}

\begin{proof}[Proof of Lemma \ref{Lem:RUB}]
Applying \eqref{eq:BarrSplit2} from Lemma \ref{Lem:BarrSplit} for $z_1 = w_r-W_r$, $z_2 =y_0$, $I = [0,t]$, $t_0 = r$, $Z_{\cdot} = X_{\cdot}+h^{\smallfrown}_t(\cdot)-W_{\cdot}$ and $x_0 = w_r+h_t^{\smallfrown}(r)-W_r$  yields that
\begin{align*}
p_{t,h^{\smallfrown}}^{(y_0-1)}(w_r-W_r)+p_{t,h^{\smallfrown}}^{(y_0)}(w_r-W_r) &\ge \mathbb{P}_{w_r-W_r}\left[\mathcal{B}_{[0,r],h^{\smallfrown}_t(\cdot)-W_\cdot}^{J_{w_r+h^{\smallfrown}_t(r)-W_r}}(X_\cdot)\right]\cdot\\
&\hskip3.5cm\cdot \mathbb{P}_{w_r}\left[\mathcal{B}_{[r,t],h^{\smallfrown}_t(\cdot)-W_\cdot}^{J_{y_0}}(X_\cdot-X_r)\right].
\end{align*}
Using $(X_s-X_r)_{s\ge r} \stackrel{d}{=} (X_s)_{s\ge 0}$ and the definition of $m_r$ in \eqref{Def:mr} yields that 
\[
 \mathbb{P}_{w_r}\left[\mathcal{B}_{[r,t],h^{\smallfrown}_t(\cdot)-W_\cdot}^{J_{y_0}}(X_\cdot-X_r)\right] = \mathbb{P}_{w_r}\left[\mathcal{B}_{[0,t-r],-m_r}^{J_{y_0}}(X_\cdot)\right]
\]
and thus it remains to prove that
\begin{equation}
\mathbb{P}_{w_r-W_r}\left[\mathcal{B}_{[0,r],h^{\smallfrown}_t(\cdot)-W_\cdot}^{J_{w_r+h^{\smallfrown}_t(r)-W_r}}(X_\cdot)\right] \ge \mathbb{P}_{-1}\left[\mathcal{B}_{[0,r], -W_\cdot}^{J_{0}}(X_\cdot)\right]. \label{eq:B0rhW}
\end{equation}
Equation \eqref{eq:B0rhW} is an immediate consequence of the facts that on the one hand $w_r-W_r\le -1$ by the definition of $\ClogW $ and $w_r$ in \eqref{Ass:WBMlog} and \eqref{Def:wr} and on the other hand $h^{\smallfrown}_t(s)\le 0$ for all $s\le 0$.
\end{proof}

\begin{proof}[Proof of Lemma \ref{Lem:UseStartMove}]
Since $r\le \chi(t)$, recalling \eqref{Def:wr} and applying \eqref{Ass:WBMlog} yields that 
\begin{align*}
|w_r-W_r|&\le \log(t)^{2}+(\ClogW )(1+\log(t)^{1.5})+\log(t)^{1.5}\ClogW \le 2\log(t)^{2}
\end{align*}
by assumption \eqref{eq:Limgut}. Furthermore, using the definition of $\ClogW $ and $\eta<0$, yields that 
\[
w_r-W_r = y-2\ClogW (1+r)^{\frac{2}{3}}+W_t r/t+\eta-W_r\le y.
\]

Thus we can apply Theorems \ref{Theo:NebOS} and \ref{Theo: MainOS} to get that there is a $t_0(\lambda)$ such that for $t\ge t_0(\lambda)$
\begin{align*}
p_{t,h^{\smallfrown}}^{(y_0)}(w_r-W_r)&\le p_{t,h^{\smallfrown}}^{(y_0)}(y)\ConstWnW^{-2}e^{4\Coh }4^{\ExpPos+2}\left(|w_r-W_r|\right)^{4\ExpPos+3\Cnull }+C\ConstWnW^{-1}\log(t)^{4\ExpPos}t^{-\lambda}p_{t,h^{\smallfrown}}^{(y_0)}(y)^{-1},
\end{align*}
which is what we wanted to prove.
\end{proof}

\begin{proof}[Proof of Proposition \ref{Prop:Endchi}]
Define $Y_s := X_{t-s}-X_t$. Using the symmetry of $h_t^{\smallfrown}$ 
\begin{align}
&\mathbb{P}_y\left[\widetilde{A}_{t,h^{\smallfrown}}^{(y_0)}\cap \widetilde{A}_{t,h^{\smallfrown},\text{end}}^c\right]\notag\\
&\le \mathbb{P}\Big[\mathcal{B}_{[0,t], -\widehat{W}_\cdot}^{y_0-1,J_{y}}\left(Y_\cdot+h^{\smallfrown}_t(s)+\frac{s}{t}\widehat{W}_t\right), \exists_{s\in [0,\chi(t)]} Y_s\le-\kappa(1+s)^\frac{2}{3}-h^{\smallfrown}_t(s)+y_0 +\eta \Big] \label{eq:Zeitumkehr1}
\end{align}
Furthermore, we have that
\begin{align}
\widetilde{p}_{t,h^{\smallfrown}}^{(y_0)}(y)&= \mathbb{P}\left[\mathcal{B}_{[0,t],-W_\cdot}^{y,J_{y_0}}\left(X_s+h^{\smallfrown}_t(s)+\frac{s}{t}W_t\right)\right]\ge \mathbb{P}\left[\mathcal{B}_{[0,t], -\widehat{W}_\cdot}^{y_0,J_{y+1}}\left(Y_\cdot+h^{\smallfrown}_t(s)+\frac{s}{t}\widehat{W}_t\right)\right] \label{eq:Zeitumkehr2}
\end{align}
We notice that $(Y_s)_{s\in [0,t]}\stackrel{d}{=} (X_s)_{s\in [0,t]}$. Thus rerunning the argument proving Proposition \ref{Prop:Startchi}, while considering equations \eqref{eq:Zeitumkehr1} and \eqref{eq:Zeitumkehr2}, finishes the proof.
\end{proof}

\subsection{Bounds Used in the Lower Bound of the Right Tail} \label{Sec:htLBRT}
In this section we consider $h_t^{\smallsmile} : [0,t]\to \R, s\mapsto \left(((1+s)\wedge (1+t-s))^{1/6}-1\right)$ and $y\in[-\log(t)^{2},0]$, $y_0\in [-\log(t)^{2},0]$.

We prove the following theorem.

\begin{theorem}\label{Sa:GirsConcH}
For $\eta$, $\kappa$ defined analogously to Theorem \ref{Sa:GirsConvH}, but using $h^{\smallsmile}$ instead of $h^{\smallfrown}$, we have that there are $c,C>0$ such that for arbitrary $\lambda>0$ there is a $t_0(\lambda)>0$ such that for $t\ge t_0(\lambda)$ with \eqref{eq:Limgut} and for $y,y_0 \in [-\log(t)^2,0]$,
\begin{align*}
&\frac{p_{t,h^{\smallsmile}}^{(y_0)}(y)}{p_t^{(y_0)}(y)-\frac{\delta}{1-\delta}\sum_{j=0}^2 p_t^{(y_0-1)}(y-j)} \\
&\ge Ce^{-c(\ClogW+\kappa-\eta) }\left(\frac{1}{2}-Ct^{-\lambda}\left((\ptast)^{-1}+(\ptasthat)^{-1}+1\right)p_t^{(y_0)}(y)^{-1}\right).
\end{align*}
\end{theorem}

\begin{proof}[Sketch of the proof of Theorem \ref{Sa:GirsConcH}]
We recall \eqref{Def:Zt-1}. Furthermore, define
\begin{align}
B_{t,h^{\smallsmile},\chi(t),2} &:= \left\{\min_{\chi(t)\le s\le t/2} \frac{X_s^{h^{\smallsmile}}+s/tW_t-y}{(1+s)^\frac{2}{3}}\ge -2\right\}\notag\\
&\hskip3.5cm\cap\left\{\min_{t/2\le s\le t-\chi(t)}\frac{X_s^{h^{\smallsmile}}-(t-s)/tW_t-y_0}{(1+t-s)^\frac{2}{3}}\ge -2\right\}. \label{eq:Athsideline}
\end{align}
As in \eqref{eq:BoundZt} we can infer from \eqref{Def:Zt-1} and \eqref{eq:Athsideline}  that there are $C,c>0$ such that
\begin{align}
\widetilde{p}_{t,h^{\smallsmile}}^{(y_0)}(y) &= \mathbb{E}_{\mathbb{Q}_y}\left[(Z_t)^{-1}\mathbf{1}_{\widetilde{A}_{t,h^{\smallsmile}}}\right] \ge \mathbb{E}_{\mathbb{Q}_y}\left[(Z_t)^{-1}\mathbf{1}_{\widetilde{A}_{t,h^{\smallsmile},\chi(t)}}\mathbf{1}_{B_{t,h^{\smallsmile},\chi(t),2}}\right] \notag\\
&\ge Ce^{-c(\ClogW+\kappa-\eta)}\exp\left(\frac{W_t^2}{2t}\right)\mathbb{Q}_y\left[\widetilde{A}_{t,h^{\smallsmile},\chi(t)}\cap{B_{t,h^{\smallsmile},\chi(t),2}}\right]. \label{eq:pttildlb}
\end{align}
We have that $\mathbb{Q}_y\left[\widetilde{A}_{t,h^{\smallsmile},\chi(t)}\cap{B_{t,h^{\smallsmile},\chi(t),2}}\right]\ge \mathbb{Q}_y\left[\widetilde{A}_{t,h^{\smallsmile},\chi(t)}\right]-\mathbb{Q}_y\left[\widetilde{A}_{t,h^{\smallsmile},\chi(t)}\cap B_{t,h^{\smallsmile},\chi(t),2}^c\right]$. As in Lemma \ref{Lem:Bound} we have that
\begin{align}
\mathbb{Q}_y\left[\widetilde{A}_{t,h^{\smallsmile},\chi(t)}\cap B_{t,h^{\smallsmile},\chi(t),2}^c\right]&\le Ct^{-\lambda} \label{eq:AthchiUB}
\end{align}
for $t$ big enough, depending on $\lambda$, but not on $W$.

Equations \eqref{eq:pttildlb} and \eqref{eq:AthchiUB} imply that for $t$ big enough depending on $\lambda$, but not on $W$,
\begin{equation}
\widetilde{p}_{t,h^{\smallsmile}}^{(y_0)}(y)\ge  Ce^{-c(\ClogW+\kappa-\eta)}\exp\left(\frac{W_t^2}{2t}\right)\left(\mathbb{Q}_y\left[\widetilde{A}_{t,h^{\smallsmile},\chi(t)}\right]-t^{-\lambda}\right). \label{eq:ptildQyBound}
\end{equation}

Furthermore, we have that
\begin{align*}
\mathbb{Q}_y\left[\widetilde{A}_{t,h^{\smallsmile},\chi(t)}\right] &= \mathbb{Q}_y\left[\widetilde{A}_{t,h^{\smallsmile}}^{(y_0)}\right]-\mathbb{Q}_y\left[\widetilde{A}_{t,h^{\smallsmile}}^{(y_0)}\cap\widetilde{A}_{t,h^{\smallsmile},\text{start}}^c\right]-\mathbb{Q}_y\left[\widetilde{A}_{t,h^{\smallsmile}}^{(y_0)}\cap\widetilde{A}_{t,h^{\smallsmile},\text{end}}^c\right]\\
&= p_t^{(y_0)}(y)-\mathbb{Q}_y\left[\widetilde{A}_{t,h^{\smallsmile}}^{(y_0)}\cap\widetilde{A}_{t,h^{\smallsmile},\text{start}}^c\right]-\mathbb{Q}_y\left[\widetilde{A}_{t,h^{\smallsmile}}^{(y_0)}\cap\widetilde{A}_{t,h^{\smallsmile},\text{end}}^c\right].
\end{align*}

From here one proceeds as in the proof of Proposition \ref{Prop:Remh2}. We need to calculate the tails of $\tau$ with regard to $\mathbb{Q}_y$ instead of $\mathbb{P}_y$. This change results in the calculation
\begin{align*}
\mathbb{Q}_y\left[\tau\in [l-1,l]\right] &\le \mathbb{Q}_y\left[\min_{s\in [0,l]} X_s^{h^{\smallsmile}}+s/t W_t-y\le -2\ClogW (1+l)^\frac{2}{3}+h^{\smallsmile}_t(l)+l^\frac{2}{3}\ClogW +\eta\right]\\
&\le \mathbb{Q}_y\left[\min_{s\in [0,l]}X_s^{h^{\smallsmile}}+s/tW_t-y\le -\ClogW  (1+l)^\frac{2}{3}+\eta\right]\\
&\le 2e^{-\frac{(\ClogW )^2l^{\frac{4}{3}-1}+\eta^2l^{-1}}{2}}.
\end{align*} 

The rest of the proof of Proposition \ref{Prop:Remh2} needs only very minor changes\footnote{Since we work with regard to $\mathbb{Q}_y$ our Brownian motion is $(X_s+h^{\smallsmile}_t(s)+s/t W_t)_{s\in [0,t]}$. Thus the analogue to $m_r$ of \eqref{Def:mr} won't have an $h^{\smallsmile}_t$ term, which only simplifies the situation.} to give that

\begin{align*}
\mathbb{Q}_y\left[\widetilde{A}_{t,h^{\smallsmile}}^{(y_0)}\cap\widetilde{A}_{t,h^{\smallsmile},\text{start}}^c\right]&\le \delta\left(p_t^{(y_0-1)}(y)+p_t^{(y_0)}(y)\right)+Ct^{-\lambda}(\ptast)^{-1}
\end{align*}
as well as
\[
\mathbb{Q}_y\left[\widetilde{A}_{t,h^{\smallsmile}}^{(y_0)}\cap\widetilde{A}_{t,h^{\smallsmile},\text{end}}^c\right]\le \delta\left(p_t^{(y_0-1)}(y-2)+p_t^{(y_0-1)}(y-1)\right)+Ct^{-\lambda}(\ptasthat)^{-1}
\]
which yields that
\begin{align}
\frac{\mathbb{Q}_y\left[\widetilde{A}_{t,h^{\smallsmile},\chi(t)}^{(y_0)}\right]}{p_t^{(y_0)}(y)-\frac{\delta}{1-\delta}\sum_{j=0}^2 p_t^{(y_0-1)}(y-j)}\ge \frac{1}{2}-Ct^{-\lambda}\left((\ptast)^{-1}+(\ptasthat)^{-1}\right)p_t^{(y_0)}(y)^{-1}. \label{eq:yptyBound}
\end{align}

Combining \eqref{eq:linUB}, \eqref{eq:ptildQyBound} and \eqref{eq:yptyBound} and writing
\begin{align*}
\frac{p_{t,h^{\smallsmile}}^{(y_0)}(y)}{p_t^{(y_0)}(y)-\frac{\delta}{1-\delta}\sum_{j=0}^2 p_t^{(y_0-1)}(y-j)}&= \frac{p_{t,h^{\smallsmile}}^{(y_0)}(y)}{\widetilde{p}_{t,h^{\smallsmile}}^{(y_0)}(y)}\frac{\widetilde{p}_{t,h^{\smallsmile}}^{(y_0)}(y)}{\mathbb{Q}_y[\widetilde{A}_{t,h^{\smallsmile},\chi(t)}]}\frac{\mathbb{Q}_y[\widetilde{A}_{t,h^{\smallsmile},\chi(t)}]}{p_t^{(y_0)}(y)-\frac{\delta}{1-\delta}\sum_{j=0}^2 p_t^{(y_0-1)}(y-j)}
\end{align*}
gives Theorem \ref{Sa:GirsConcH}.
\end{proof}

\section{Crude lower bound on $p_{t,h}^{(y_0)}(y)$}\label{Sec:Crude}
\begin{figure}
\centering
\includegraphics[scale=0.7]{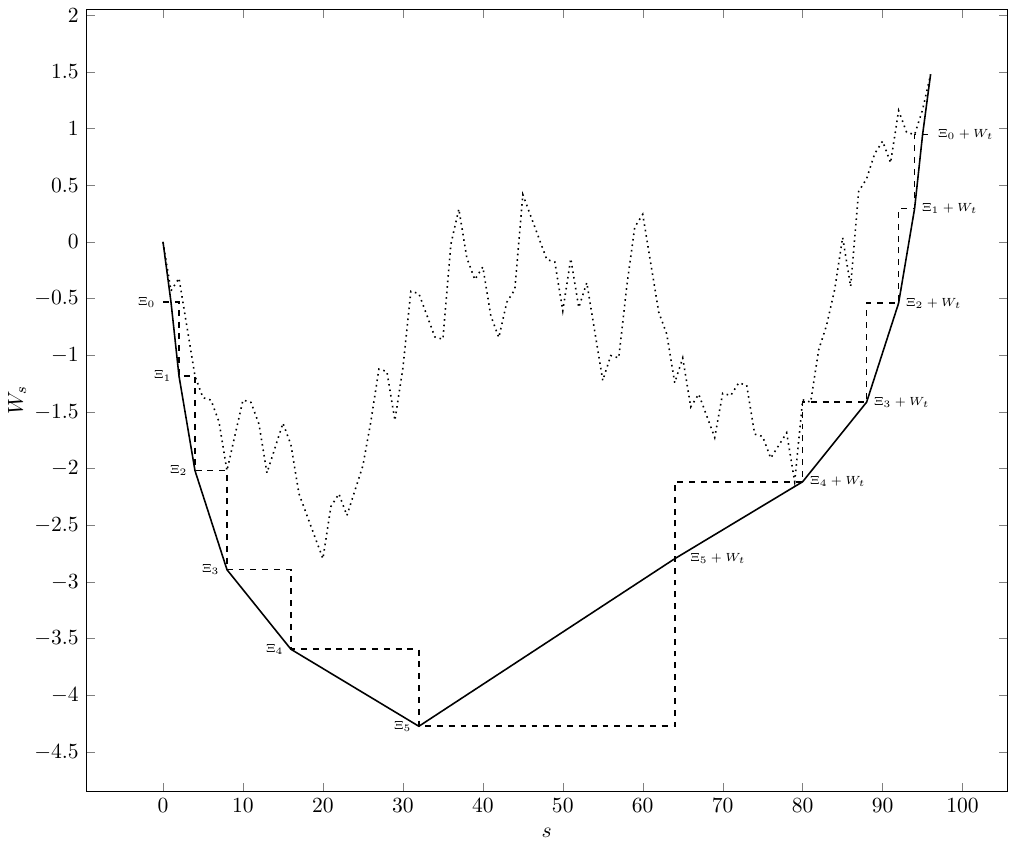}
\caption{The event in \eqref{eq:Sec61}. Drawn are $(W_s)_{0\le s\le t}$ (dotted line), $\Xi_j(W)$, $j\in\{0,\dots, 5\}$ (dashed lines) and $g_t$ (solid piecewise linear curve) with parameters $h=0$, $t=96$, $t_1 =32$, $t_2 =64$, $k_1=5$. We note that it is less likely to stay below $(g_s)_{s\le t}$ than below $(W_s)_{s\le t}$. }
\label{fig:Defgt}
\end{figure}

The aim of this section is to provide a lower bound on $\smash{p_{t,h}^{(y_0)}(y)}$ of Definition \ref{Def:pthAllg}. This will be used to deal with the $p_{t,h}^{(y_0)}(x)^{-1}$ terms occurring in Theorems \ref{Theo: MainOS}, \ref{Sa:GirsConvH} and \ref{Sa:GirsConcH} and additionally be the main ingredient in the proof of Lemma \ref{Lem:Boundpnastweak}. We recall Definition \ref{Def:pthAllg}. The idea of the proof is to replace the event $\{\forall_{s\le t}\, X_s\le W_s-h_t(s)\}$ with $\{\forall_{s\le t} X_s \le g_t(s)\}$ with a $g_t$ for which we can apply Girsanov's theorem, to reduce the situation to a ballot theorem for Brownian motion. The $g_t(s)$ will roughly be a piecewise linear approximation of the running minimum of $W_s-h_t(s)$. The next definition defines $g_t$ precisely as well as some related quantities, see Figure \ref{fig:Defgt}.

\begin{definition}\label{Def:Delta}
Fix $t>0$. Set $k_1(t) := \lfloor\log_2(t/3)\rfloor$, $t_1(t) := 2^{k_1}$, $t_2(t) := t-t_1$, the dependence on $t$ will be omitted from notation. For $j\le k_1-1$ define
\begin{align*}
\Xi_j^{\text{start}}(W) &:= \min_{s\in [0,2^{j+1}]} W_s,\\
\Xi_j^{\text{end}}(W) &:= \min_{s\in [t-2^{j+1}, t]} W_s-W_t,\\
\Xi_j(W) &:= \min\{\Xi_j^{\text{start}}(W), \Xi_j^{\text{end}}(W)\},\\
\Xi_{k_1}(W)&:= \min\{\min_{s\in [0,t]} W_s, \min_{s\in [0,t]} W_s-W_t,0\},\\
\Delta_j(W) &:= \frac{\Xi_{j+1}(W)-\Xi_j(W)}{2^j}.
\end{align*}
For $j\in\{1,\dots, k_1-1\}$ define $\Xi_j(h)$ (and analogous variables) as above, substituting $-h_t(s)$ for $W_s$, $s\in[0,t]$. Furthermore, choose $\Xi_0(h)$ such that $s\Xi_0(h)\le -h_t(s)$ for all $s\in [0,1]$ and $(t-s)\Xi_0(h)\le -h_t(s)$ for all $s\in [t-1,t]$.

Set $\Xi_j := \Xi_j(W)+\Xi_j(h)$ and $\Delta_j = \Delta_j(W)+\Delta_j(h)$. Finally, set
\[
g_t(s) = \begin{cases}
s\Xi_0,&\ s\in [0,1),\\
\Xi_j+(s-2^j)\Delta_j,&\ s \in [2^j, 2^{j+1}]\\
\Xi_{k_1}+(s-t_1)/(t_2-t_1)(W_t-h_t(t)),&\ s\in [t_1,t_2],\\
\Xi_{j}+(t-2^j-s)\Delta_j+W_t-h_t(t),&\ s\in [t-2^{j+1}, t-2^j],\\
(t-s)\Xi_0+W_t,&\ s\in (t-1,t].
\end{cases}
\]
\end{definition}

\begin{theorem}\label{Theo:CrudeLB}
For $W$ piecewise linear and $y,y_0\in [-\sqrt{t},0]$, there is a $\gamma_0>0$ deterministic such that for $t\ge 2$,
\begin{align*}
&\mathbb{P}_y\left[\mathcal{B}_{[0,t],h_t(\cdot)-W_\cdot}^{J_{y_0}}(X_\cdot)\right]\\
&\ge t^{-\gamma_0-\log(t)^{-1}\left(\sum_{j=0}^{k_1\hskip-0.035cm-\hskip-0.035cm1}\Delta_j^22^j+2(1-2^{3/2})\sum_{j=0}^{k_1\hskip-0.035cm-\hskip-0.035cm1}\Delta_j2^{j/2}+\Xi_0^2-4\Xi_0+3\frac{W_t^2+h_t(t)^2}{t_2-t_1}+(\sqrt{t_1}+|y-y_0|)\frac{|W_t|+|h_t(t)|}{t_2-t_1}\right)}. 
\end{align*}
\end{theorem}
\begin{proof}
As illustrated in Figure \ref{fig:Defgt} we have that by definition $g_t(s)\le W_s$ for all $s\in [0,t]$ as well as $g_t(t) = W_t = W_t-h_t(t)$ and thus
\begin{align}
\mathbb{P}_y\left[\mathcal{B}_{[0,t],h_t(\cdot)-W_\cdot}^{J_{y_0}}(X_\cdot)\right]&\ge \mathbb{P}_{y}\left[\forall_{s\le t}\, X_s \le g_t(s), X_t-g_t(t)\in [y_0-1,y_0] \right]. \label{eq:Sec61}
\end{align}
We apply Girsanov to remove the $g_t$. For controlling the Radon-Nikodym derivative appearing in the application of Girsanov, we need control over $X_{2^{j+1}}-X_{2^j}$, $j\in\{0,\dots, k_1-1\}$. Thus define
\begin{align*}
A_j^{(1)} &:= \left\{ X_{2^{j+1}}-g_t(2^{j+1}) \in [y-2\sqrt{2^{j+1}},y-\sqrt{2^{j+1}}]\right\},\\
A_j^{(2)} &:= \Big\{X_{t-2^{j+1}}-g_t(t-2^{j+1})\in [y_0-2\sqrt{2^{j+1}},y_0-\sqrt{2^{j+1}}]\Big\}.
\end{align*}
Furthermore, set
\begin{align*}
A_{\text{start}} &:= \left\{y-2\le X_1-g_t(1)\le y-1\right\}\\
A_{\text{end}} &:= \left\{y_0-2\le X_{t-1}-g_t(t-1) \le y_0-1 \right\}.
\end{align*}
By monotonicity and \eqref{eq:Sec61}
\begin{align*}
&\mathbb{P}_y\left[\mathcal{B}_{[0,t],h_t(\cdot)-W_\cdot}^{J_{y_0}}(X_\cdot)\right]\\
&\ge \mathbb{P}_{y}\left[A_{\text{start}}, \bigcap_{i\in\{1,2\}}\bigcap_{j=0}^{k_1-1} A_j^{(i)}, A_{\text{end}},\forall_{s\le t}\, X_s \le g_t(s), X_t-g_t(t)\in [y_0-1,y_0]  \right]:= \mathbb{P}_{y}[A(y_0)].
\end{align*}

Consider $\mathbb{Q}_{y}$ with
\begin{align*}
Z_t &:= \frac{d\mathbb{Q}_{y}}{d\mathbb{P}_{y}} = \exp\left(\int_0^t g_t'(s)\;\mathrm{d}X_s-\frac{1}{2}\int_0^t (g_t'(s))^2\;\mathrm{d}s\right)\\
&= \exp\Bigg(\sum_{j=0}^{k_1-1} \big(\Delta_j\left((X_{2^{j+1}}-X_{2^j})-(X_{t-2^j}-X_{t-2^{j+1}})\right) \big)- \sum_{j=0}^{k_1-1} \Delta_j^2 2^j-\Xi_0^2\\\
&\qquad+\Xi_0\cdot\left(X_1-X_0-(X_t-X_{t-1})+(X_{t_2}-X_{t_1})\frac{W_t-h_t(t)}{t_2-t_1}-\frac{(W_t-h_t(t))^2}{2(t_2-t_1)}\right)\Bigg).
\end{align*}

On $A(y_0)$ we have that
\begin{align*}
|X_{2^{j+1}}-X_{2^j}| &\le -2^j\Delta_j+2^{\frac{j}{2}}\left(2^{\frac{3}{2}}-1\right),\\
|X_{t-2^j}-X_{t-2^{j+1}}|&\le -2^j\Delta_j+2^{\frac{j}{2}}\left(1-2^{\frac{3}{2}}-1\right),\\
X_1-X_0 &\ge \Xi_0-2,\\
X_t-X_{t-1}&\le  2-\Xi_0,\\
|X_{t_2}-X_{t_1}| &\le \sqrt{t_1}+|W_t-h_t(t)|+|y-y_0|.
\end{align*}

Thus we have
\begin{align*}
\mathbb{P}_y[A(y_0)] &= \mathbb{Q}_y[Z_t^{-1};A(y_0)]\\
&\ge \exp\Bigg(\begin{aligned}[t]&-\sum_{j=0}^{k_1-1}\Delta_j^22^j-\left(2-2^\frac{5}{2}\right)\sum_{j=0}^{k_1-1} \Delta_j 2^{\frac{j}{2}}-\Xi_0^2+4\Xi_0\\
&-\frac{3(W_t-h_t(t))^2}{2(t_2-t_1)}-\frac{(\sqrt{t_1}+|y-y_0|)|W_t-h_t(t)|}{t_2-t_1}\Bigg)\cdot\mathbb{Q}_y[A(y_0)].
\end{aligned}
\end{align*}
By the Girsanov theorem we have that $(X_s-g_t(s))_{s\in[0,t]}$ is a Brownian motion with regard to $\mathbb{Q}_y$. This implies that there is a constant $\gamma_0>0$ such that for $t\ge2$ arbitrary $\mathbb{Q}_y[A(y_0)]\ge t^{-\gamma_0}$. Thus we have that
\begin{align*}
&\mathbb{P}_y[A(y_0)]\\
&\ge t^{-\log(t)^{-1}\left(\sum_{j=0}^{k_1-1}\Delta_j^22^j+2(1-2^{3/2})\sum_{j=0}^{k_1-1}\Delta_j2^{j/2}+\Xi_0^2-4\Xi_0+3\frac{W_t^2+h_t(t)^2}{t_2-t_1}+(\sqrt{t_1}+|y-y_0|)\frac{|W_t|+|h_t(t)|}{t_2-t_1}\right)}\\
&\qquad\hskip11.5cm\cdot\mathbb{Q}_y[A(y_0)]\\
&\ge  t^{-\gamma_0-\log(t)^{-1}\left(\sum_{j=0}^{k_1\hskip-0.035cm-\hskip-0.035cm1}\Delta_j^22^j+2(1-2^{3/2})\sum_{j=0}^{k_1\hskip-0.035cm-\hskip-0.035cm1}\Delta_j2^{j/2}+\Xi_0^2-4\Xi_0+3\frac{W_t^2+h_t(t)^2}{t_2-t_1}+(\sqrt{t_1}+|y-y_0|)\frac{|W_t|+|h_t(t)|}{t_2-t_1}\right)}.\qedhere
\end{align*}
\end{proof}

\section{Control for the Constants Introduced in Section \ref{Subsec:DefSum}, Definition \ref{Def:ptast} and the Proof of Lemma \ref{Lem:Boundpnastweak}}\label{Sec:Constants}

We recall Definitions \ref{Def:Barr} and \ref{Def:pn}. From now on $W$ will be as in Definition \ref{Def:pn}. For $t>0$ define $\widehat{W}$ via $\widehat{W}_s:= W_{t-s}-W_t$, where the dependence on $t$ is not reflected in the notation. Furthermore, from now on we will mostly care about $h$ given by $h_t:[0,t]\to \R$, $s\mapsto \pm((1+s)\wedge(1+t-s)^{1/6}-1)$. We note here, that for these $h$ the $C_1(h)$ from \eqref{Ass:WachsH} is finite. We recall Section \ref{Subsec:DefSum}, we use $\widehat{C}$ to denote variants of constants using $\widehat{W}$ instead of $W$, for example we set $\widehat{\mathcal{C}}_{\log} := \mathcal{C}_{\log}(\widehat{W},t)$. We point out that in \eqref{Ass:QuenchedWall} the $\mathbb{P}$ will from now on be denoted by $\mathbb{P}_{\mathcal{L}}$ since the probability is taken given $W$, i.e.\@ given the environment.

Consider $f:\R_{\ge0}\to \R$. For $a,b> 0$ set
\[
\gamma_{a,b}(W,f,t) := \inf\left\{\gamma>0 : \forall_{t\ge s\ge 2}\, \mathbb{P}_{\mathcal{L}}\left[\mathcal{B}_{[0,s],f(\cdot)-W_\cdot}^{-1, [-b\sqrt{s},-a\sqrt{s}]}(B_\cdot)\right]\ge s^{-\gamma}\right\}
\]
and $\gamma_{a,b}(\widehat{W},f,t)$ analogous. 

\begin{lemma} \label{Lem:BoundGam}
Let $\beta\in \{\pm1\}$. We consider $(h_t)_{t\ge0}$ given by $h_t : [0,t]\to\R, s\mapsto \beta (((1+s)\wedge (1+t-s))^{1/6}-1)$, $b>a\ge 0$. The families of random variables $(\gamma(W,h,n))_{n\in\N}$, $(\gamma(\widehat{W},h,n))_{n\in\N}$ are tight.
\end{lemma}
\begin{proof}
On $[2,n/2]$ we have that $h_n(s) = \beta \left((1+s)^{1/6}-1\right) := f(s)$. Thus it suffices to prove  that $(\gamma_{a,b}(\widetilde{W},f,n))_{n\in\N}$, $\widetilde{W}\in \{W,\widehat{W}\}$ are tight. Both of these are proved analogously, we will write the proof for $\widetilde{W} = W$ and note that it goes through as written when replacing $W$ by $\widehat{W}$.

 By Theorem 1.10 in \cite{MM_BMaboveQRW} we know that there exists $\gamma>0$ independent of $a,b$ such that $\mathbb{P}$-almost surely
\begin{align*}
\lim_{n\to \infty}\frac{\log\left(\mathbb{P}_{\mathcal{L}}\left[\mathcal{B}_{\{0,\dots,n\}, f(\cdot)-W_\cdot}^{-1,[-b\sqrt{n},-a\sqrt{n}]}(B_\cdot)\right]\right)}{\log(n)} = -\gamma.
\end{align*}
in particular  $\mathbb{P}$-a.s.
\begin{equation}
\widetilde{\gamma} := \inf_{n\in\N} \frac{\log\left(\mathbb{P}_{\mathcal{L}}\left[\mathcal{B}_{\{0,\dots,n\}, f(\cdot)-W_\cdot}^{-1,[-b\sqrt{n},-a\sqrt{n}]}(B_\cdot)\right]\right)}{\log(n)}>-\infty. \label{eq:gammatild}
\end{equation}

We claim that there is a $c>0$ such that for all $n\ge 2$ $\mathbb{P}$-a.s.
\begin{align}
\frac{\mathbb{P}_{\mathcal{L}}\left[\mathcal{B}_{[0,n],f(\cdot)-W_\cdot}^{-1,[-b\sqrt{n},-a\sqrt{n}]}(B_\cdot)\right]}{\mathbb{P}_{\mathcal{L}}\left[\mathcal{B}_{\{0,\cdots, n\}, 2^{\beta}f(\cdot)-W_\cdot}^{-1, [-b\sqrt{n}, -a\sqrt{n}]}(B_\cdot)\right]}\ge c.\label{eq:DiskrtoCont}
\end{align} 
Combining \eqref{eq:gammatild} with \eqref{eq:DiskrtoCont} allows us to conclude the statement of the lemma, it remains to prove that \eqref{eq:DiskrtoCont} holds. For this purpose, we define
\begin{align*}
A_1 &:= \mathcal{B}_{[0,n], f(\cdot)-W_\cdot}^{-1, [-b\sqrt{n},-a\sqrt{n}]}(B_\cdot),\quad A_2 := \mathcal{B}_{\{0,\cdots, n\}, 2^{\beta}f(\cdot)-W_\cdot}^{-1, [-b\sqrt{n}, -a\sqrt{n}]}(B_\cdot).
\end{align*}
We need to bound $\mathbb{P}_{\mathcal{L}}[A_1]/\mathbb{P}_{\mathcal{L}}[A_2]$. In essence this calculation boils down to the fact, that if at times $\{0,\dots,n\}$ $B$ is below a (curved) barrier, then it is unlikely that $B$ crosses a multiple of the barrier in-between these times. We have that
\begin{align}
\mathbb{P}_{\mathcal{L}}\left[A_1\right] &=\mathbb{E}_{\mathcal{L}}\left[\mathbb{P}_{\mathcal{L}}\left[A_1| (B_1,\dots,  B_{n})\right]\right]\ge \mathbb{E}_{\mathcal{L}}\left[\mathbf{1}_{A_2}\mathbb{P}_{\mathcal{L}}\left[A_1| (B_1,\dots, B_{n})\right]\right] \notag\\
&= \mathbb{E}_{\mathcal{L}}\left[\mathbf{1}_{A_2}\prod_{j=1}^{n} \mathbb{P}_{\mathcal{L}}\left[\mathcal{B}_{[j-1,j],f(\cdot)-W_\cdot}^{-1}(B_\cdot)| (B_{j-1},B_j)\right]\right] . \label{eq:A1toA2}
\end{align}

To get a lower bound on $\mathbf{1}_{A_2}\prod_{j=1}^{n} \mathbb{P}_{\mathcal{L}}\left[\mathcal{B}_{[j-1,j],f(\cdot)-W_\cdot}^{-1}(B_\cdot)| (B_{j-1},B_j)\right]$ it suffices to prove that there are $c,C>0$ independent of $t$ such that $\mathbb{P}$-a.s.
\begin{align}
\mathbf{1}_{A_2}\mathbb{P}_{\mathcal{L}}\left[\mathcal{B}_{[j-1,j],f(\cdot)-W_\cdot}^{-1}(B_\cdot)| (B_{j-1},B_j)\right]&\ge c\mathbf{1}_{A_2},\qquad \forall_{j\in \{2,\dots, n\}}, \label{eq:Dtc1}\\
\mathbf{1}_{A_2}\sum_{j=2}^{n} \left(1-\mathbb{P}_{\mathcal{L}}\left[\mathcal{B}_{[j-1,j],f(\cdot)-W_\cdot}^{-1}(B_\cdot)| (B_{j-1},B_j)\right]\right) &\le C, \label{eq:Dtc3}\\
\mathbf{1}_{A_2}\mathbb{P}\left[\mathcal{B}_{[0,1], f(\cdot)-W_\cdot}^{-1} |(B_{0}, B_1)\right]&\ge c\mathbf{1}_{A_2}. \label{eq:Dtc4}
\end{align}

Equations \eqref{eq:A1toA2} and \eqref{eq:Dtc1}--\eqref{eq:Dtc4} directly imply that there is a $c>0$ such that
\begin{equation}
\frac{\mathbb{P}_{\mathcal{L}}[A_1]}{\mathbb{P}_{\mathcal{L}}[A_2]}\ge c,
\end{equation}
i.e.\@ \eqref{eq:DiskrtoCont}. Thus to finish the proof, we only need to prove that \eqref{eq:Dtc1}--\eqref{eq:Dtc4} hold.

We proceed towards this, by first proving \eqref{eq:Dtc1} and \eqref{eq:Dtc3} in one calculation. Set $x_j := W_j+1-2^{\beta}f(j)$ for $j\in\{1,\dots, n\}$. On $A_2$ and given $B_{j-1}$, $B_j$ we can couple $(B_s)_{s\in [j-1,j]}$ with a Brownian bridge $(E_s^j)_{s\in [0,1]}$ with $E_0^j = x_{j-1}$, $E_1^j = x_j$ and $B_s\le E_{s-j+1}^j$ for all $j\in \{2,\dots,n\}$, $s\in [j-1,j]$. This coupling yields that
\begin{align}
&\mathbf{1}_{A_2}\mathbb{P}_{\mathcal{L}}\left[\exists_{s\in [j-1,j]}\, -1+B_s+f(s)> W_s| (B_{j-1},B_j)\right] \notag\\
&\le \mathbf{1}_{A_2}P\left[ \exists_{s\in[0,1]} E_s^j > W_{s+j-1}-f(s+j-1)+1|W\right] \label{eq:BMzBB}
\end{align}
Let $E_s'$ be a standard Brownian bridge. We have that
\begin{equation}
(E_s^j)_{j\in [0,1]} \stackrel{d}{=} (E_s'+sx_j+(1-s)x_{j-1})_{s\in [0,1]}. \label{eq:BBzSBB}
\end{equation}
Plugging \eqref{eq:BBzSBB} into \eqref{eq:BMzBB} and simplifying, using that $W$ is piecewise linear, gives that
\begin{align}
&\mathbf{1}_{A_2}\mathbb{P}_{\mathcal{L}}\left[\exists_{s\in [j-1,j]}\, -1+B_s+f(s)> W_s| (B_{j-1},B_j)\right]\notag\\
&\le \mathbf{1}_{A_2}P\left[\exists_{s\in[j-1,j]}\, E_s'>-f(j-\mathbf{1}_{\beta =-1})+2^{\beta}f(j-\mathbf{1}_{\beta =1}) \right]. \label{eq:fbeta}
\end{align}
Direct calculation yields that there is a $c'>1/2$ such that for $j\ge 1$
\begin{align}
\frac{f(j)}{f(j+1)} \ge c'. \label{eq:fquot}
\end{align}
Plugging \eqref{eq:fquot} into  \eqref{eq:fbeta} yields that there is a $c''>0$ such that
\begin{align}
&\mathbf{1}_{A_2}\mathbb{P}_{\mathcal{L}}\left[\exists_{s\in [j-1,j]}\, -1+B_s+f(s)> W_s| (B_{j-1},B_j)\right]\notag\\
&\hskip4.5cm\le \mathbf{1}_{A_2}\mathbb{P}_{\mathcal{L}}\left[\exists_{s\in [j-1,j]}\, E_s' > c''((1+j)^{1/6}-1)\right]. \label{eq:c''}
\end{align}
Since the maximum of a standard Brownian bridge is stochastically dominated by the maximum of a standard Brownian motion, which can be proved using e.g. Slepian's lemma, \eqref{eq:c''} yields that for $j\in \{2,\dots,n\}$
\begin{equation}
\mathbf{1}_{A_2}\mathbb{P}_{\mathcal{L}}\left[\exists_{s\in [j-1,j]}-1+B_s+f(s)> W_s| (B_{j-1},B_j)\right]\le 2\mathbf{1}_{A_2}\int_{c((1+j)^{1/6}-1)}^{\infty} \frac{1}{\sqrt{2\pi}}e^{-\frac{t^2}{2}}\;\mathrm{d}t. \label{eq:DzcZw1}
\end{equation}
The equations \eqref{eq:Dtc1} and \eqref{eq:Dtc3} are directly implied by \eqref{eq:DzcZw1}.

One can prove \eqref{eq:Dtc4} similarly taking $x_0 = 0$ and showing by calculation, that $(1-s)-f(s)+s2^\beta f(1)>0$ is bounded away from zero.
\end{proof}
\begin{lemma}\label{Lem:Clogtight}
We have that $(\ClogW(W,n))_{n\in\N}$  and $(\ClogW(\widehat{W},n))_{n\in\N}$ are tight.
\end{lemma}
\begin{proof}
By construction, we have that
\[
\ClogW \le \sup_{k\in \N } \frac{|W_{k+1}|}{\sqrt{k}\sqrt{\log(k)}}<\infty\quad\mathbb{P}\text{-a.s.}
\]
 For the last inequality, we recall that by Remark \ref{Rem:LiL} we can apply the law of iterated logarithms for the random walk $(W_k)_{k\in\N}$. The second statement of the lemma follows from the first, since $\mathcal{C}_{\log}(\widehat{W},n)\stackrel{d}{=}\mathcal{C}_{\log}(W,n)$ for all $n\in\N$, since $(W_k)_{k\in [0,n]} \stackrel{d}{=} (-\widehat{W}_k)_{k\in [0,n]}$.
\end{proof}
\begin{lemma}
For all $\lambda>0$ $(\Cnull(W,h,\lambda,n))_{n\in\N}$ and $(\Cnull(\widehat{W},h,\lambda,n))_{n\in\N}$ are tight.
\end{lemma}
\begin{proof}
This is an immediate consequence of Lemmata \ref{Lem:BoundGam},  \ref{Lem:Clogtight} and the Definition of $\Cnull$ in \eqref{Ass:AssC0W}.
\end{proof}

\begin{lemma}\label{Lem:CurlyC1}
We have that $(\ConstWnW(W,h,n)^{-1})_{n\in\N}$ and $(\ConstWnW(\widehat{W},h,n)^{-1})_{t\ge0}$ are tight.
\end{lemma}
\begin{proof}
Follows from Lemma \ref{Lem:Clogtight}, \eqref{Ass:WachsH}, and the definition of $\ConstWnW$ in \eqref{Def:WeirdC1}.
\end{proof}

\paragraph{Proof of Lemma \ref{Lem:Boundpnastweak}}
 We recall Definitions \ref{Def:pthAllg}, \ref{Def:Delta} and Theorem \ref{Theo:CrudeLB}. We prove a more general statement than Lemma \ref{Lem:Boundpnastweak}, for this we need one additional definition.
\begin{definition}
We call $(h_t)_{t\in I}$, $I\in \{\N, \R_{\ge0}\}$, $h_t:[0,t]\to \R$, nice, if there is a $\eta<1/2$ such that for all $t\in I$ 
\begin{align*}
\max_{s\in [0,t]\cap I} |h_t(s)| &\le t^{\eta},\\
\max_{s\in [0,t/3]\cap I} \frac{|h_t(s)|}{(1+s)^\eta}&\le C(h),\\
\max_{s\in [0,t/3]\cap I} \frac{|h_t(t)-h_t(t-s)|}{(1+s)^\eta} &\le C(h).
\end{align*}
\end{definition}

\begin{lemma}\label{Lem:Boundpnast}
For $(h_n)_{n\in\N}$  nice and $y(n),y_0(n)\in [-\sqrt{n},0]$ (the arguments are omitted in what follows), there is a deterministic constant $C_2(h)>0$ such that $\mathbb{P}$-a.s.
\begin{equation}
\limsup_{n\to\infty} \frac{|\log(p_{n,h}^{(y_0)}(y))|}{\log(n)} \le C_2(h). \label{eq:pnastas}
\end{equation}
\end{lemma}
\begin{remark}
Lemma \ref{Lem:Boundpnastweak} is an immediate corollary of Lemma \ref{Lem:Boundpnast} by taking $h\equiv 0$ and $y(n) = y_0(n) = \xi_0$ for all $n\in\N$.
\end{remark}
\begin{proof}
Since $p_{n,h}^{(y_0)}(y)$ is a probability it is smaller than $1$, thus $\log(p_{n,h}^{(y_0)}(y)))\le 0$, and we only need to prove a lower bound on $p_{n,h}^{(y_0)}(y)$. In the following, we will use $C(h)$ for a positive constant depending on $h$, but on neither $W$ nor $n$, which may change from line to line.

Use the notation from Definition \ref{Def:Delta}. By Theorem \ref{Theo:CrudeLB} we have that  for $n\ge2$ $\mathbb{P}$-a.s.
\begin{align*}
&p_{n,h}^{(y_0)}(y) \\
&\ge n^{-\gamma_0-\log(n)^{-1}\hskip-0.07cm\left(\hskip-0.05cm\sum_{j=0}^{k_1\hskip-0.05cm-\hskip-0.05cm1}\hskip-0.07cm\Delta_j^22^j+2(1-2^{3/2})\hskip-0.05cm\sum_{j=0}^{k_1\hskip-0.05cm-\hskip-0.05cm1}\hskip-0.07cm\Delta_j2^{j/2}+\Xi_0^2+4\Xi_0+3\frac{W_n^2+h_n(n)^2}{t_2-t_1}+(\sqrt{t_1}+|y-y_0|)\frac{|W_n|+|h_n(n)|}{t_2-t_1}\hskip-0.05cm\right)},
\end{align*}
where $k_1$, $t_1$ and $t_2$ are $n$ dependent, which we omit from notation. Using that $(h_n)_{n\in\N}$ is nice, $y,y_0 \in[-\sqrt{n},0]$ and using the law of iterated logarithm for $(W_k)_{k\in\N}$, recall Remark \ref{Rem:LiL}, yields that there is a $C(h)\in (0,\infty)$ such that $\mathbb{P}$-a.s.
\begin{align*}
&\limsup_{n\to\infty} (\log(n))^{-1}\Bigg(\Xi_0^2+4\Xi_0+\frac{3W_n^2+h_n(n)^2}{2(t_2-t_1)}+(\sqrt{t_1}+|y-y_0|)\frac{|W_n|+|h_n(n)|}{t_2-t_1}\\
&\hskip4cm+\Delta_{k_1-1}^2 2^{k_1-1}+(2-2^{5/2})\Delta_{k_1-1}2^{(k_1-1)/2}\Bigg)\le C(h).
\end{align*} Thus it suffices to prove that there are $C, C(h)>0$ such that $\mathbb{P}$-a.s.
\begin{align}
\limsup_{n\to\infty} \log(n)^{-1}\sum_{j=0}^{k_1-2}\Delta_j(h)^22^j+|\Delta_j(h)| 2^{j/2} \le C(h), \label{eq:ptasth}\\
\limsup_{n\to\infty} \log(n)^{-1}\sum_{j=0}^{k_1-2}\Delta_j(W)^22^j+|\Delta_j(W)|2^{j/2} \le C, \label{eq:ptastW}
\end{align}
Since $(h_n)_{n\in\N}$ is nice we have that for all $j\le k_1-2$
\begin{align*}
|\Delta_j(h)| &\le ct^{\eta-1}+2^{\eta (j+2)-1},
\end{align*}
which gives \eqref{eq:ptasth}. It is left to prove \eqref{eq:ptastW}. For this purpose, we notice, that for all $n\in\N$, $j\le k_1-2$
\begin{align}
2^j|\Delta_j(W)| &\le \max\{ \underbrace{\max_{s\in [2^{j+1},2^{j+2}]}|W_s-W_{2^{j+1}}|}_{=:a_{j,n,1}}, \underbrace{\max_{s\in [n-2^{j+2},n-2^{j+1}]} |W_s-W_{n-2^{j+1}}|}_{:=a_{j,n,2}}\}.\label{eq:DjW1}
\end{align}
We have that $a_{j,n,1}$ does not depend on $n$ and $(a_{j,n,1})_{j\ge 1}$ is a sequence of independent random variables. With for example Doobs maximal inequality one can show that there is a $c>0$ such that for $j\le k_1-2$
\begin{align*}
\mathbb{E}[a_{j,n,1}^22^{-j}] &\le c,  \\
\mathbb{E}[|a_{j,n,1}|2^{-j/2}] &\le c. 
\end{align*}
Thus the strong law of large numbers gives that there is a $C>0$ such that $\mathbb{P}$-a.s.
\begin{align}
\limsup_{n\to\infty} \log(n)^{-1}\sum_{j=0}^{k_1-2} a_{j,n,1}^22^{-j}+|a_{j,n,1}|2^{-j/2} &\le C. \label{eq:Boundajn1}
\end{align}

For the sum involving $a_{j,n,2}$ we want to use Borel-Cantelli, and thus want to bound $\mathbb{P}[\sum_{j=0}^{k_1-2} a_{j,n,2}2^{j/2}>x]$. For this, we use the exponential Chebyshev inequality to get that for $x>0$ 
\begin{align}
\mathbb{P}\left[\sum_{j=0}^{k_1-2} |a_{j,n,2}|2^{-j/2} >x\cdot \log(n)\right] &\le e^{-x\log(n)}\prod_{j=0}^{k_1-2}\mathbb{E}\left[e^{|a_{j,n,2}|2^{-j/2} }\right].\label{eq:ExpCheb}
\end{align}
Thus we want an upper bound on the tails of $|a_{j,n,2}|$. A union bound gives that for $(\lambda_k)_{k\in\N}>0$
\begin{align*}
\mathbb{P}\left[|a_{j,n,2}|2^{-j/2} > x\right] &\le \sum_{k=1}^{2^j} \mathbb{P}\left[\left|\sum_{r=1}^{k} W_r\right|> 2^{j/2}x\right] \\
&\le \sum_{k=1}^{2^j} e^{- \lambda_k 2^{j/2}x} \left(\mathbb{E}\left[e^{\lambda_k W_1}\right]^k+\mathbb{E}\left[e^{-\lambda_k W_1}\right]^k\right).
\end{align*}
Direct calculation yields that $W_1$ is sub-Gaussian, in particular there are $C,c>0$ such that $\mathbb{E}\left[e^{\pm \lambda_k W_1}\right]\le Ce^{c\lambda_k^2}$. This implies that
\begin{align*}
\mathbb{P}\left[|a_{j,n,2}|2^{-j/2} > x\right]  \le \sum_{k=1}^{2^j} Ce^{- \lambda_k 2^{j/2}x+c\lambda_k^2k}\stackrel{\lambda_k = 2^{j/2}x(2ck)^{-1}}{\le}  \sum_{k=1}^{2^j}Ce^{-\frac{2^j x^2}{4ck}} &\le 2^jC e^{-\frac{x^2}{4c}}
\end{align*}

This implies that for $x\ge 1$
\[
\mathbb{P}\left[e^{|a_{j,n,2}|2^{-j/2}}> x \right] = \mathbb{P}\left[|a_{j,n,2}| > 2^{j/2}\log(x)\right] \le 2^jCe^{-\frac{\log(x)^22^j}{4c}}.
\]
and since $ |a_{j,n,2}|2^{-j/2}\ge 0$ we have that $\mathbb{P}\left[e^{\lambda |a_{j,n,2}|2^{-j/2}}> x \right] =1$ for $x\in [0,1]$. Thus
\begin{equation}
\mathbb{E}\left[e^{|a_{j,n,2}|2^{-j/2}}\right] \le 1+\frac{2^jc^24}{2^j-4c} \le e^{\frac{2^jc^24}{2^j-4c}}. \label{eq:Eajn2}
\end{equation}
For $j$ big enough, we have $\frac{2^jc^24}{2^j-4c} \le 8c^2$, especially by changing $c$ we get from \eqref{eq:ExpCheb} and \eqref{eq:Eajn2} that
\[
\mathbb{P}\left[\sum_{j=0}^{k_1-2} |a_{j,n,2}|2^{-j/2} >x\cdot \log(n)\right] \le e^{-x\log(n)}e^{c^2\log(n)},
\]
which is summable in $n$ for all $x$ big enough, deterministic, which implies that there is a $C>0$, such that $\mathbb{P}$-a.s.
\begin{equation}
\limsup_{n\to \infty} \log(n)^{-1}\sum_{j=0}^{k_1-2} |a_{j,n,2}|2^{-j/2} \le C. \label{eq:Resaj21}
\end{equation}
Very similarly one can prove that there is a $C>0$ such that $\mathbb{P}$-a.s.
\begin{equation}
\limsup_{n\to \infty} \log(n)^{-1}\sum_{j=0}^{k_1-2} a_{j,n,2}^22^{-j} \le C. \label{eq:Resaj22}
\end{equation}
Equations \eqref{eq:Resaj21} and \eqref{eq:Resaj22} imply that there is a $C>0$ such that $\mathbb{P}$-a.s.
\begin{equation}
\limsup_{n\to\infty} \log(n)^{-1}\sum_{j=0}^{k_1-2} a_{j,n,2}^22^{-j}+|a_{j,n,2}|2^{-j/2} \le C. \label{eq:Boundajn2}
\end{equation}
Equation \eqref{eq:DjW1}, \eqref{eq:Boundajn1} and \eqref{eq:Boundajn1} yield \eqref{eq:ptastW}, which finishes the proof of \eqref{eq:pnastas}.
\end{proof}

We come back to controlling three more constants which pop up in the bounds we derived, recall \eqref{Ass:Assdoof}, Definition \ref{Def:ptast} and that for given $t$ we set  $\widehat{W}_s := W_{t-s}-W_t$
\begin{lemma}\label{Lem:C1Wc2W}
In \eqref{Ass:Assdoof} we can choose $(C_4(W,n))_{n\in\N}$, $(c_2(W,n)^{-1})_{n\in\N}$, $(C_4(\widehat{W},n))_{n\in\N}$ and $(c_2(\widehat{W},n)^{-1})_{n\in\N}$ tight.
\end{lemma}
We omit the proof of Lemma \ref{Lem:C1Wc2W}, since its statement can be inferred from \ref{Theo:CrudeLB} very similarly to the proof of Lemma \ref{Lem:Boundpnast}.

\begin{lemma}\label{Lem:kappaEta}
For $y = y_0 = \xi_0$ and $\delta\in (0,1/8)$ fixed in Theorems \ref{Sa:GirsConvH} and \ref{Sa:GirsConcH} the constants $(\eta(W,\widehat{W},n,h,\delta,y,y_0))_{n\in\N}$ and $(\kappa(W,\widehat{W},n,h,y,y_0))_{n\in\N}$ can be chosen to be tight.
\end{lemma}
\begin{proof}
As mentioned in Remark \ref{Rem:Rem1} we can choose $\kappa$ and $\eta$  as a continuous function of $\ExpPos$, $\Cnull $, $\ConstWnW^{-1}$, $\ClogW $, $\ExpPoshat $, $\Cnullhat$, $\ConstWnWhat^{-1}$, $\ClogWhat $, $y$, $y_0$, $\ExpStart$, $\ConstStart^{-1}$, $\ExpStarthat$, $\ConstStarthat^{-1}$ which together with Lemmata \ref{Lem:BoundGam}--\ref{Lem:CurlyC1} and Lemma \ref{Lem:C1Wc2W} yields the claim.
\end{proof}
\section{Barrier Computations -- Proof of Propositions \ref{Prop:BarrComp}, \ref{Prop:barrprobLBLB} and \ref{Prop:BarrierUBLB} and Lemmata \ref{Lem:Ga}, \ref{Lem:Handlepstart}}\label{Sec:BarrComps}

\subsection{Proof of Proposition \ref{Prop:BarrComp}}\label{Subsec:PP12}
We will change the start- and endpoint, remove the $h^{\smallfrown}$ and the logarithmic drift one by one and start by defining all the probabilities, which we encounter on the way. We recall Definition \ref{Def:Barr} and $h^{\smallfrown}$, $m_{t,h}^{\smallfrown}$ from Definition \ref{Def:tau} and $p_t^{(x)}(y)$, $p_{t,h}^{(x)}(y)$ from Definition \ref{Def:pttildes}. We note, that in this section we use the $W$ defined in Definition \ref{Def:pn} and that the probability measure used is $\mathbb{P}_{\mathcal{L}}$.
\begin{definition}\label{Def:OnTheWay}
Let $t>0$. For $x,y\le 0$, define
\begin{align*}
p_{t,h^{\smallfrown},\text{barr}}^{(x)}(y) &:= \mathbb{P}_{\mathcal{L}}\left[\mathcal{B}_{[0,t], m_{t,h}^{\smallfrown}}^{y,J_x}(T_{\cdot}/\vartheta^\ast)\right],
\end{align*}
and recall \eqref{Def:barrprobUB}. Furthermore, define $\widehat{W}_s := W_{t-s}-W_t$, $s\in [0,t]$, and $\widehat{T}_s := \vartheta^\ast(B_s-\widehat{W}_s)$ and
\[
\widehat{p}_{t,h^{\smallfrown}}^{(x)}(y):= \mathbb{P}_{\mathcal{L}}\left[\mathcal{B}_{[0,t], h_t^{\smallfrown}}^{y,J_x}(\widehat{T}_\cdot/\vartheta^\ast)\right].
\]
\end{definition}
We will now state the lemmata necessary to prove Proposition \ref{Prop:BarrComp}, then show the proof and after that prove the lemmata. We recall the definition of $\barrierprobUB$ in \eqref{Def:barrprobUB}.
\begin{lemma} \label{Lem:RemDiskr}
There is a $C>0$ such that for $n\in\N$ and all $x,y\le \xi_0$, $\mathbb{P}$-a.s.
\[
\frac{\barrierprobUB}{p_{n,h^{\smallfrown},\text{barr}}^{(x)}(y)} \le C.
\]
\end{lemma}
The proof of Lemma \ref{Lem:RemDiskr} is similar to the proof of \eqref{eq:Dtc1}--\eqref{eq:Dtc3} in Lemma \ref{Lem:BoundGam} for $\beta = -1$, so we won't give details. In words, Lemma \ref{Lem:RemDiskr} means, that the events to be below the barrier for all $j\in \{1,\dots, n\}$ and to be below the barrier (times two) for all $s\le n$ do only differ by a constant.
\begin{lemma} \label{Lem:RemBarr}
For all $\varepsilon>0$, there is a $C_\varepsilon>0$ such that
\[
\liminf_{n\to\infty} \mathbb{P}\left[\bigcap_{y\in [-\log(n)^2,\xi_0] \cap \Z \atop x\in [-\log(n),\xi_0] \cap \Z}\frac{p_{n,h^{\smallfrown},\text{barr}}^{(x)}(y)}{p_{n,h^{\smallfrown}}^{(x)}(y)}\le C_\varepsilon\right]\ge 1-\varepsilon.
\]
\end{lemma}
In words, Lemma \ref{Lem:RemBarr} means, that removing the $\frac{j}{n}\log(p_n)$ part of the barrier costs a multiplicative constant.

\begin{lemma}\label{Lem:Remy}
For all $\varepsilon>0$, there are $c_\varepsilon, C_\varepsilon>0$ such that 
\begin{align}
\liminf_{n\to\infty} \mathbb{P}\left[\bigcap_{y\in [-\log(n)^2,\xi_0]\cap\Z\atop x\in [-\log(n),\xi_0]\cap\Z}\left\{\frac{p_{n,h^{\smallfrown}}^{(x)}(y)}{p_{n,h^{\smallfrown}}^{(x)}(\xi_0+1)} \le C_\varepsilon|y|^{c_\varepsilon}\right\}\right]&\ge 1-\varepsilon, \label{eq:P17Movestart}\\
\liminf_{n\to \infty} \mathbb{P}\left[\bigcap_{y\in [-\log(n)^2,\xi_0+1]\cap\Z\atop x\in [-\log(n),\xi_0+1]\cap\Z}\left\{\frac{\widehat{p}_{n,h^{\smallfrown}}^{(x)}(y)}{\widehat{p}_{n,h^{\smallfrown}}^{(x)}(\xi_0)} \le C_\varepsilon|y|^{c_\varepsilon}\right\}\right]&\ge 1-\varepsilon. \label{eq:P17Moveend}
\end{align}
\end{lemma}
Lemma \ref{Lem:Remy} states that moving the startpoint or endpoint from $y$ to $\xi_0$ costs a power of $|y|$.

\begin{lemma}\label{Lem:Switch}
For all $n\in\N$, $x,y \le 0$, we have that $\mathbb{P}$-a.s.
\begin{align*}
p_{n,h^{\smallfrown}}^{(x)}(y) &\le \widehat{p}_{n,h^{\smallfrown}}^{(y)}(x-1),\\
\widehat{p}_{n,h^{\smallfrown}}^{(x)}(y)&\le  p_{n,h^{\smallfrown}}^{(y)}(x-1).
\end{align*}
\end{lemma}
Lemma \ref{Lem:Switch} states that reversing time makes (close to) no difference.

\begin{lemma}\label{Lem:Remh}
For all $\varepsilon>0$, there is a $C_\varepsilon>0$ such that
\[
\liminf_{n\to \infty} \mathbb{P}\left[\frac{p_{n,h^{\smallfrown}}^{(\xi_0)}(\xi_0)}{p_n} \le C_\varepsilon\right]\ge 1-\varepsilon.
\]
\end{lemma}

Lemma \ref{Lem:Remh} states that when starting and ending at $\xi_0$, the curve of the barrier can be removed by paying a multiplicative constant.

Now we have stated everything we need to prove Proposition \ref{Prop:BarrComp}.
\begin{proof}[Proof of Propositions \ref{Prop:BarrComp} assuming L. \ref{Lem:RemDiskr}--\ref{Lem:Remh}]
In this proof we shorten $h^{\smallfrown}$ to $h$.
We can write
\begin{align*}
&\frac{\barrierprobUB}{p_n} \\
&= \frac{\barrierprobUB}{p^{(x)}_{n,h,\text{barr}}(y)}\frac{p^{(x)}_{n,h,\text{barr}}(y)}{p^{(x)}_{n,h}(y)}\frac{p^{(x)}_{n,h}(y)}{p^{(x)}_{n,h}(\xi_0+1)}\frac{p^{(x)}_{n,h}(\xi_0+1)}{\widehat{p}^{(\xi_0+1)}_{n,h}(x-1)}\frac{\widehat{p}^{(\xi_0+1)}_{n,h}(x-1)}{\widehat{p}^{(\xi_0+1)}_{n,h}(\xi_0)}\frac{\widehat{p}^{(\xi_0+1)}_{n,h}(\xi_0)}{p^{(\xi_0)}_{n,h}(\xi_0)}\frac{p^{(\xi_0)}_{n,h}(\xi_0)}{p_n}.
\end{align*}
Applying Lemmata \ref{Lem:RemDiskr} to \ref{Lem:Remh} then immediately yields the claim.
\end{proof}

Now we prove Lemmata \ref{Lem:RemBarr}--\ref{Lem:Remh}. We note that these proofs will make heavy use of Sections \ref{Sec:MoveStart} and \ref{Sec:Moveh}.

\begin{proof}[Proof of Lemma \ref{Lem:RemBarr}]
Applying Theorem \ref{Theo: LinTerm} yields that for all $n\in\N$, $y\in [-\log(n)^{2},0]$, $x\in [-\log(n),0]$, $\mathbb{P}$-a.s.
\[
\frac{p_{n,h^{\smallfrown},\text{barr}}^{(x)}(y)}{p_{n,h^{\smallfrown}}^{(x)}(y)} \le e^{-\frac{\log(p_n)^2}{2n(\vartheta^\ast)^2}}e^{-3\log(n)^{2}\frac{\log(p_n)}{\vartheta^\ast n}}e^{ W_n \frac{\log(p_n)}{\vartheta^\ast n}}.
\]
Since by Lemma \ref{Lem:Boundpnast} $\mathbb{P}$-a.s.\@ $\log(p_n)/\log(n)^2\to 0$ and furthermore $\log(n)W_n/n\to 0$ $\mathbb{P}$-a.s., the statement of the lemma follows.
\end{proof}

\begin{proof}[Proof of Lemma \ref{Lem:Remy}]
Applying Theorems \ref{Theo:NebOS} and \ref{Theo: MainOS} yields that for all $\lambda>0$
\begin{align*}
1 &= \liminf_{n\to \infty}\mathbb{P}\Bigg[\bigcap_{y\in [-\log(n)^2,\xi_0]\atop x\in [-\log(n),\xi_0]} \Big\{\frac{p_{n,h^{\smallfrown}}^{(x)}(y)}{p_{n,h^{\smallfrown}}^{(x)}(\xi_0+1)}\le \ConstWnW^{-2}4^{\ExpPos+2}e^{4\Coh }|y|^{4\ExpPos+3\Cnull }\\
&\hskip6cm+Ce^{-\lambda \log(n)}p_{n,h^{\smallfrown}}^{(x)}(\xi_0+1)^{-1}(2\log(n)^{2})^{2\ExpPos}\ConstWnW^{-1}\Big\}\Bigg].
\end{align*}

By Lemma \ref{Lem:BoundGam} to \ref{Lem:CurlyC1} the sequences $(\ExpPos(W,h^{\smallfrown},n))_{n\in\N}$, $(\Cnull(W,h^{\smallfrown},\lambda,n))_{n\in\N} $ and $(\ConstWnW(W,h^{\smallfrown},n)^{-1})_{n\in\N}$ are tight, and Lemma \ref{Lem:Boundpnast} implies that there is a $\lambda>0$ such that $\mathbb{P}$-a.s.\@
\[
\lim\limits_{n\to\infty} \sup_{x\in [-\log(n),\xi_0]}Ce^{-\lambda \log(n)}p_{n}^{(x)}(\xi_0+1)^{-1}(2\log(n)^{2})^{2\ExpPos}\ConstWnW^{-1} = 0.
\]  This allows us to conclude \eqref{eq:P17Movestart}.

Equation \eqref{eq:P17Moveend} follows the same way, only replacing $W$ by $\widehat{W}$ in the above and using that by Lemma \ref{Lem:BoundGam} to \ref{Lem:CurlyC1} $(\gamma(\widehat{W},h^{\smallfrown},n))_{n\in\N}$,  $(C_3(\widehat{W},h^{\smallfrown},\lambda,n))_{n\in\N}$, $(\ConstWnW(\widehat{W},h^{\smallfrown},n)^{-1})_{n\in\N}$ are tight.
\end{proof}

\begin{proof}[Proof of Lemma \ref{Lem:Switch}]
We have that
\begin{align}
&\mathbb{P}_{\mathcal{L}}\left[\mathcal{B}_{\{0,\dots,n\}, h^{\smallfrown}_n(\cdot)-W_\cdot}^{y,J_x}(B_\cdot)\right]\notag\\
&= \int_{x-1}^{x} \mathbb{P}_{\mathcal{L}}\left[B_n = \mathrm{d}u+ W_n-y \right]\mathbb{P}_{\mathcal{L}}\left[\mathcal{B}_{\{0,\dots,n\},h^{\smallfrown}_n(\cdot)-W_\cdot}^{y}(B_\cdot)|B_n = u+W_n-y\right]\label{eq:Switch1}
\end{align}
Set $\widehat{W}_j := W_{n-j}-W_n$, $\widehat{B}_j := B_{n-j}-B_n$. Since $h^{\smallfrown}_n$ is symmetric, reversing time in \eqref{eq:Switch1} gives
\begin{align}
&\mathbb{P}_{\mathcal{L}}\left[\mathcal{B}_{\{0,\dots,n\}, h^{\smallfrown}_n(\cdot)-W_\cdot}^{y,J_x}(B_\cdot)\right]\notag\\
&=\int_{x-1}^{x} \mathbb{P}_{\mathcal{L}}\left[B_n = \mathrm{d}u+ W_n-y \right]\mathbb{P}_{\mathcal{L}}\left[\mathcal{B}_{\{0,\dots,n\},h^{\smallfrown}_n(\cdot)-\widehat{W}_\cdot}^{u}(\widehat{B}_\cdot)|B_n = u+W_n-y\right]\notag\\
&=\int_{-x}^{-x+1} \mathbb{P}_{\mathcal{L}}\left[\widehat{B}_n = \mathrm{d}u+\widehat{W}_n+y \right]\mathbb{P}_{\mathcal{L}}\left[\mathcal{B}_{\{0,\dots,n\},h^{\smallfrown}_n(\cdot)-\widehat{W}_\cdot}^{u}(\widehat{B}_\cdot)|\widehat{B}_n = u+\widehat{W}_n+y\right],\label{eq:Switch2}
\end{align}
where the last step used that $\widehat{B}_n =-B_n$ and $\widehat{W}_n =-W_n$.

Since $\mathbb{P}_{\mathcal{L}}\left[\mathcal{B}_{\{0,\dots,n\},h^{\smallfrown}_n(\cdot)-\widehat{W}_\cdot}^{\mathbf{u}}(\widehat{B}_\cdot)|\widehat{B}_n = u-W_n+y\right]$ is monotone in the bolded $u$ we can replace it by $-x+1$. Furthermore, we can shift the region of integration by $y-1+x$ in \eqref{eq:Switch2}, which yields that 
\begin{align}
&\mathbb{P}_{\mathcal{L}}\left[\mathcal{B}_{\{0,\dots,n\}, h^{\smallfrown}_n(\cdot)-W_\cdot}^{y,J_x}(B_\cdot)\right]\notag\\
&\le \int_{y-1}^{y}\mathbb{P}_{\mathcal{L}}\left[\widehat{B}_n = \mathrm{d}r+\widehat{W}-x+1\right]\mathbb{P}_{\mathcal{L}}\left[\mathcal{B}_{\{0,\dots,n\}, h^{\smallfrown}_n(\cdot)-\widehat{W}_\cdot}^{x-1}(\widehat{B}_\cdot)|\widehat{B}_n = r+\widehat{W}-x+1\right]\notag\\
&= \mathbb{P}_{\mathcal{L}}\left[\mathcal{B}_{\{0,\dots,n\},h^{\smallfrown}_n(\cdot)-\widehat{W}_\cdot}^{x-1,J_y}(\widehat{B}_\cdot)\right]. \label{eq:Switch3}
\end{align}
Since $(\widehat{B}_j)_{j\in\{0,\dots,n\}}\stackrel{d}{=} (B_j)_{j\in\{0,\dots,n\}}$ equation \eqref{eq:Switch3} implies the claim of the lemma.
\end{proof}

\begin{proof}[Proof of Lemma \ref{Lem:Remh}]
Denote $A_n := \left\{\log(n)^{-0.5}\max\{\ClogW ,\ClogWhat \}\le \frac{1}{3}\right\}$. By Lemma \ref{Lem:Clogtight} we have that $\log(n)^{-0.5}\max\{\ClogW ,\ClogWhat \}\to 0$ in $\mathbb{P}$-probability and thus $\mathbf{1}_{A_n}\to 1$ in $\mathbb{P}$-probability. 
On $A_n$ condition \eqref{eq:Limgut} holds, and thus we can apply Theorem \ref{Sa:GirsConvH} to get that there are $c,C, n_0(\lambda)>0$ such that for $n\ge n_0(\lambda)$
\begin{align}
&\mathbf{1}_{A_n}\frac{p_{n,h^{\smallfrown}}^{(\xi_0)}(\xi_0)-\frac{\delta}{1-\delta}\sum_{j=0}^2p_{n,h^{\smallfrown}}^{(\xi_0-1)}(\xi_0-j)}{p_n} \notag\\
&\le \mathbf{1}_{A_n}\Bigg(C e^{c(\ClogW+\kappa-\eta) }\left(1+n^{-\lambda}p_n^{(y_0)}(y)^{-1}\right)+C n^{-\lambda}((\pnast)^{-1}+(\pnasthat)^{-1})p_n^{(y_0)}(y)^{-1}\Bigg).\label{eq:L115}
\end{align}
By Lemma \ref{Lem:Remy} we know that there is a $C_{\varepsilon,1}>0$ such that
\begin{equation}
1-\varepsilon/2\le \liminf_{n\to\infty} \mathbb{P}\left[\frac{\sum_{j=0}^2p_{n,h^{\smallfrown}}^{(\xi_0-1)}(\xi_0-j)}{p_{n,h^{\smallfrown}}^{(\xi_0)}(\xi_0)} \le C_{\varepsilon,1}\right]. \label{eq:Approxpxixi}
\end{equation}
Choosing $\delta_\varepsilon := \frac{1}{2C_{\varepsilon,1}+1}\wedge \frac{1}{8}$ in \eqref{eq:L115}, applying \eqref{eq:Approxpxixi} and using $\mathbf{1}_{A_n}\to 1$ in $\mathbb{P}$-probability implies that for all $\lambda>0$
\begin{align}
1-\varepsilon/2 &\le \liminf_{n\to\infty}\mathbb{P}\Bigg[\frac{p_{n,h^{\smallfrown}}^{(\xi_0)}(\xi_0)}{p_n}\le \Bigg(C e^{c(\ClogW+\kappa(W,\widehat{W},n,h^{\smallfrown},\xi_0,\xi_0)-\eta(W,\widehat{W},n,h^{\smallfrown},\delta_\varepsilon,\xi_0,\xi_0)) }\cdot\notag\\
&\hskip2.9cm\cdot\left(1+n^{-\lambda}p_n^{(y_0)}(y)^{-1}\right)+C n^{-\lambda}((\pnast)^{-1}+(\pnasthat)^{-1})p_n^{(y_0)}(y)^{-1}\Bigg].\label{eq:L115'}
\end{align}

By Lemma \ref{Lem:Clogtight} we have that $(\ClogW(W,n))_{n\in\N}$ is tight. By Lemma \ref{Lem:kappaEta} the sequences  $(\eta(W,\widehat{W},n,h^{\smallfrown},\delta_\varepsilon,\xi_0,\xi_0))_{n\in\N}$ and $(\kappa(W,\widehat{W},n,h^{\smallfrown},\xi_0,\xi_0))_{n\in\N}$ are tight. By Lemma \ref{Lem:Boundpnast} and Definition \ref{Def:ptast} there is a $\lambda>0$ such that the terms $Cn^{-\lambda}((\pnast)^{-1}+(\pnasthat)^{-1})p_n^{(y_0)}(y)^{-1}$ and $n^{-\lambda}p_n^{-1}$ go to zero in $\mathbb{P}$-probability. Thus we can conclude the statement of Lemma \ref{Lem:Remh} from \eqref{eq:L115'}.
\end{proof}

\subsection{Proof of Proposition \ref{Prop:barrprobLBLB}.}
This section parallels Section \ref{Subsec:PP12}, although it is slightly simpler since we do not need to consider endpoints other than the (negative) constant $\xi_0$.  Thus recall $h^{\smallsmile}$ and $m_{t,h}^{\smallsmile}$ from Definition \ref{Def:BasicsLB} as well as Definition \ref{Def:OnTheWay}.
The proof of Proposition \ref{Prop:barrprobLBLB} is split into several lemmata, we will next state those and then prove the Proposition assuming the lemmata. Since the proofs of the lemmata are very similar to Section \ref{Subsec:PP12}, we won't repeat them, but instead just reference the corresponding lemmata in Section \ref{Subsec:PP12}.

\begin{lemma} \label{Lem:BarrRemLB}
For all $\varepsilon>0$, there is a $C_\varepsilon>0$ such that for $y\le \xi_0$,
\[
\liminf_{n\to \infty}\mathbb{P}\left[\frac{p_{n,h^{\smallsmile},\text{barr}}^{(\xi_0)}(y)}{p_{n,h^{\smallsmile}}^{(\xi_0)}(y)}\ge C_\varepsilon\right]\ge 1-\varepsilon.
\]
\end{lemma}
The proof of Lemma \ref{Lem:BarrRemLB} parallels the proof of Lemma \ref{Lem:RemBarr} and thus is based on Theorem \ref{Theo: LinTerm}.
\begin{lemma}\label{Lem:yRemLB}
For all $\varepsilon>0$, there are $C_\varepsilon, c_\varepsilon>0$ such that for $y\le \xi_0$,
\[
\liminf_{n\to \infty}\mathbb{P}\left[\frac{p_{n,h^{\smallsmile}}^{(\xi_0)}(y)}{p_{n,h^{\smallsmile}}^{(\xi_0)}(\xi_0)} \ge C_\varepsilon |y|^{-c_\varepsilon}\right]\ge 1-\varepsilon.
\]
\end{lemma}
The proof of Lemma \ref{Lem:yRemLB} is parallel to the proof of Lemma \ref{Lem:Remy}, but uses Theorem \ref{Theo:MainThLB} instead of Theorem \ref{Theo: MainOS}.
\begin{lemma}\label{Lem:hRemLB}
For all $\varepsilon>0$, there is a $C_\varepsilon>0$ such that 
\[
\liminf_{n\to \infty}\mathbb{P}\left[\frac{p_{n,h^{\smallsmile}}^{(\xi_0)}(\xi_0)}{p_n}\ge C_\varepsilon\right]\ge 1-\varepsilon.
\]
\end{lemma}
The proof of Lemma \ref{Lem:hRemLB} is parallel to the proof of Lemma \ref{Lem:Remh} using Theorem \ref{Sa:GirsConcH} in place of Theorem \ref{Sa:GirsConvH}.

\begin{proof}[Proof of Proposition \ref{Prop:barrprobLBLB} assuming Lemma \ref{Lem:BarrRemLB}--\ref{Lem:hRemLB}]
We have that
\begin{align*}
\frac{\barrierprobLB}{p_n} &= \frac{\barrierprobLB}{p_{n,h^{\smallsmile},\text{barr}}^{(\xi_0)}(y)}\frac{p_{n,h^{\smallsmile},\text{barr}}^{(\xi_0)}(y)}{p_{n,h^{\smallsmile}}^{(\xi_0)}(y)}\frac{p_{n,h^{\smallsmile}}^{(\xi_0)}(y)}{p_{n,h^{\smallsmile}}^{(\xi_0)}(\xi_0)}\frac{p_{n}^{(\xi_0)}(\xi_0)}{p_n}\\
&\ge \frac{p_{n,h^{\smallsmile},\text{barr}}^{(\xi_0)}(y)}{p_{n,h^{\smallsmile}}^{(\xi_0)}(y)}\frac{p_{n,h^{\smallsmile}}^{(\xi_0)}(y)}{p_{n,h^{\smallsmile}}^{(\xi_0)}(\xi_0)}\frac{p_{n}^{(\xi_0)}(\xi_0)}{p_n},
\end{align*}
where we used that by monotonicity $\barrierprobLB/p_{n,h^{\smallsmile},\text{barr}}^{(\xi_0)}(y)\ge1$. Now, applying Lemmata \ref{Lem:BarrRemLB} to \ref{Lem:hRemLB} immediately yields the claim of the Proposition.
\end{proof}

\subsection{Proof of Proposition \ref{Prop:BarrierUBLB}}\label{Subsec:BarrierUBLB}
We recall $h^{\smallsmile}$ and $m_{t,h}^{\smallsmile}$ from Definition \ref{Def:BasicsLB} as well as Definition \ref{Def:OnTheWay}. We won't give proofs for the lemmata used in the proof of Proposition \ref{Prop:BarrierUBLB}, but instead reference lemmata, which are proved similarly.

\begin{lemma}\label{Lem:remlin}
For all $\varepsilon>0$, there is a $C_\varepsilon>0$ such that
\[
\liminf_{n\to \infty}\mathbb{P}\left[\bigcap_{y\in [-2\log(n)^{2},\xi_0]\cap \Z}\left\{\frac{\barrierprobLB}{p_{n,h^{\smallsmile}/2,\text{barr}}^{(\xi_0)}(y)}\le C_\varepsilon\right\}\right]\ge 1-\varepsilon.
\]
\end{lemma}
The proof is similar to the proof of Lemma \ref{Lem:RemDiskr}.

\begin{lemma}\label{Lem:rembarr}
For all $\varepsilon>0$, there is a $C_\varepsilon>0$ such that 
\[
\liminf_{n\to \infty}\mathbb{P}\left[\bigcap_{y\in [-2\log(n)^2,\xi_0]\cap\Z}\left\{\frac{p_{n,h^{\smallsmile}/2,\text{barr}}^{(\xi_0)}(y)}{p_{n,h^{\smallsmile}/2}^{(\xi_0)}(y)}\le C_\varepsilon\right\}\right]\ge 1-\varepsilon.
\]
\end{lemma}

The proof is parallel to the proof of Lemma \ref{Lem:RemBarr}.

\begin{lemma}\label{Lem:remy}
For all $\varepsilon>0$, there are $C_\varepsilon, c_\varepsilon>0$ such that
\[
\liminf_{n\to \infty}\mathbb{P}\left[\bigcap_{y\in[-2\log(n)^2,\xi_0]\cap\Z}\left\{\frac{p_{n,h^{\smallsmile}/2}^{(\xi_0)}(y)}{p_{n,h^{\smallsmile}/2}^{(\xi_0)}(\xi_0)}\le C_\varepsilon |y|^{c_\varepsilon}\right\}\right]\ge 1-\varepsilon.
\]
\end{lemma}
The proof is analogous to the proof of Lemma \ref{Lem:Remy}.

\begin{proof}[Proof of Proposition \ref{Prop:BarrierUBLB} assuming Lemmata \ref{Lem:remlin}, \ref{Lem:rembarr} and \ref{Lem:remy}]
By monotonicity $p_{n,h/2}^{(\xi_0)}(\xi_0)/p_n\le 1$. Thus writing
\[
\frac{\barrierprobLB}{p_n} = \frac{\barrierprobLB}{p_{n,h^{\smallsmile}/2,\text{barr}}^{(\xi_0)}(y)}\frac{p_{n,h^{\smallsmile}/2,\text{barr}}^{(\xi_0)}(y)}{p_{n,h^{\smallsmile}/2}^{(\xi_0)}(y)}\frac{p_{n,h^{\smallsmile}/2}^{(\xi_0)}(y)}{p_{n,h^{\smallsmile}/2}^{(\xi_0)}(\xi_0)}\frac{p_{n,h^{\smallsmile}/2}^{(\xi_0)}(\xi_0)}{p_n}
\]
and applying Lemmata \ref{Lem:remlin}, \ref{Lem:rembarr} and \ref{Lem:remy} yields the claim of Proposition \ref{Prop:BarrierUBLB}. 
\end{proof}
\subsection{Proof of Lemma \ref{Lem:Ga}}

We recall Definitions \ref{Def:Barr}, \eqref{Def:pt} and \eqref{Def:qkend}. Before we can proceed with the proof, we need one additional definition.
\begin{definition}
Define
\begin{align*}
W^{(k)}_j &:= W_{n-k-1+j}-W_{n-k-1},\\
\mathcal{C}_{\text{lin},k}(W) &:= e^{-(k+1)\frac{\left(\frac{1}{k+1}((2+k)^{1/6}-1)+\log(p_n)/(\vartheta^\ast n)\right)^2}{2}} e^{-2|\xi_0-3|\cdot \left|  \frac{(2+k)^{1/6}-1}{k+1}+\frac{\log(p_n)}{n\vartheta^\ast}\right|}\\
&\qquad e^{\left((2+k)^{1/6}-1\right)\left(\frac{(2+k)^{1/6}-1}{k+1}+\frac{\log(p_n)}{n\vartheta^\ast}\right)} e^{W_{k+1}^{(k)}\left(-\frac{(2+k)^{1/6}-1}{k+1}-\frac{\log(p_n)}{n\vartheta^\ast}\right)}.
\end{align*}
\end{definition}

\begin{proof}[Proof of Lemma \ref{Lem:Ga}]
Shorten $((2+k)^{1/6}-(2+k-s)^{1/6}) =: \eta_k(s)$.

By applying Theorem \ref{Theo: LinTerm} for $t = k+1$, $c_{k+1} = -\frac{(2+k)^{1/6}-1}{k+1}-\frac{\log(p_n)}{n\vartheta^\ast}$ we get for all $k\le \log(n)^7$
\begin{align*}
\inf_{x\in [\xi_0-2,-1]}q_{k,\text{end}}(x)&\ge \mathcal{C}_{\text{lin},k}(W)\cdot\inf_{x\in [\xi_0-2,-1]}\mathbb{P}_{\mathcal{L}}^{n-k-1}\left[\mathcal{B}_{[0,k+1],\eta_k(s)-s/(k+1)\eta_k(k+1)}^{x,J_{\xi_0}}(T_{\cdot}/\vartheta^\ast)\right].
\end{align*}
Since $\eta(s)$ is convex and $\eta_k(0) = 0$ we have that $\eta_k(s)-\frac{s}{k+1}\eta_k(k+1)\le 0$ for all $s\in[0,k+1]$ and can just drop that term only making the probability smaller.  Thus we have that for all $k\le \log(n)^7$
\begin{align}
\inf_{x\in [\xi_0-2,-1]}q_{k,\text{end}}(x) &\ge \mathcal{C}_{\text{lin},k}(W) \inf_{x\in [\xi_0-2,-1]}\mathbb{P}_{\mathcal{L}}^{n-k-1}\left[\mathcal{B}_{[0,k+1]}^{x,J_{\xi_0}}(T_{\cdot}/\vartheta^\ast)\right] \label{eq:qkendLB}
\end{align}
We continue by providing lower bounds for both factors in \eqref{eq:qkendLB}. Applying the law of iterated logarithms for $(W_k)_{k\in\N}$, which has been justified in Remark \ref{Rem:LiL}, allows us to conclude that there is a $c_\varepsilon>0$ such that
\begin{equation}
\liminf_{n\to \infty}\mathbb{P}\left[\bigcap_{k\le \lfloor\log(n)^7\rfloor}\left\{\frac{|W_{k+1}^{(k)}|}{(k+1)^{3/4}}\le c_\varepsilon \right\}\right]\ge 1-\varepsilon/4. \label{eq:WKClin}
\end{equation}
 Furthermore, by Lemma \ref{Lem:Boundpnast} we have that $\mathbb{P}$-a.s.
\begin{equation}
\frac{|\log(p_n)|}{\log(n)^2} \to 0. \label{eq:pnClin}
\end{equation} Plugging \eqref{eq:WKClin} and \eqref{eq:pnClin} into the definition of $\mathcal{C}_{\text{lin},k}(W)$ and using that we only care about $k\le \log(n)^7$ yields that there is a $C_\varepsilon>0$ such that 
\begin{equation}
\liminf_{n\to \infty}\mathbb{P}\left[\bigcap_{k\le \lfloor\log(n)^7\rfloor}\left\{\mathcal{C}_{\text{lin},k}(W) \ge C_\varepsilon\right\}\right]\ge 1-\varepsilon/2. \label{eq:ClinBound}
\end{equation}

By Lemma \ref{Lem:Boundpnast} for all $n\in\N$ there are random variables $\mathcal{C}_1^{(n)}$, $\mathcal{C}_2^{(n)}$ such that for all $k\le n$
\begin{align*}
\inf_{x\in [\xi_0-2,-1]}\mathbb{P}_{\mathcal{L}}^{n-k}\left[\mathcal{B}_{[0,k]}^{x,J_{\xi_0}}(T_{\cdot}/\vartheta^\ast)\right] \ge \mathcal{C}_1^{(n)}k^{-\mathcal{C}_2^{(n)}}
\end{align*}
and since the environment is i.i.d.\@ we can choose them such that for $i\in\{1,2\}$ the distribution of $\mathcal{C}_i^{(n)}$ does not depend on $n$. This implies that there are $C_\varepsilon, c_\varepsilon>0$ such that
\begin{equation}
\liminf_{n\to\infty} \mathbb{P}\left[\bigcap_{k\le \lfloor\log(n)^7\rfloor} \left\{\inf_{x\in [\xi_0-2,-1]}\mathbb{P}_{\mathcal{L}}^{n-k-1}\left[\mathcal{B}_{[0,k+1]}^{x,J_{\xi_0}}(T_{\cdot}/\vartheta^\ast)\right] \ge C_\varepsilon (k+1)^{-c_\varepsilon}\right\}\right]\ge 1-\varepsilon/2. \label{eq:pkBound}
\end{equation}

Plugging equations \eqref{eq:ClinBound} and \eqref{eq:pkBound} into \eqref{eq:qkendLB} and the observation that 
\[
\sup_n \sum_{k=1}^{\lfloor \log(n)^7\rfloor} e^{\vartheta^\ast h_n^{\smallfrown}(k)/2}C_\varepsilon (k+1)^{c_\varepsilon} <\infty
\]
finishes the proof.
\end{proof}
\subsection{Proof of Lemma \ref{Lem:Handlepstart}}
We recall Definition \ref{Def:pstark}.

Before we can proceed with the proof of Lemma \ref{Lem:Handlepstart} we need one additional definition.
\begin{definition}
Set
\begin{align*}
\mathcal{C}_{\text{lin},k}'(W) &:= e^{-\frac{k\log(p_n)^2}{2(n\vartheta^\ast)^2}}e^{-h_n^{\smallsmile}(k)\frac{\log(p_n)}{n\vartheta^\ast}}e^{-2(\xi_0-2)\frac{\log(p_n)}{n\vartheta^\ast}}e^{-\frac{\log(p_n)W_{k}}{n\vartheta^\ast}}.
\end{align*}
\end{definition}

\begin{proof}[Proof of Lemma \ref{Lem:Handlepstart}]
The proof is parallel to the proof of Corollary \ref{Lem:Ga}. Applying Theorem \ref{Theo: LinTerm} for $t = k$, $c_{k} = -\frac{\log(p_n)}{n\vartheta^\ast}$ yields that $\mathbb{P}$-a.s.\@ for all $k\le \log(n)^7$
\begin{align}
p_{\text{start},k} &\ge \mathcal{C}_{\text{lin},k}'(W)\inf_{x\in[\xi_0-1,0]} \mathbb{P}_{\mathcal{L}}\left[\mathcal{B}_{\{0,\dots,k\},h_n^{\smallsmile}}^{\xi_0,J_x}(T_{\cdot}/\vartheta^\ast)\right]\label{eq:Seedpnstart} 
\end{align}
As in the proof of Corollary \ref{Lem:Ga} combining Lemma \ref{Lem:Boundpnast} with the fact that there is a $c_\varepsilon>0$ such that
\begin{equation*}
\liminf_{n\to \infty}\mathbb{P}\left[\bigcap_{k\le \lfloor\log(n)^7\rfloor}\left\{\frac{|W_{k+1}|}{(k+1)^{3/4}}\le c_\varepsilon \right\}\right]\ge 1-\varepsilon/4,
\end{equation*} yields that there is a $C_\varepsilon>0$ such that 
\begin{equation}
\liminf_{n\to \infty}\mathbb{P}\left[\bigcap_{k\le \lfloor\log(n)^7\rfloor}\left\{\mathcal{C}_{\text{lin},k}'(W) \ge C_\varepsilon\right\}\right]\ge 1-\varepsilon/2. \label{eq:Clin'Bound}
\end{equation}

Furthermore, by setting $g_k : [0,k]\to \R, s\mapsto (1+s)^{1/6}-1$ we have that
\[
\mathbb{P}_{\mathcal{L}}\left[\mathcal{B}_{\{0,\dots,k\},h_n^{\smallsmile}}^{\xi_0,J_x}(T_{\cdot}/\vartheta^\ast)\right] = p_{k,g_k}^{(x)}(\xi_0)
\] and applying Lemma \ref{Lem:Boundpnast} yields that there are $c_\varepsilon>0$, $C_\varepsilon>0$ such that 
\begin{equation}
\liminf_{n\to\infty} \mathbb{P}\left[\bigcap_{k\le \lfloor\log(n)^7\rfloor} \left\{\inf_{x\in [\xi_0-2,-1]}\mathbb{P}_{\mathcal{L}}\left[\mathcal{B}_{\{0,\dots,k\},h_n^{\smallsmile}}^{\xi_0,J_x}(T_{\cdot}/\vartheta^\ast)\right] \ge C_\varepsilon k^{-c_\varepsilon}\right\}\right]\ge 1-\varepsilon/2. \label{eq:pkBound'}
\end{equation}
Plugging \eqref{eq:Clin'Bound} and \eqref{eq:pkBound'} into \eqref{eq:Seedpnstart} together with the observation that
\[
\sup_n\sum_{k=0}^{\lfloor\log(n)^7\rfloor} e^{-\vartheta^\ast h_n(k)^{\smallsmile}} C_\varepsilon k^{c_\varepsilon} <\infty
\]
finishes the proof.
\end{proof}
\printbibliography
\end{document}